\documentclass[AMA,STIX1COL]{WileyNJD-v2}

\articletype{Article Type}%

\usepackage{xcolor}
\usepackage{graphicx}
\usepackage{tikz}
\usepackage{filecontents}
\usepackage{pgfplots}

\usetikzlibrary{external}
\usepackage{amsmath}
\usepackage{amssymb}
\usepackage{epsfig}
\usepackage{hyperref}
\usepackage{stackengine}
\newcommand{\eps}{\varepsilon}
\newcommand{\beq}{\begin{equation}}
\newcommand{\eeq}{\end{equation}}
\newcommand{\bea}{\begin{eqnarray}}
\newcommand{\eea}{\end{eqnarray}}

\newcommand{\half}{\mbox{$\frac{1}{2}$}}
\newcommand{\dif}{{\textrm d}}
\newcommand{\dx}{{\rm d}x}
\received{26 April 2016}
\revised{6 June 2016}
\accepted{6 June 2016}

\begin{document}

\title{Efficient  formulation of a geometrically nonlinear beam element}

\author[1]{Milan Jir\'{a}sek}

\author[1,2]{Emma La Malfa Ribolla*}

\author[1]{Martin  Hor\'{a}k}

\authormark{Jir\'{a}sek \textsc{et al}}

\address[1]{\orgdiv{Faculty of Civil Engineering}, \orgname{Czech
Technical University in Prague}, \orgaddress{\state{Prague}, \country{Czech Republic}}}

\address[2]{\orgdiv{Department of Engineering}, \orgname{University of Palermo}, \orgaddress{\state{Palermo}, \country{Italy}}}

\corres{*Emma La Malfa Ribolla, Viale delle Scienze, Ed. 8. 90128 Palermo (PA). \email{emma.lamalfaribolla@unipa.it}}

%\presentaddress{This is sample for present address text this is sample for present address text}

\abstract[Summary]{The paper presents a two-dimensional geometrically nonlinear formulation of a beam element that can accommodate arbitrarily large rotations of cross sections. The formulation is based on the integrated form of equilibrium equations, which are  combined with the kinematic equations and generalized material equations, leading to a set 
of three first-order differential equations. 
These equations are then discretized by finite differences
and the boundary value problem is converted into an initial value problem using a technique inspired by the
shooting method. Accuracy of the numerical approximation is conveniently increased by 
refining the integration scheme on the element level 
while the number of global degrees of freedom is kept constant, which leads to high computational efficiency. The element has been implemented into an open-source finite element code. 
Numerical examples show a favorable comparison with standard beam elements formulated in the finite-strain framework and with  analytical solutions.}


%How to cite this article:
%\jnlcitation{\cname{%
%\author{Williams K.}, 
%\author{B. Hoskins}, 
%\author{R. Lee}, 
%\author{G. Masato}, and 
%\author{T. Woollings}} (\cyear{2016}), 
%\ctitle{A regime analysis of Atlantic winter jet variability applied to evaluate HadGEM3-GC2}, \cjournal{Q.J.R. Meteorol. Soc.}, \cvol{2017;00:1--6}.}

\maketitle

%\footnotetext{\textbf{Abbreviations:} ANA, anti-nuclear antibodies; APC, antigen-presenting cells; IRF, interferon regulatory factor}

\section{Introduction}\label{sec1}

Highly  slender  fiber-  or  rod-like  components  represent  essential  constituents  of  mechanical  systems  in  many  fields of application  such as civil, mechanical  and  biomedical  engineering. It is widely recognized that slender bodies can be efficiently modeled applying a beam theory instead of a three-dimensional continuum mechanics theory. Kirchhoff proposed the first beam formulation which includes large three-dimensional deformations \cite{kirchhoff1859,dill1992}, and Reissner completed the theory for two-dimensional  \cite{reissner1972} as well as three-dimensional cases \cite{reissner1981} with two additional deformation measures representing the shear distortion of beam segments.

Reissner's finite-strain beam theory is one of the most important geometrically nonlinear models, subsequently extended and used by many other authors for two- and three-dimensional analysis of static as well as dynamic problems. Simo developed a dynamic formulation for Reissner's beam \cite{simo1985,Auricchio} and together with Vu-Quoc \cite{simo1986} initiated the finite element implementation. He also introduced the useful concept of a \emph{geometrically  exact  beam}, based on recasting Reissner's theory in a form which is valid for any magnitude of displacements and rotations.

In this paper, a geometrically nonlinear beam model is formulated in the two-dimensional setting. This model is applicable when the rotations of beam sections become arbitrarily large and it properly accounts for the effect of curvature on the change of distance between end sections measured along the chord. Note that the axial strain is computed in a geometrically exact way while the cross section is assumed to remain rigid, planar, and perpendicular to the deformed beam axis. The material
is described by Hooke's law (an extension to a nonlinear material law would be relatively straightforward). 

In contrast to standard displacement-based finite element approaches, the proposed formulation exploits the equilibrium equations in their strong form and does not need any a priori chosen shape functions for the kinematic approximation. Based on equilibrium, the relevant internal forces (normal force and bending moment) are expressed in terms of the left-end forces and moment and the displacement and rotation functions, and then linked to the deformation variables (axial extension and curvature) using generalized material equations that describe the behavior of an infinitesimal beam segment. Substitution into the kinematic equations then leads to a set of three first-order differential equations for two displacement components and the rotation. These equations are integrated numerically, using an explicit finite difference scheme. 

On the global level, the governing equations are assembled using the standard procedure and the element is treated as a standard beam element with six degrees of freedom that represent joint displacements and rotations.
For given values of these degrees of freedom, the unknown left-end forces and moment that enter the numerical scheme are determined by local, element-level iterations that lead to satisfaction of the compatibility conditions. After that, the contribution
of the element to the nodal equilibrium equations as well as the tangent element stiffness matrix are evaluated. As a result, the beam element can remain arbitrarily long and accuracy is increased not by a reduction of the element size but by reduction of the finite difference integration step,  while keeping the number of global degrees of freedom  fixed and low.

The paper is structured as follows. Section~\ref{sec2} presents the basic assumptions and the derivation of the fundamental equations describing the beam model, which have the form of three first-order
ordinary differential equations.
An analytical solution for
the limit case of axial inextensibility
is briefly summarized (the details are provided in Appendix~\ref{app:analytical}).
Section~\ref{sec:numerics} shows how to treat
the fundamental equations numerically in a general case by an efficient procedure that
exploits the idea of the shooting method.
Numerical examples encompassing simple
one-beam problems, several frames and a honeycomb lattice are studied in Section~\ref{sec:numexamples},
and the accuracy and efficiency of the proposed
approach are evaluated by comparison to
analytical solutions and to numerical results
from the literature.
 Finally, the conclusions are summarized
 and possible extensions are discussed
 in Section~\ref{sec:conclusions}.

\section{Beam model: Governing equations and analytical solution}\label{sec2}
\subsection{Basic assumptions and variables}

\begin{figure}[h]
\centerline{\includegraphics[width=0.5 \linewidth]{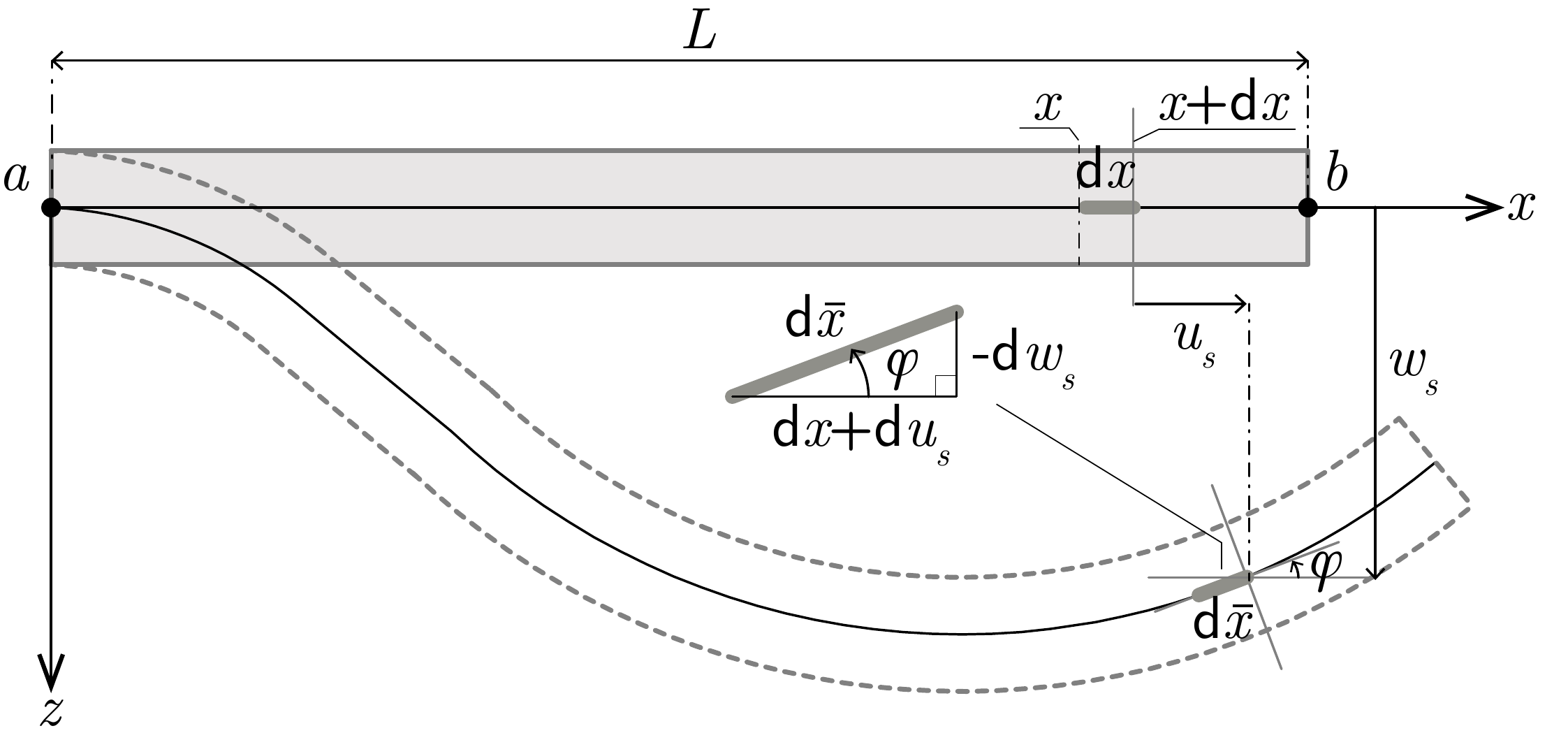}}
\caption{Kinematics of deformation of the presented nonlinear beam model.} 
\label{Fig1}
\end{figure}

Let us consider an initially straight beam of length $L$, deforming in a plane. A local coordinate system is constructed such that the origin is placed at the centroid of the ``left'' end section, the $x$-axis passes through the centroid of the ``right'' end section, and the $z$ axis is rotated by 90$^\circ$ clockwise, see Fig. \ref{Fig1}. Of course, which end is considered as the left one is a matter of choice, but once this choice is made, it is considered as fixed. The left end will be referred to by subscript $a$ and the right end by subscript $b$.

In-plane displacements and rotations of the end sections lead to deformation of the beam. During the deformation process, all cross sections are assumed to remain planar and perpendicular to the deformed beam centerline. The motion of each
section is characterized by displacements of its centroid, $u_s$ and $w_s$, and by the rotation, $\varphi$, which is taken as positive if the section rotates counterclockwise. 

The displacements $u$ and $w$ at a generic point with initial coordinates $(x,z)$ can be expressed in terms of the sectional rotation and centerline displacements as
\bea \label{eqn_aa}
u&=&u_s+z\,\sin \varphi \\ \label{eqn_ab}
w&=&w_s-z\,(1-\cos \varphi)
\eea 
These relations are nonlinear and remain valid for arbitrarily large rotations. In the spirit of the
standard beam theory, the change of distance from
the centerline caused by transversal strains is neglected, i.e., coordinate $z$ that corresponds to the signed distance from the centerline is not
adjusted.

Equations (\ref{eqn_aa})--(\ref{eqn_ab}) follow from the assumption that the cross sections remain planar. The additional assumption of perpendicularity to the deformed centerline
leads to relations
\bea\label{eqn_bb}
\sin\varphi &=&  -\frac{w'_s}{\sqrt{(1+u_s')^2+w_s'^2}} = -\frac{w'_s}{\lambda_s}
\\
\label{eqn_cc}
\cos\varphi &=&  \frac{1+u_s'}{\sqrt{(1+u_s')^2+w_s'^2}} =  \frac{1+u_s'}{\lambda_s}
\eea
which can be deduced from the geometry of the infinitesimal triangle shown in  Fig.~\ref{Fig1}.
Primes denote derivatives with respect to the axial coordinate $x$, and 
\beq\label{eqn_lams}
\lambda_s = \sqrt{(1+u_s')^2+w_s'^2}
\eeq
is the centerline stretch.
Based on (\ref{eqn_aa})--(\ref{eqn_lams}), the stretch of a generic fiber with coordinate $z$ is evaluated as
\bea\nonumber
\lambda &=& \sqrt{(1+u')^2+w'^2}=\sqrt{(1+u_s'+z\varphi'\cos\varphi)^2+(w_s'-z\varphi'\sin\varphi)^2}= 
\\ &=&
\sqrt{(\lambda_s\cos\varphi+z\varphi'\cos\varphi)^2+(-\lambda_s\sin\varphi-z\varphi'\sin\varphi)^2}= \lambda_s+z\varphi'
\label{eq6}
\eea
This means that the stretch, and thus also the Biot strain, defined as $\eps_B=\lambda-1$, varies across the depth of the section in a linear fashion. In contrast to that, the Green-Lagrange strain, $\eps_{GL}=(\lambda^2-1)/2$, would be described by a quadratic function of $z$.

The description in terms of Biot strain leads to simpler equations, and so, as a prototype linear elastic model, we will develop the governing equations based on the strain energy density
\beq 
\mathcal E_{int}(\lambda) = \half E(\lambda-1)^2
\eeq 
considered as a quadratic function of the Biot strain, with parameter $E$ representing the Young modulus. Conceptually, there is no problem with
replacement of this assumption by another hyperelastic law or even by an inelastic stress-strain law, if needed.

It is worth noting that since shear distortion
is neglected here and the stress-strain law on the level of each fiber is essentially uniaxial, the role of the stress that is work-conjugate to the Biot strain is played
by the first Piola-Kirchoff stress. In what follows, the Biot strain will be denoted
simply by 
\beq 
\eps=\lambda-1
\eeq

\subsection{Variational derivation of equilibrium equations}

The equilibrium state can be found by exploiting the principle of minimum potential energy. The total potential energy,
\beq\label{eq:Ep} 
 E_{p} = E_{int}+E_{ext}
\eeq 
is the sum of the strain energy, $E_{int}$, and the energy of external forces, $E_{ext}$. In this work, we neglect body forces and we consider the beam structure to be loaded only at its joints.
The strain energy of one beam,
\beq 
E_{int} = \int_0^L \int_A \mathcal E_{int}(\lambda) \,\mathrm{d} A\,\mathrm{d} x
\eeq 
is calculated by integrating the strain energy density over the volume of the beam.

For prescribed joint displacements and rotations, the energy of external forces vanishes and the equilibrium state is found by minimizing functional $E_{int}$ over all kinematically admissible states that are characterized by
functions $u_s$, $w_s$ and $\varphi$ satisfying the kinematic boundary conditions and the perpendicularity constraint  expressed by equations (\ref{eqn_bb})--(\ref{eqn_cc}).
The first variation of the beam strain energy is evaluated as
 \beq \label{eq11}
\delta E_{int} =  \int_0^L \int_A \frac{\mathrm{d}\mathcal E_{int}}{\mathrm{d}\lambda}\delta\lambda \,\mathrm{d} A\,\mathrm{d} x =
\int_0^L \int_A 
\sigma\,(\delta\lambda_s+z\,\delta\varphi') \,\mathrm{d}A\,\mathrm{d}x =
\int_0^L (N\,\delta\lambda_s+M\,\delta\kappa) \,\mathrm{d}x
\eeq
where 
\beq 
\sigma = \frac{\mathrm{d}\mathcal{E}_{int}}{ \mathrm{d}\eps}= \frac{\mathrm{d}\mathcal{E}_{int}}{ \mathrm{d}\lambda}
\eeq
is the stress work-conjugate with the Biot strain,
and
\bea 
\label{eee1}
N&=&\int_A \sigma\,\mathrm{d}A \\
\label{eee2}
M&=&\int_A z\sigma\,\mathrm{d}A 
\eea 
are the normal force and the bending moment.

For the material model based on strain energy density taken as a quadratic function of the Biot strain, the stress can be expressed as
\beq 
\sigma =E(\lambda-1) = E\,(\lambda_s-1+z\varphi') = E\eps_s+zE\kappa
\eeq
where
\beq\label{eq:epss}
\eps_s = \lambda_s-1 = \sqrt{(1+u_s')^2+w_s'^2}-1
\eeq 
is the strain at the centerline and
\beq \label{eq:kappa}
\kappa=\varphi'
\eeq
is the curvature. Since the stress is linearly distributed across the depth of the section, one can evaluate the integrals in (\ref{eee1})--(\ref{eee2}) analytically and derive the standard relations between the internal forces and the deformation variables:
\bea 
\label{e1a}
N&=&\int_A (E\eps_s+zE\kappa)\,\mathrm{d}A = EA\eps_s \\
\label{e2a}
M&=&\int_A (zE\eps_s+z^2E\kappa)\,\mathrm{d}A =EI\kappa
\eea 
Here, $A$ is the sectional area and $I$ is the sectional moment of inertia.

To proceed from the first variation $\delta E_{int}$ described by (\ref{eq11}) to the stationarity conditions for functional $E_{int}$, we have to realize that $\lambda_s$ and $\kappa$ are not the primary independent fields, and so their variations need to be expressed in terms of the centerline displacement variations, $\delta u_s$ and $\delta w_s$. It turns out to be convenient to keep working for a while with the rotation, $\varphi$, and its variation, $\delta\varphi$, as  auxiliary fields which can later be expressed in terms of the primary ones.
From (\ref{eqn_lams}) and (\ref{eq:kappa}), one gets
\bea
\delta\lambda_s &=&
\frac{(1+u_s')\,\delta u_s' + w_s'\,\delta w_s'}{\lambda_s} = \cos\varphi\,\delta u_s' -\sin\varphi\,\delta w_s' \\
\delta\kappa &=&\delta\varphi'
\eea
Substitution of these expressions into (\ref{eq11}) and integration by parts leads to
\bea\nonumber
\delta E_{int} &=& 
 \int_0^L \left(N\,\cos\varphi\,\delta u_s'-N\,\sin\varphi\,\delta w_s'+M\,\delta\varphi'\right) \,\dx = \\
 &=& 
 \left[N\,\cos\varphi\,\delta u_s-N\,\sin\varphi\,\delta w_s+M\,\delta\varphi \right]_0^L
 -\int_0^L \left((N\,\cos\varphi)'\delta u_s-(N\,\sin\varphi)'\delta w_s+M'\,\delta\varphi\right) \,\dx
 \label{eq22}
\eea

If the displacements and rotations of the end sections are considered as prescribed, their variations are zero and the boundary terms in (\ref{eq22}) vanish.
Inside the beam, the variation of the rotation is not independent of the displacement variations, because of the perpendicularity constraint. Based on (\ref{eqn_bb})--(\ref{eqn_cc}), it is possible to show that
\beq 
\delta\varphi = -\frac{\sin\varphi\,\delta u_s'+\cos\varphi\,\delta w_s'}{\lambda_s}
\eeq 
and, consequently,
\bea\nonumber
-\int_0^L M'\,\delta\varphi\,\dx &=&  \int_0^L M'\,\frac{\sin\varphi\,\delta u_s'+\cos\varphi\,\delta w_s'}{\lambda_s}\,\dx = 
\\ &=&
\left[ M'\,\frac{\sin\varphi\,\delta u_s+\cos\varphi\,\delta w_s}{\lambda_s}\right]_0^L
-\int_0^L \left(\frac{M'}{\lambda_s}\sin\varphi\right)'\delta u_s\,\dx
-\int_0^L \left(\frac{M'}{\lambda_s}\cos\varphi\right)'\delta w_s\,\dx
\eea
Taking into account that the boundary terms vanish and
substituting back into (\ref{eq22}), we finally obtain
\beq
\delta E_{int} = -\int_0^L \left((N\,\cos\varphi)'\delta u_s-(N\,\sin\varphi)'\delta w_s+\left(\frac{M'}{\lambda_s}\sin\varphi\right)'\delta u_s+ \left(\frac{M'}{\lambda_s}\cos\varphi\right)'\delta w_s\right) \,\dx
 \label{eq22x}
\eeq

\begin{figure}[h]
\centering
\begin{tabular}{ccc}
(a) & (b) & (c) \\
\includegraphics[width=0.32 \linewidth]{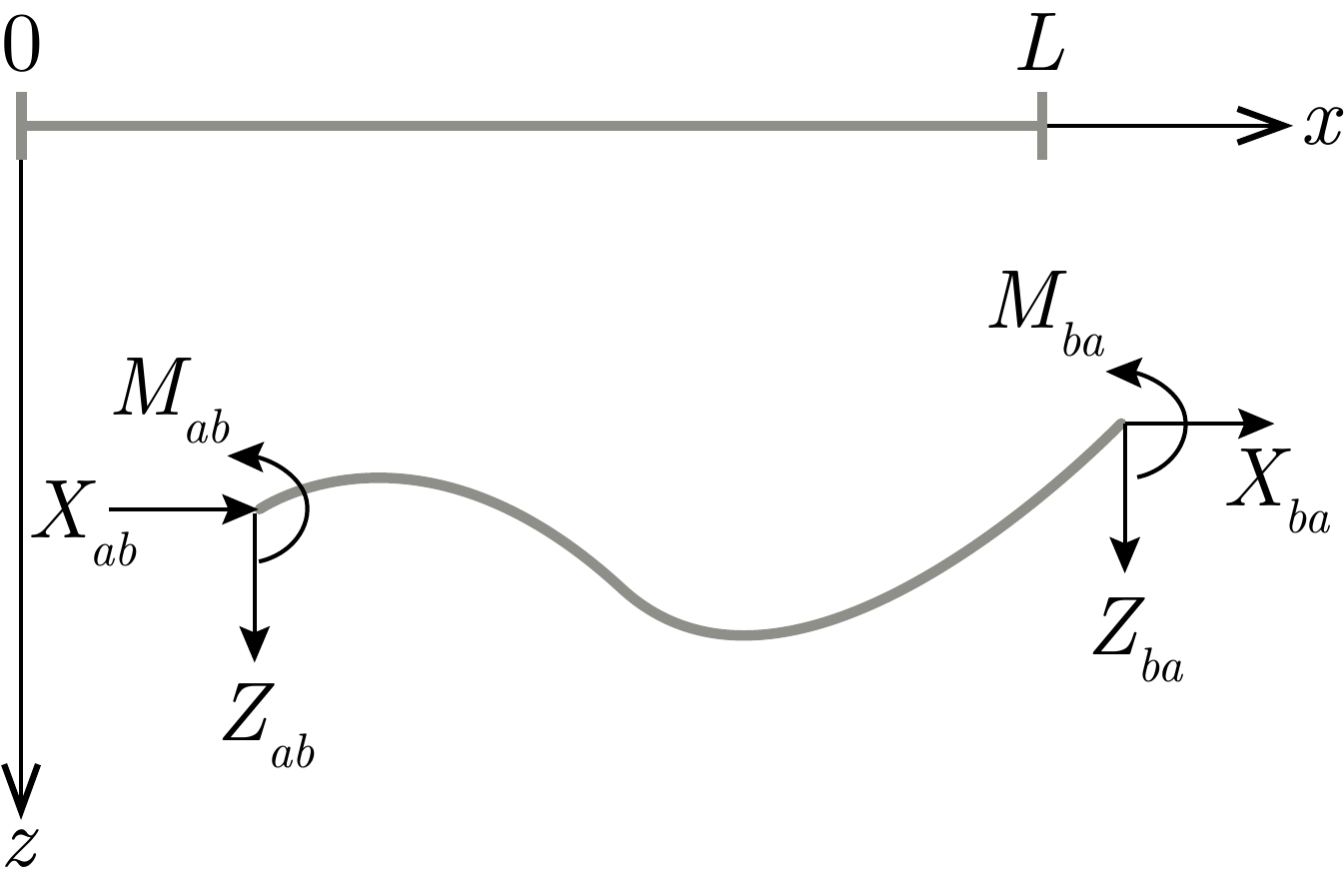}
&
\includegraphics[width=0.32 \linewidth]{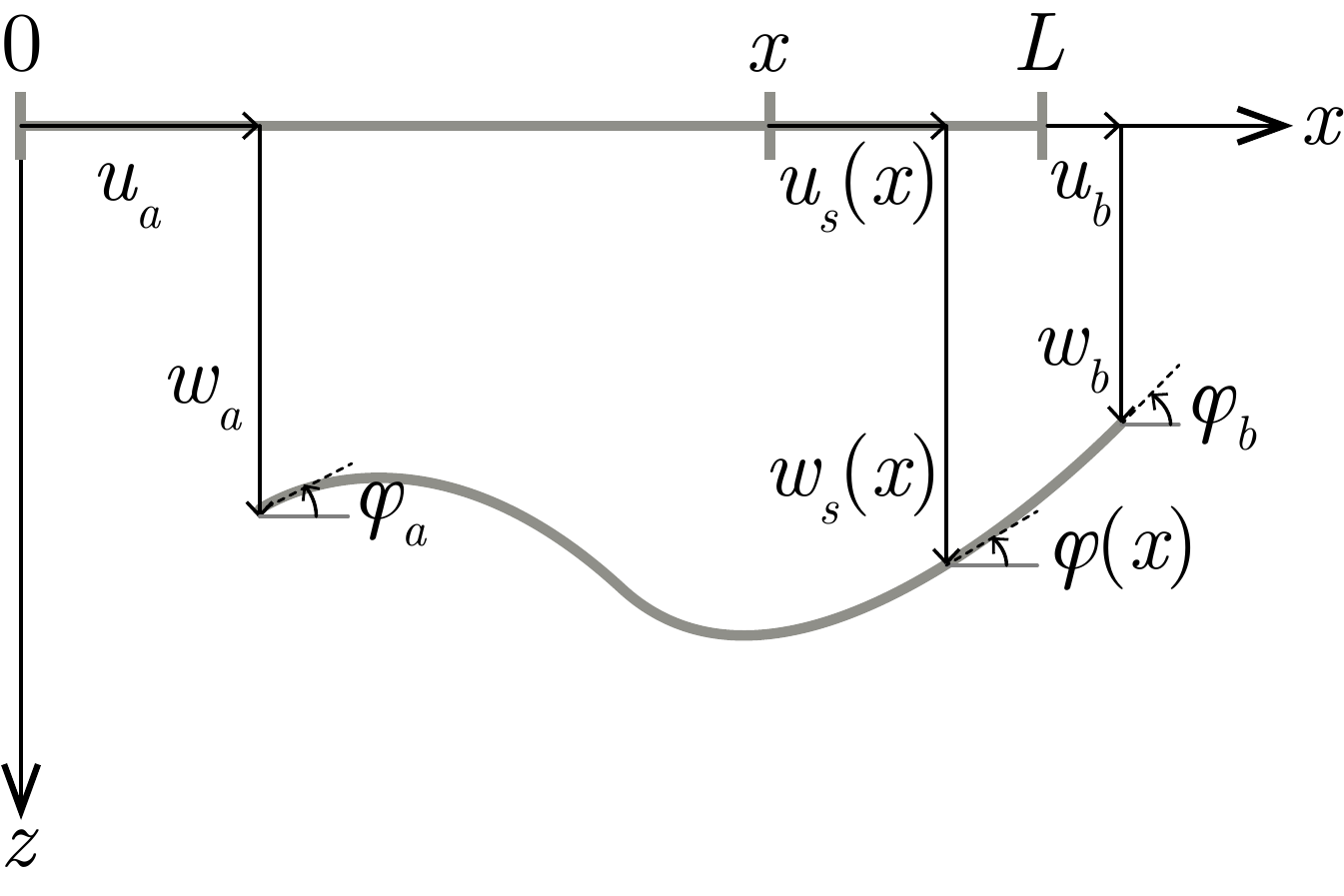}
&
\includegraphics[width=0.32 \linewidth]{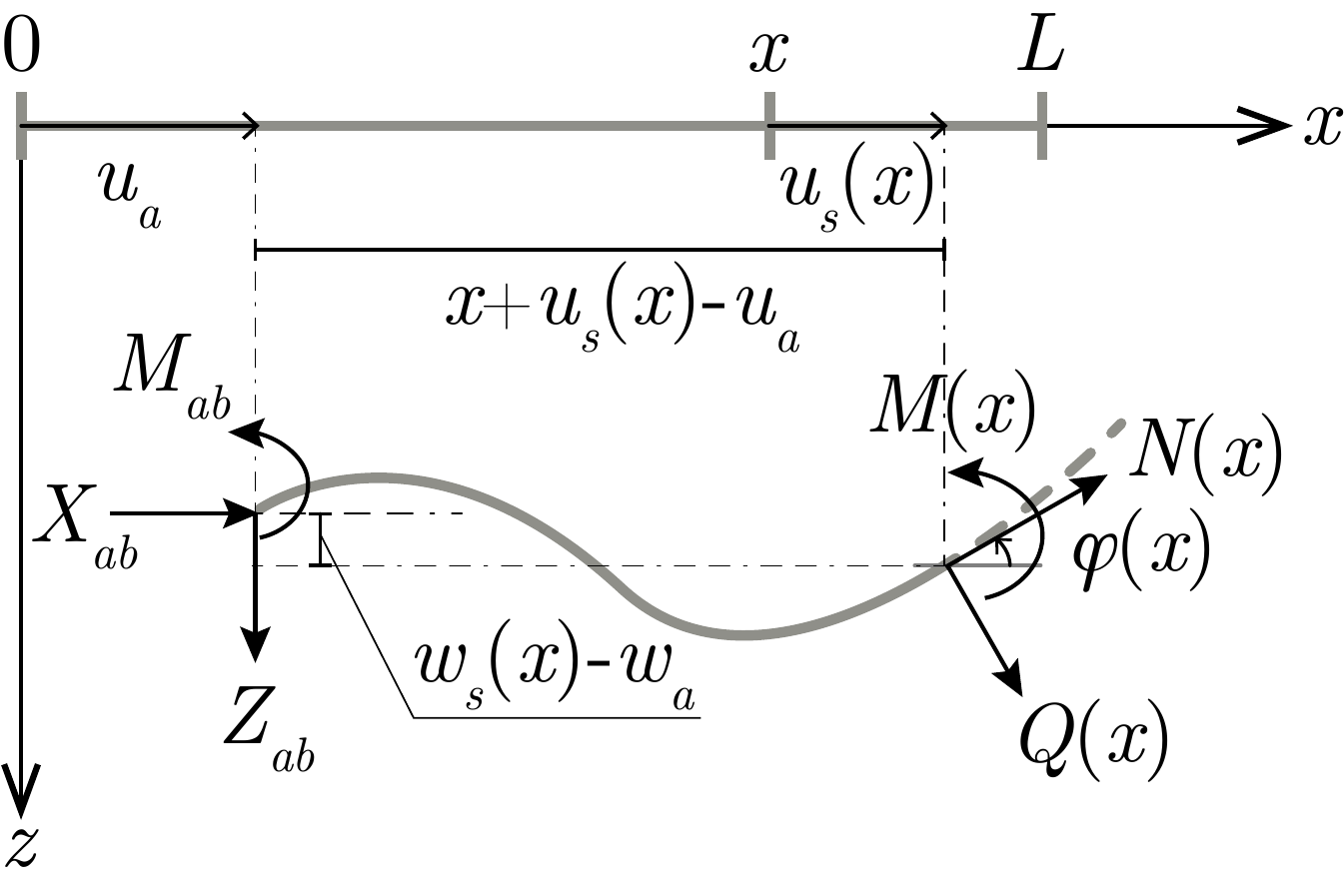}
\end{tabular}
\caption{Deformed beam: (a) end forces and moments, (b) displacements and rotations, (c) free-body diagram with left-end forces and internal forces at a generic section.} 
\label{Fig2}
\end{figure}

The stationarity condition $\delta E_{int}=0$ for all admissible variations $\delta u_s$ and $\delta w_s$ yields  differential equations
\bea \label{eq26}
-(N\,\cos\varphi)'-\left(\frac{M'}{\lambda_s}\sin\varphi\right)' &=& 0
\\
(N\,\sin\varphi)'-\left(\frac{M'}{\lambda_s}\cos\varphi\right)'&=&0
\label{eq27}
\eea 
which represent the strong form of equilibrium equations. Due to their special form, it is possible to perform closed-form integration and write
\bea \label{eq28}
-N\cos\varphi - \frac{M'}{\lambda_s}\sin\varphi &=& X_{ab} \\
N\sin\varphi - \frac{M'}{\lambda_s}\cos\varphi &=& Z_{ab}
\label{eq29}
\eea 
where $X_{ab}$ and $Z_{ab}$ are integration constants, which physically correspond to the components of the left-end force, i.e., the force acting between the left end section of the beam and the joint to which this section is attached (Fig.~\ref{Fig2}a).

From relations (\ref{eq28})--(\ref{eq29}), it is possible to express
\bea \label{eq30}
N &=& -X_{ab}\cos\varphi + Z_{ab}\sin\varphi \\
\frac{M'}{\lambda_s} &=& -X_{ab}\sin\varphi - Z_{ab}\cos\varphi
\label{eq31}
\eea 
Based on (\ref{eqn_bb})--(\ref{eqn_cc}), equation (\ref{eq31}) can be recast as
\beq 
M' = X_{ab}w_s' - Z_{ab}(1+u_s')
\eeq 
and yet another integration leads to
\beq \label{eq33}
M = -M_{ab}+X_{ab}(w_s-w_a) - Z_{ab}(x+u_s-u_a)
\eeq 
where the integration constant $M_{ab}$
represents the left-end moment (Fig.~\ref{Fig2}a), and
$u_a=u_s(0)$ and $w_a=w_s(0)$ are the displacements
at the left end (Fig.~\ref{Fig2}b).

The fraction on the left-hand side of (\ref{eq31}) physically corresponds to the shear force, $Q$, which plays here only an auxiliary role and is not linked to any deformation variable by a constitutive law, because the shear distortion is neglected. 
Of course, equations (\ref{eq30})--(\ref{eq31}) and (\ref{eq33}) could be constructed as equilibrium equations from a free-body diagram, as illustrated in Fig.~\ref{Fig2}c. The present derivation shows that they can be consistently derived by closed-form integration of stationarity conditions obtained from the principle of minimum potential energy.
These equations properly take into account geometric effects and remain accurate for arbitrarily large rotations.

\subsection{Fundamental equations of the present approach}

\subsubsection{General case}
\label{sec:fund1}

In the standard displacement-based approach, relations (\ref{e1a})--(\ref{e2a}) that link the internal forces to the deformation variables, combined with an expression for the rotation derived from (\ref{eqn_bb}) or (\ref{eqn_cc})
and with relations (\ref{eq:epss})--(\ref{eq:kappa}) that express the deformation variables in terms of the centerline displacement functions, would be substituted into the 
differential equations of equilibrium (\ref{eq26})--(\ref{eq27}). As an alternative, one can start from equations
\bea 
\varphi' &=& \kappa \\
\label{eq35}
u_s' &=& \lambda_s\cos\varphi - 1\\
\label{eq36}
w_s' &=& -\lambda_s\sin\varphi
\eea 
which easily follow from (\ref{eqn_bb})--(\ref{eqn_cc}) and (\ref{eq:kappa}).
The centerline displacements $u_s$ and $w_s$ as well as the rotation $\varphi$ are considered here as primary unknown functions that will be computed by integration of the above
first-order differential equations. To this end, we must express the deformation variables $\kappa$ and $\lambda_s$ on the right-hand sides in terms of the primary variables,
which can be achieved by combining the inverted form of equations (\ref{e1a})--(\ref{e2a}) that link the internal forces to the deformation variables
with the integrated equilibrium equations (\ref{eq30}) and (\ref{eq33}).
The resulting equations read
\bea \label{e212}
\varphi' &=& \frac{-M_{ab}+X_{ab} (w_s-w_a) - Z_{ab}(x+u_s-u_a)}{EI}\\
\label{e213}
u_s' &=& \left(1+\frac{-X_{ab}\cos\varphi + Z_{ab}\sin\varphi}{EA}\right)\cos\varphi-1 \\
\label{e214}
w_s' &=& -\left(1+\frac{-X_{ab}\cos\varphi + Z_{ab}\sin\varphi}{EA}\right)\sin\varphi
\eea 
 and they indeed form a set of three first-order differential  equations for three unknown functions.
 
Equations (\ref{e212})--(\ref{e214}) are considered as the fundamental equations of the present approach.
Interestingly,  they could be reduced to a single second-order differential equation  for the unknown rotation. Differentiating (\ref{e212}) and substituting from (\ref{e213})--(\ref{e214}), we obtain
\beq \label{e221}
EI\varphi'' + \left(1+\frac{-X_{ab}\cos\varphi + Z_{ab}\sin\varphi}{EA}\right)
\left(X_{ab}\sin\varphi + Z_{ab}\cos\varphi\right)=0
\eeq

In addition to the primary unknown functions and given sectional stiffnesses $EA$ and $EI$, the fundamental equations (\ref{e212})--(\ref{e214})
contain constants $X_{ab}$, $Z_{ab}$ and $M_{ab}$,  which are  usually unknown. Integration of (\ref{e212})--(\ref{e214}) generates three additional integration constants. In total,
we have six unknown constants that can be determined from six boundary conditions (three at each end section). In problems that involve analysis of a single beam, the structure of boundary conditions depends on the way the beam is supported. This is illustrated by an example in Appendix~\ref{sec:examplecantilever}. On the other hand, in the context of structural analysis of a frame, the joint displacements and rotations play the role of global unknowns that are determined by iteratively solving the joint equilibrium equations. 
On the beam element level, the basic tasks are (1) to compute the end forces and moments that correspond to prescribed values of the end displacements and rotations, and (2) to evaluate the corresponding element tangent stiffness matrix. Numerical procedures for task 1 will be elaborated in Section~\ref{sec:endforces} and for task 2 in Section~\ref{sec:stiff}. 

\subsubsection{Special cases: axial inextensibility and moderate rotations}

A special case is the {\bf axially inextensible/incompressible} beam model,
characterized by $EA\to\infty$ and $\eps_s=0$. The fundamental equations (\ref{e212})--(\ref{e214})
then reduce to
\bea \label{e215}
\varphi' &=& \frac{-M_{ab}+X_{ab} (w_s-w_a) - Z_{ab}(x+u_s-u_a)}{EI}\\
u_s' &=& \cos\varphi-1 \\
\label{e217}
w_s' &=& -\sin\varphi
\eea 
and equation (\ref{e221}) reduces to
\beq \label{e221x}
EI\varphi'' + 
X_{ab}\sin\varphi + Z_{ab}\cos\varphi=0
\eeq
The axially inextensible model can be treated analytically; see Appendix~\ref{app:analytical} and Section~\ref{sec:analytical}. However, in the context of general frame analysis, the analytical approach would lead to numerical problems, because the combinations of end displacements would not be completely arbitrary (they would be restricted by an inequality resulting from the incompressibility constraint) and the joint displacements could not be considered as unconstrained unknowns.

For {\bf small or moderate rotations}, the exact equations could be approximated. For instance, keeping only terms up to the first order in $\varphi$, we can replace  (\ref{e213})--(\ref{e214}) by
\bea 
\label{e213x}
u_s' &=& \frac{-X_{ab} + Z_{ab}\,\varphi}{EA} \\
\label{e214x}
w_s' &=& -\left(1-\frac{X_{ab} }{EA}\right)\varphi
\eea 
and (\ref{e221}) by
\beq \label{e221y}
EI\varphi'' + \left(1-\frac{X_{ab}}{EA}\right)Z_{ab}
+\left(X_{ab}-\frac{X_{ab}^2-Z_{ab}^2}{EA}\right)\varphi
=0
\eeq

One needs to be careful when combining {\bf small rotations with axial inextensibility}. Setting $EA\to\infty$ in (\ref{e213x}), we would obtain $u_s'=0$. However, if $EA\to\infty$ is used in (\ref{e213}), the equation reduces to $u_s'=\cos\varphi-1$ and the approximation for {\bf moderate rotations} should keep a quadratic term. The resulting equations are then
\bea 
\label{e213z}
u_s' &=& -\half\varphi^2 \\
\label{e214z}
w_s' &=& -\varphi \\
\label{e221z}
EI\varphi'' + Z_{ab}
+X_{ab}\varphi &=& 0
\eea 
Differentiating (\ref{e221z}) and substituting $\varphi=-w_s'$ according to (\ref{e214z}), we end up with the well-known equation describing buckling of an axially compressed straight beam,
\beq 
EIw_s^{IV} +X_{ab}w_s'' = 0
\eeq 
Equation (\ref{e213z}) can then be used to estimate the relative displacement of the beam ends caused by second-order effects,
\beq 
\Delta L = u(L)-u(0) = \int_0^L u_s'\,{\rm d}x = -\half \int_0^L w_s'^2\,{\rm d}x
\eeq 

\subsection{Analytical solution}\label{sec:analytical}
Interestingly, the fundamental equations in their reduced form (\ref{e215})--(\ref{e221x}), valid for the inextensible case, admit an analytical solution in terms of elliptic functions and elliptic integrals.
The derivation of this solution is presented in Appendix~\ref{app:analytical}. For a beam segment without an inflexion point, the resulting expressions for the rotation and displacements have the form
\bea\label{eq:284z} 
\varphi(x) &=& 2\,\arcsin(\tilde{k}\,{\rm sn}(\tilde{a}+\tilde{b}x,\tilde{k}))-\alpha 
\\
\label{eq288z}
u_s(x) &=& C_u-(1+\cos\alpha)\,x-\frac{2\tilde{k}\sin\alpha}{\tilde{b}}{\rm cn}(\tilde{a}+\tilde{b}x,\tilde{k})+\frac{2\cos\alpha}{\tilde{b}}\,E_J({\rm am}(\tilde{a}+\tilde{b}x,\tilde{k}),\tilde{k})
\\
\label{eq287z}
w_s(x) &=& C_w-(\sin\alpha)\,x+ \frac{2\tilde{k}\cos\alpha}{\tilde{b}}{\rm cn}(\tilde{a}+\tilde{b}x,\tilde{k})+\frac{2\sin\alpha}{\tilde{b}}\,E_J({\rm am}(\tilde{a}+\tilde{b}x,\tilde{k}),\tilde{k})
\eea 
where ``sn'' and ``cn'' are the elliptic sine and cosine, ``am'' is the Jacobi amplitude function, and $E_J$ is the incomplete elliptic integral of the second kind. The relation of constants $\tilde k$, $\tilde a$, $\tilde b$, $\alpha$, $C_w$, and $C_u$
to the beam properties ($EI$ and $L$) and boundary conditions is described in detail in Appendix~\ref{app:analytical}.

Based on the general solution, it is possible to derive analytical expressions for a cantilever loaded at its free end by an arbitrarily inclined force. As shown in Appendix~\ref{app:analytical}, the applied
force, $F$, and the displacements of the left end of a cantilever fixed at its right end, $u_a$ and $w_a$, can be expressed in terms of the left-end rotation, $\varphi_a$; see formulae (\ref{eq:ex298})--(\ref{eq:ex300}).
These expressions will later be used as benchmarks.
Nevertheless, the analytical or semi-analytical approach is applicable only to simple cases, and general frame analysis needs to be based on numerical methods, which will be developed in the following section.

\section{Numerical procedures}
\label{sec:numerics}

\subsection{Evaluation of end forces and moments}\label{sec:endforces}

Analytical formulae such as
(\ref{eq:284z})--(\ref{eq287z}) are useful only if the elliptic functions and elliptic integrals are already implemented by efficient algorithms.
Moreover, these analytical solutions are valid only under the restrictive assumption of axial inextensibility.
A more flexible and straightforward approach is to construct approximate solutions of the fundamental differential equations using standard numerical procedures. 
Numerical treatment will be based on the full form of fundamental equations (\ref{e212})--(\ref{e214}), because inextensibility would lead to numerical problems (e.g., infinite axial stiffness for a straight beam under tension) and the assumption of small or moderate rotations would induce a large error if the beams deform substantially.

Suppose that the displacements and rotations of the end section of a beam element are prescribed. It is convenient to decompose the motion of the beam into (A) the rigid-body motion dictated by the displacements and rotation of the left end and (B) the deformation of the beam (stretching and bending of the beam centerline) during which the left end remains fixed.
Phase A is easy to handle as a simple geometric transformation, and so we focus first on phase B, leaving the implementation of phase A to Section~\ref{sec:trans}.

\subsubsection{Shooting method}
\label{sec:shoot}

The first partial task is to evaluate the right-end displacements $u_b$, $w_b$ and $\varphi_b$ if the left-end displacements $u_a$, $w_a$ and $\varphi_a$ and the left-end forces $X_{ab}$, $Z_{ab}$ and $M_{ab}$ are given. 
In phase B, 
the rotation $\varphi$ and displacements $u_s$ and $w_s$ in (\ref{e212})--(\ref{e214})  are taken with respect to a co-rotational coordinate system attached to the left end section, and so the 
conditions to be imposed at the left end read
\bea 
\varphi(0) &=& 0 \\
u_s(0) &=& 0\\
w_s(0) &=& 0
\eea 
They can be understood as initial conditions that make the solution of differential equations (\ref{e212})--(\ref{e214}) unique, provided that the left-end forces $X_{ab}$ and $Z_{ab}$ and moment $M_{ab}$ are known. The solution can be constructed numerically, using a suitable finite difference scheme.

The interval $[0,L]$ is divided into $N$ numerical segments of length $\Delta x=L/N$, with grid points $x_i=i\,\Delta x$, $i=0,1,2,\ldots N$, and approximate values of the rotation, centerline displacements and internal forces at these grid points 
are denoted as $\varphi_i$, $u_i$, $w_i$, $N_i$ and $M_i$, $i=0,1,2,\ldots N$.
The derivatives in (\ref{e212})--(\ref{e214}) are replaced by finite differences. 
The simplest approach is based on the following explicit scheme:
\begin{enumerate}
    \item Set initial values $\varphi_0=0$, $u_0=0$ and $w_0=0$.
    \item For $i=1,2,\ldots N$ evaluate
\bea \label{e234}
M_{i-1} &=& -M_{ab}+X_{ab}w_{i-1}-Z_{ab}(x_{i-1}+u_{i-1}) \\ \label{e235}
\varphi_{i-1/2} &=& \varphi_{i-1} + \frac{M_{i-1}}{EI}\frac{\Delta x}{2} \\
N_{i-1/2} &=&  -X_{ab}\cos\varphi_{i-1/2} + Z_{ab}\sin\varphi_{i-1/2} \\ \label{e238a}
u_i &=& u_{i-1}+\left[\left(1+\frac{N_{i-1/2}}{EA}\right)\cos\varphi_{i-1/2}-1\right] \,\Delta x\\ \label{e238b}
w_i &=& w_{i-1}-\left(1+\frac{N_{i-1/2}}{EA}\right)\sin\varphi_{i-1/2}\,\Delta x \\ \label{e239}
M_i &=& -M_{ab}+X_{ab}w_i-Z_{ab}(x_i+u_i) \\ 
\varphi_i &=& \varphi_{i-1/2} + \frac{M_i}{EI}\frac{\Delta x}{2}\label{e240}
\eea 
\item The resulting displacement and rotation values at the right end are $u(L)=u_N$, $w(L)=w_N$ and $\varphi(L)=\varphi_N$.
\end{enumerate}
As indicated in (\ref{e235}) and (\ref{e240}), the rotation is integrated in two half-steps, one of them based on the curvature at $x_{i-1}$ and the other at $x_i$. The first half-step allows to get an approximation of the rotation at midstep, $\varphi_{i-1/2}$, which is then exploited for evaluation of the normal force and centerline strain at midstep and to integration of the centerline displacement in one single step based on the central difference scheme. This allows evaluation of the curvature at the end of the step, and thus the second half-step for the integration of the rotation remains explicit, even though it is based on the backward finite difference formula.
For simplicity, the sectional stiffnesses $EI$ and $EA$ are considered as constant, but it would be straightforward to extend the algorithm to beams with variable section. In this case, $EA$ in (\ref{e238a})--(\ref{e238b}) would be replaced by $EA_{i-1/2}$ while $EI$ in (\ref{e235}) and (\ref{e240}) would be replaced by $EI_{i-1}$ and $EI_i$, respectively.

Of course, the left-end forces and moment,  $X_{ab}$, $Z_{ab}$ and $M_{ab}$, which are needed to run the algorithm, are not known in advance. If we somehow estimate their values and prescribe zero initial values of the kinematic quantities (as specified in step 1), we can run the algorithm and determine the values of right-end displacements and rotation, $u_s(L)$, $w_s(L)$ and $\varphi(L)$.
The values of the left-end forces and moment then need to be adjusted such that the resulting kinematic quantities at the right end satisfy the yet unused boundary conditions
\bea 
\varphi(L) &=& \varphi_b \\
u_s(L) &=& u_b \\
w_s(L) &=& w_b
\eea 
in which $u_b$, $w_b$ and $\varphi_b$ are prescribed displacements and rotation of the right end with respect to the left end that arise during phase B of the deformation process (after rigid-body motion A during which the whole beam translates and rotates with its left end).

In fact, the suggested approach is a special version of the shooting method. For a given set of end displacements and rotations, the initial estimate of $M_{ab}$, $X_{ab}$ and $Z_{ab}$ can be constructed based on linear beam theory, or on the values at the end of the previous step if the calculation is done in the context of an incremental iterative structural analysis.

The foregoing algorithm defines a certain mapping of the left-end forces and moment on the right-end displacements and rotation. Formally we can write
\beq\label{e241}
\boldsymbol{u}_b = \boldsymbol{g}(\boldsymbol{f}_{ab})
\eeq 
where
\beq 
%\boldsymbol{u}_a = \left(\begin{array}{c} u_a \\ w_a \\ \varphi_a \end{array}\right), \hskip 10mm
\boldsymbol{u}_b = \left(\begin{array}{c}  u_b \\  w_b \\ \varphi_b \end{array}\right), \hskip 10mm
\boldsymbol{f}_{ab} = \left(\begin{array}{c} X_{ab} \\ Z_{ab} \\ M_{ab} \end{array}\right)
\eeq 
For a given column matrix $\boldsymbol{u}_b$, equation (\ref{e241}) represents a set
of three nonlinear equations for unknowns
collected in column matrix $\boldsymbol{f}_{ab}$. The solution is found by the
Newton-Raphson method, using the recursive formula
\beq 
\delta\boldsymbol{f}_{ab}^{(k)} = \boldsymbol{G}^{-1}\left(\boldsymbol{f}_{ab}^{(k)}\right)\left(\boldsymbol{u}_b-\boldsymbol{g}\left(\boldsymbol{f}_{ab}^{(k)}\right)\right) , \hskip 10mm \boldsymbol{f}_{ab}^{(k+1)}=\boldsymbol{f}_{ab}^{(k)}+\delta\boldsymbol{f}_{ab}^{(k)}, \hskip 10mm k=0,1,2,\ldots
\eeq 
where 
\beq 
\boldsymbol{G} = \frac{\partial \boldsymbol{g}}{\partial \boldsymbol{f}_{ab}}
\eeq 
is the Jacobi matrix of mapping $\boldsymbol{g}$.
%with respect to variables $\boldsymbol{f}_{ab}$.

The entries of the Jacobi matrix are evaluated numerically using the differentiated version of the computational scheme.
Suppose that the input values $X_{ab}$, $Z_{ab}$ and $M_{ab}$ are changed by infinitesimal increments
$\dif X_{ab}$, $\dif Z_{ab}$ and $\dif M_{ab}$.
Linearization of equations (\ref{e234})--(\ref{e240}) around the currently considered solution leads to 
\bea \label{e234x}
\dif \varphi_{i-1/2} &=& \dif \varphi_{i-1} + \frac{\Delta x}{2EI}\left[-\dif M_{ab}+\dif X_{ab}w_{i-1}+X_{ab}\dif w_{i-1}-\dif Z_{ab}(x_{i-1}+u_{i-1})-Z_{ab}\dif u_{i-1}\right] \\
\dif N_{i-1/2} &=&  -\dif X_{ab}\cos\varphi_{i-1/2}+ X_{ab}\sin\varphi_{i-1/2}\dif \varphi_{i-1/2}+ \dif Z_{ab}\sin\varphi_{i-1/2}+Z_{ab}\cos\varphi_{i-1/2}\dif \varphi_{i-1/2} \\
\dif u_i &=& \dif u_{i-1}+\frac{\dif N_{i-1/2}}{EA}\cos\varphi_{i-1/2}\,\Delta x-\left(1+\frac{N_{i-1/2}}{EA}\right)\sin\varphi_{i-1/2}\dif\varphi_{i-1/2} \,\Delta x\\
\dif w_i &=& \dif w_{i-1}-\frac{\dif N_{i-1/2}}{EA}\sin\varphi_{i-1/2}\,\Delta x-\left(1+\frac{N_{i-1/2}}{EA}\right)\cos\varphi_{i-1/2}\dif\varphi_{i-1/2}\,\Delta x \\
\dif\varphi_i &=& \dif\varphi_{i-1/2} + \frac{\Delta x}{2EI}\left[-\dif M_{ab}+\dif X_{ab}w_i+X_{ab}\dif w_i-\dif Z_{ab}(x_i+u_i)-Z_{ab}\dif u_i\right]
\label{e240x}
\eea 
The values of $\dif u_0$, $\dif w_0$ and $\dif\varphi_0$ are set to zero, because the initial zero values of $u_0$,
$w_0$ and $\varphi_0$ are fixed and remain 
unaffected by changes of $X_{ab}$, $Z_{ab}$ and $M_{ab}$.

If we set $\dif X_{ab}=1$ and $\dif Z_{ab}=\dif M_{ab}=0$,
the resulting values of $\dif u_N$, $\dif w_N$ and $\dif\varphi_N$ will correspond to the first column of the Jacobi matrix. They are evaluated using the adapted scheme
\bea \label{e234y}
\dif \varphi_{i-1/2} &=& \dif \varphi_{i-1} + \frac{\Delta x}{2EI}\left(w_{i-1}+X_{ab}\dif w_{i-1}-Z_{ab}\dif u_{i-1}\right) \\
\dif N_{i-1/2} &=&  -\cos\varphi_{i-1/2}+ X_{ab}\sin\varphi_{i-1/2}\dif \varphi_{i-1/2}+ Z_{ab}\cos\varphi_{i-1/2}\dif \varphi_{i-1/2} \\
\dif u_i &=& \dif u_{i-1}+\frac{\dif N_{i-1/2}}{EA}\cos\varphi_{i-1/2}\,\Delta x-\left(1+\frac{N_{i-1/2}}{EA}\right)\sin\varphi_{i-1/2}\dif\varphi_{i-1/2}  \,\Delta x\\
\dif w_i &=& \dif w_{i-1}-\frac{\dif N_{i-1/2}}{EA}\sin\varphi_{i-1/2}\,\Delta x-\left(1+\frac{N_{i-1/2}}{EA}\right)\cos\varphi_{i-1/2}\dif\varphi_{i-1/2}\,\Delta x \\
\dif\varphi_i &=& \dif\varphi_{i-1/2} + \frac{\Delta x}{2EI}\left(w_i+X_{ab}\dif w_i-Z_{ab}\dif u_i\right)\label{e240y}
\eea 
The additional two columns of the Jacobi matrix are obtained in an analogous fashion,
setting $\dif Z_{ab}=1$ or $\dif M_{ab}=1$.

\subsubsection{Transformation to global coordinates}
\label{sec:trans}

Suppose that the shooting method described in the previous section  has been implemented.
The computed displacements $\boldsymbol{u}$
as well as the end forces $\boldsymbol{f}_{ab}$
are expressed in an auxiliary coordinate system $xz$ with the origin located at the left end of the beam in the deformed configuration 
%(point $\tilde a$ in Fig.~\ref{f:xtilde})
and with the $x$ axis in the direction of the tangent to the
deformed centerline at the left end. Now we would like
to link them to the components expressed with respect
to the global coordinate system, which
will be denoted by a superscript $G$. 

The initial geometry is described by global coordinates
of the joints connected by the beam, i.e., 
$x_{a}^G$ and $z_{a}^G$
at the left end and $x_{b}^G$ and $z_{b}^G$ at the right end,
from which we can compute the beam length
\beq 
L_{ab} = \sqrt{(x_{b}^G-x_{a}^G)^2+(z_{b}^G-z_{a}^G)^2}
\eeq
and the angle
\beq\label{mj208} 
\alpha_{0,ab} = \arctan \frac{z_{b}^G-z_{a}^G}{x_{b}^G-x_{a}^G}
\eeq 
that indicates how the undeformed beam axis deviates (clockwise) from
the global axis $x^G$.
Strictly speaking, formula (\ref{mj208}) gives the correct result
only if $x_{Gb}>x_{Ga}$ and the rule for evaluation
of $\alpha_{0,ab}$ would need to be split into several cases
if the whole range had to be covered. However, we will
not really use the angle $\alpha_{0,ab}$ as such but rather
its sine and cosine, which are conveniently expressed as
\bea 
\cos\alpha_{0,ab} &=& \frac{x_{b}^G-x_{a}^G}{L_{ab}}
\\
\sin\alpha_{0,ab} &=& \frac{z_{b}^G-z_{a}^G}{L_{ab}}
\eea 

In the deformed configuration,
the auxiliary coordinate system $xz$ is rotated with respect
to the global system $x^Gz^G$ clockwise by angle 
\beq 
\alpha_{ab} = \alpha_{0,ab}-\varphi_a^G
\eeq 
We can imagine that, during phase A, the beam first moves as a rigid body
such that it gets translated by $u_{a}^G$ and $w_{a}^G$
and then rotated about the left end by $\varphi_a^G$ counterclockwise. Only then, 
during phase B, the 
right end is moved to its actual position in the
deformed configuration and the right
end section is rotated by $\varphi_b=\varphi_b^G-\varphi_a^G$. 
The displacements of the right end experienced
during phase B and expressed with respect to
the auxiliary axes are
\bea \label{mj215}
u_b &=& (u_{b}^G-u_{a}^G)\cos\alpha_{ab} + (w_{b}^G-w_{a}^G)\sin\alpha_{ab}+L_{ab}(\cos\varphi_a^G-1)\\
 w_b &=& -(u_{b}^G-u_{a}^G)\sin\alpha_{ab} + (w_{b}^G-w_{a}^G)\cos\alpha_{ab}+L_{ab}\sin\varphi_a^G
\eea 
and the rotation is
\beq \label{mj217}
\varphi_b = \varphi_b^G-\varphi_a^G
\eeq 
Therefore, if the global displacements are prescribed,
the local displacements with respect to the auxiliary
coordinate system can be evaluated---they represent
components of column matrix ${\boldsymbol{u}}_b$.
The corresponding column matrix ${\boldsymbol{f}}_{ab}$,
formally evaluated as ${\boldsymbol{g}}^{-1}({\boldsymbol{u}}_b)$,
has components ${X}_{ab}$, ${Z}_{ab}$ and
${M}_{ab}$. Here, ${M}_{ab}$ is directly
the end moment acting at the left end, while the end forces must be transformed to the global coordinate
system, which leads to
\bea 
X_{ab}^G &=& {X}_{ab}\cos{\alpha}_{ab} - {Z}_{ab}\sin{\alpha}_{ab} \\
Z_{ab}^G &=& {X}_{ab}\sin{\alpha}_{ab} + {Z}_{ab}\cos{\alpha}_{ab} 
\label{mj219}
\eea 
Finally, the forces at the right end,
\bea \label{e142}
X_{ba}^G &=& -X_{ab}^G \\
\label{e143}
Z_{ba}^G &=& -Z_{ab}^G
\eea 
are easily obtained from equilibrium,
and the moment at the right end is
\beq \label{e348}
{M}_{ba} = -{M}_{ab} + {X}_{ab}w_b - {Z}_{ab}(L_{ab}+u_b)
\eeq 

It is convenient to rewrite transformation rules (\ref{mj215})--(\ref{mj219}) in the matrix notation as
\bea\label{e145}
{\boldsymbol{u}}_b &=& \boldsymbol{T}(\varphi_a^G)\,(\boldsymbol{u}_b^G-\boldsymbol{u}_a^G) + \boldsymbol{l}(\varphi_a^G)\\
\label{e146}
\boldsymbol{f}_{ab}^G &=& \boldsymbol{T}^T(\varphi_a^G)\,{\boldsymbol{f}}_{ab}
\eea
where
\beq 
\boldsymbol{T}(\varphi_a^G) = \left(\begin{array}{ccc}
\cos(\alpha_{0,ab}-\varphi_a^G) & \sin(\alpha_{0,ab}-\varphi_a^G) & 0 \\
-\sin(\alpha_{0,ab}-\varphi_a^G) & \cos(\alpha_{0,ab}-\varphi_a^G) & 0 \\
0 & 0 & 1
\end{array}\right), 
\hskip 10mm
\boldsymbol{l}(\varphi_a^G) = L_{ab}\left(\begin{array}{ccc}
\cos\varphi_a^G-1 \\ \sin\varphi_a^G \\ 0 \end{array} \right)
\eeq 
Combining this with equation
\beq \label{e148}
{\boldsymbol{f}}_{ab} = {\boldsymbol{g}}^{-1}({\boldsymbol{u}}_b)
\eeq 
that formally describes the evaluation of the left-end forces $\boldsymbol{f}_{ab}$ by iterative solution of the set of nonlinear equations $\boldsymbol{g}(\boldsymbol{f}_{ab})=\boldsymbol{u}_b$,
we get
\beq 
\boldsymbol{f}_{ab}^G = \boldsymbol{T}^T(\varphi_a^G)\,{\boldsymbol{g}}^{-1}(\boldsymbol{T}(\varphi_a^G)(\boldsymbol{u}_b^G-\boldsymbol{u}_a^G) + \boldsymbol{l}(\varphi_a^G))
\eeq
This is the relation between the global components of joint displacements
and global components of end forces on beam $ab$.
To make it more readable, we rewrite it as
\beq \label{e150}
\boldsymbol{f}_{ab}^G = \boldsymbol{T}^T\,{\boldsymbol{g}}^{-1}(\boldsymbol{T}(\boldsymbol{u}_b^G-\boldsymbol{u}_a^G) + \boldsymbol{l})
\eeq
bearing in mind that matrices $\boldsymbol{T}$ and $\boldsymbol{l}$ depend on the left-end rotation, $\varphi_a^G$.

\subsection{Stiffness matrix}
\label{sec:stiff}

In the simplified notation, the differentiated form of equations (\ref{e145})--(\ref{e146}) reads
\bea
{\rm d}{\boldsymbol{u}}_b &=& \boldsymbol{T}\,({\rm d}\boldsymbol{u}_b^G-{\rm d}\boldsymbol{u}_a^G) + \left[\boldsymbol{T}'\,(\boldsymbol{u}_b^G-\boldsymbol{u}_a^G) +\boldsymbol{l}'\right]\,{\rm d}\varphi_a^G\\
{\rm d}\boldsymbol{f}_{ab}^G &=& \boldsymbol{T}^T\,{\rm d}{\boldsymbol{f}}_{ab} + \boldsymbol{T}'^T{\boldsymbol{f}}_{ab}\,{\rm d}\varphi_a^G
\eea
where 
\bea
\boldsymbol{T}'(\varphi_a^G) &=&\frac{\partial\boldsymbol{T}(\varphi_a^G)}{\partial\varphi_a^G} = \left(\begin{array}{ccc}
\sin(\alpha_{0,ab}-\varphi_a^G) & -\cos(\alpha_{0,ab}-\varphi_a^G) & 0 \\
\cos(\alpha_{0,ab}-\varphi_a^G) & \sin(\alpha_{0,ab}-\varphi_a^G) & 0 \\
0 & 0 & 0
\end{array}\right)
\\
\boldsymbol{l}'(\varphi_a^G) &=&\frac{\partial\boldsymbol{l}(\varphi_a^G)}{\partial\varphi_a^G} = L_{ab}\left(\begin{array}{ccc}
-\sin\varphi_a^G \\ \cos\varphi_a^G \\ 0 \end{array} \right)
\eea
Combining this with the differentiated form of (\ref{e148}),
\beq 
{\rm d}{\boldsymbol{f}}_{ab} = {\boldsymbol{G}}^{-1}\,{\rm d}{\boldsymbol{u}}_b
\eeq 
we get
\bea\nonumber
{\rm d}\boldsymbol{f}_{ab}^G  &=& \boldsymbol{T}^T{\boldsymbol{G}}^{-1}\left[\boldsymbol{T}({\rm d}\boldsymbol{u}_b^G-{\rm d}\boldsymbol{u}_a^G) + \left[\boldsymbol{T}'(\boldsymbol{u}_b^G-\boldsymbol{u}_a^G) +\boldsymbol{l}'\right]\,{\rm d}\varphi_a^G\right] + \boldsymbol{T}'^T{\boldsymbol{f}}_{ab}\,{\rm d}\varphi_a^G = \\
&=& \boldsymbol{T}^T{\boldsymbol{G}}^{-1}\boldsymbol{T}({\rm d}\boldsymbol{u}_b^G-{\rm d}\boldsymbol{u}_a^G)+
\left[\boldsymbol{T}^T{\boldsymbol{G}}^{-1}\left[\boldsymbol{T}'(\boldsymbol{u}_b^G-\boldsymbol{u}_a^G) +\boldsymbol{l}'\right]+ \boldsymbol{T}'^T{\boldsymbol{f}}_{ab}\right]\,{\rm d}\varphi_a^G
\label{e156}
\eea
which is the differentiated form of
(\ref{e150}).

Based on (\ref{e156}), we can set up the first three rows
of the element tangent stiffness matrix (in global coordinates).
The fourth row is minus the first row,
and the fifth row is minus the second row, because of relations (\ref{e142})--(\ref{e143}).
The sixth row is a bit more difficult to compute,
one needs to differentiate the expression for the right-end moment, ${M}_{ba}$.
From the moment equilibrium condition written with respect to the
centroid of the right end section in the deformed state,
we get\footnote{Equation (\ref{eMba}) is equivalent with (\ref{e348}),
just written here in terms of the global components.}
\beq\label{eMba} 
{M}_{ba} = -{M}_{ab} + X_{ab}^G(L_{ab}\sin\alpha_{0,ab}+w_b^G-w_a^G) - Z_{ab}^G(L_{ab}\cos\alpha_{0,ab}+u_b^G-u_a^G)
\eeq 
and the infinitesimal increment 
can be expressed as
\bea\nonumber
{\rm d}{M}_{ba} &=& -{\rm d}{M}_{ab} + (L_{ab}\sin\alpha_{0,ab}+w_b^G-w_a^G)\,{\rm d}X_{ab}^G - (L_{ab}\cos\alpha_{0,ab}+u_b^G-u_a^G)\,{\rm d}Z_{ab}^G +
\\ &&+ X_{ab}^G\,({\rm d}w_b^G-{\rm d}w_a^G) - Z_{ab}^G\,({\rm d}u_b^G-{\rm d}u_a^G) 
\eea
Consequently, the sixth row can be constructed as a linear
combination of the first, second and third row with coefficients
$L_{ab}\sin\alpha_{0,ab}+w_b^G-w_a^G$, $-L_{ab}\cos\alpha_{0,ab}-u_b^G+u_a^G$ and $-1$, resp.,
added to the row  $(Z_{ab}^G,-X_{ab}^G,0,-Z_{ab}^G,X_{ab}^G,0)$.
However, this does not even have to be done, since we know that the
stiffness matrix must be symmetric and we already know its sixth
column, except for the last (i.e., diagonal) entry.
So it is sufficient to copy the entries from the sixth column into
the sixth row and put 
\beq 
k_{66} = (L_{ab}\sin\alpha_{0,ab}+w_b^G-w_a^G)\,k_{16} - (L_{ab}\cos\alpha_{0,ab}+u_b^G-u_a^G)\,k_{26} - k_{36}
\eeq 
on the diagonal.

\section{Numerical examples}
\label{sec:numexamples}

A nonlinear beam element based on the proposed approach has been implemented into OOFEM \cite{patzak2001, patzak2012}, an object-oriented finite element code. To verify the implementation
and demonstrate the potential of the suggested
approach, several problems involving beams and frames
will be solved. 

\subsection{Pure bending of a cantilever beam}
\label{sec:3.1}

\begin{figure}[h!]
    \centering
    \begin{tabular}{cc}
    (a) & (b) 
    \\
    \includegraphics[width=0.35 \linewidth]{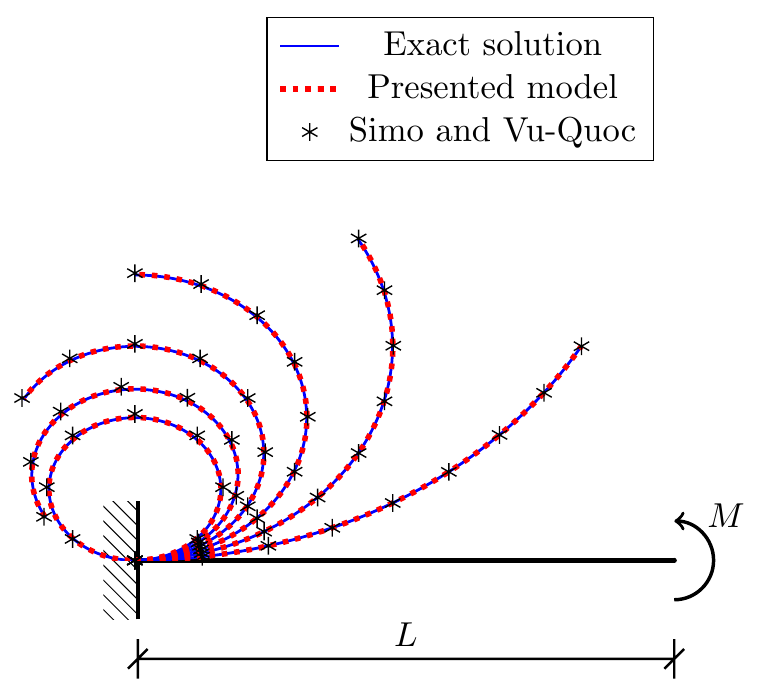}
    &  \includegraphics[width=0.35 \linewidth]{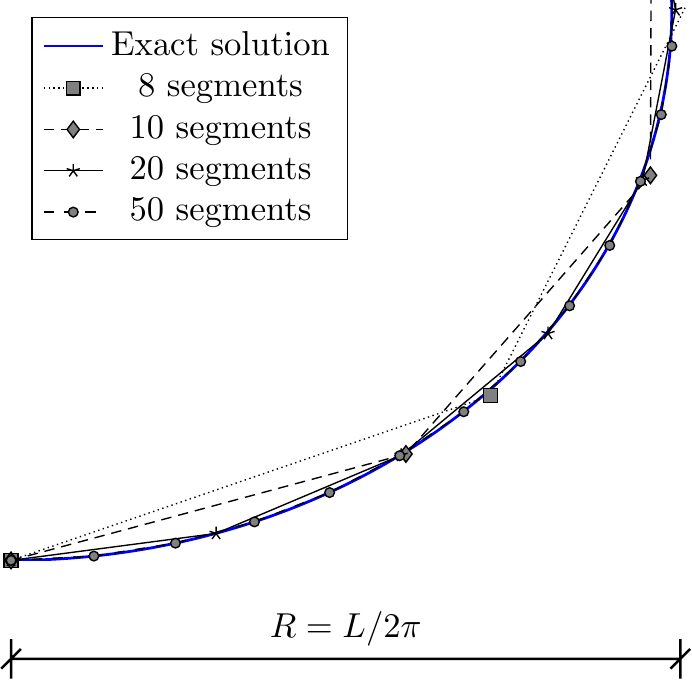}
    \end{tabular}
\caption{Pure bending of a cantilever beam subjected to end moment: (a) deformed shapes of the cantilever beam subjected to tip moment $M= 2 \pi EI/L$ obtained in six load steps, (b) close-up view of the solution at the end of the sixth step using 8, 10, 20, and 50 integration segments.}
\label{figex1}
\end{figure}

The first test, serving as a benchmark, deals with a cantilever of length $L$ and bending stiffness $EI$ loaded by a concentrated end moment $M$ on its right-end. The exact solution to this problem is a circular arc with radius $R=EI/M$. To deform the rod into a full closed circle, an end moment $M=2\pi EI/L$ needs to be applied. In this example, the loading is increased in six load steps, making the rod wind around itself at the end of the sixth step. The deformed shape of the beam at the end of each step is depicted in Fig.~\ref{figex1}a. The solution of the present model is compared with the one obtained by employing the geometrically exact finite beam element by Simo and Vu-Quoc \cite{simo1986} with a mesh of eight elements. The exact solution is reported as well. The overall agreement is good, and a detailed inspection reveals that the
simulation based on the present model, which uses
only one  two-noded element (i.e., only 3 global unknowns), is closer to the analytical solution. 

The example demonstrates that the present model
allows for a dramatic reduction of the number of global degrees of freedom, but
of course the number of segments $N$
used for numerical integration of the governing
equations (\ref{e212})--(\ref{e214}) must be chosen high enough to provide a good approximation. The results presented 
graphically in Fig.~\ref{figex1}a have been obtained using 100 segments. A close-up view of a part of the sixth step circle is showed in Fig.~\ref{figex1}b for calculations in which 8, 10, 20, and 50 numerical segments are employed. To ease the interpretation of the results, we connect the displaced grid points by straight segments, even though the curvature is constant along the beam and one could easily construct a more realistic visual representation. 
%Visually, it turns out that 20 segments provide an acceptable solution, which is further improved and gets close to the analytical one if 50 segments are used. It is also noted that even with just 10 segments the method provides an accurate estimation of displacements at the integration points. 
In contrast to  standard finite elements, for which shape functions allow interpolation of the displacement field on the basis of nodal values, here the displacement field is uniquely defined exclusively at the grid points. The values
at those points are sufficiently accurate even for a coarse grid.
%but it can be not sufficiently representative for the whole structure.\\

Considering that the exact ratio between the normalized moment $ML/EI$ and the the dimensionless curvature $L/R$ is unitary, we have calculated the dimensionless ratio $MR/EI$ based on the radius of curvature at the mid-span of the beam ($L=L/2$) and its relative error with respect to the exact solution. The results for the state at the end of the sixth load step are reported in Table~\ref{tab1}. They
illustrate how the integration grid refinement reduces the  error. When a traditional finite
element simulation with 8 elements is replaced by
the present method with 8 integration segments located within one single finite element,
the accuracy remains the same. 
The error is proportional to the square of the 
integration grid spacing,
and high accuracy can be achieved without changing
the number of the global degrees of freedom. 

\begin{center}
 \begin{table}[htb]
      \caption{Pure bending of a cantilever beam subjected to end moment: evaluation of errors caused by numerical integration along the beam element}
     \label{tab1}
     \centering
          \begin{tabular}{rllll}
\toprule
        \textbf{ Model} & \textbf{$MR/EI$} & \textbf{error $[\%]$}  \\
\midrule
         Exact & 1 & -\\
         Simo and Vu-Quoc \cite{simo1986} & 1.0262 & 2.617 \\
         8 segments & 1.0262 & 2.617 \\
         10 segments & 1.0166 & 1.664 \\
         20 segments & 1.0041 & 0.412 \\
         40 segments & 1.0010 & 0.103 \\
         80 segments &  1.0003 & 0.026\\
         \bottomrule
     \end{tabular}
 \end{table}
 \end{center}
 
\subsection{Williams toggle}
\label{sec:3.2}

\begin{figure}[h!]
    \centering
    \begin{tabular}{cc}
    (a) & (b) 
    \\
\includegraphics[width=0.4 \linewidth]{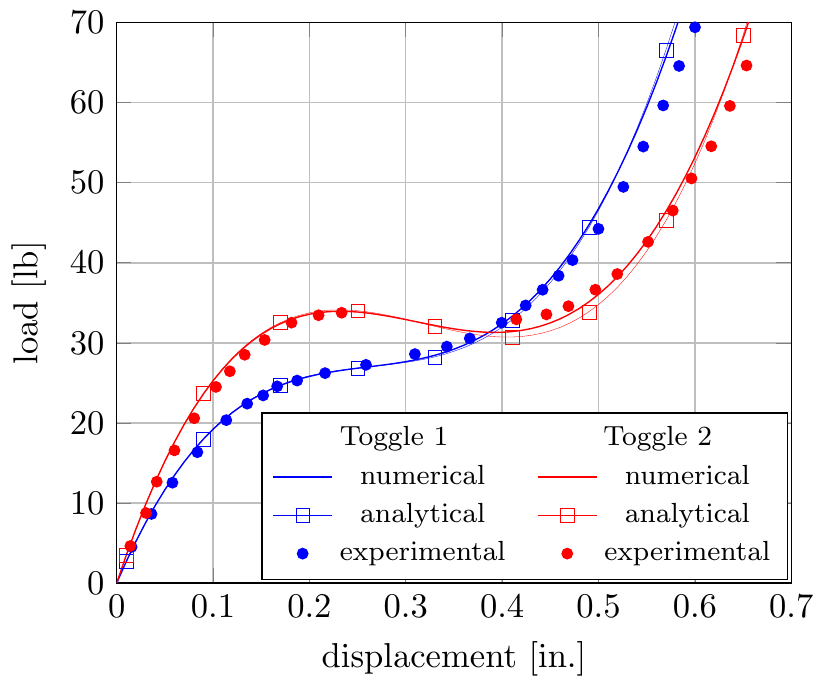}
    &  \includegraphics[width=0.22 \linewidth]{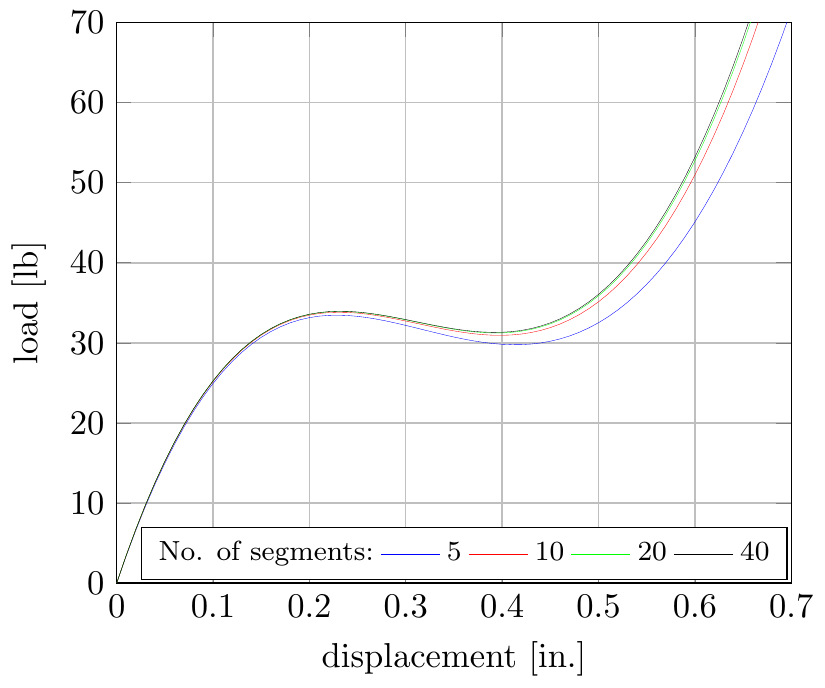}
    \end{tabular}
\caption{(a) Williams toggle and (b) its computational model that makes use of symmetry.} 
\label{f:will}
\end{figure}

Another relatively simple yet much more interesting problem is the so-called Williams toggle, for which Williams \cite{Williams64} 
provided experimental data as well as an approximate analytical solution. Physically, the toggle consists
of two symmetrically placed and rigidly connected straight beams whose axes
slightly deviate from the horizontal direction, see Fig.~\ref{f:will}a. The small angle between the beam axis and
the horizontal direction is denoted as $\psi$ and the
initial length of each beam as $L$. The toggle is loaded by a vertical force $P$. 

Owing to symmetry (including expected symmetry of the solution), it is sufficient to model the toggle by
a single element clamped at one end and vertically sliding
at the other end, with zero rotation and zero horizontal displacement, see Fig.~\ref{f:will}b. 
The resulting model has only one degree of freedom---the 
vertical displacement $w_2^G$. In fact, if the load control
is replaced by direct displacement control, which is perfectly legitimate here, the model has no global unknowns
and the equilibrium diagram can be constructed simply by
evaluating the end forces for a series of prescribed displacements at the right end. In the local coordinate
system of the beam, the displacement components are
$u_2=-w_2^G\sin\psi$ and $w_2=w_2^G\cos\psi$. Once
the end forces are computed, the applied force  $P=-X_{21}\sin\psi+Z_{21}\cos\psi$ is readily evaluated.

\begin{figure}[h!]
    \centering
    \begin{tabular}{cc}
    (a) & (b) 
    \\
\includegraphics[width=0.49 \linewidth]{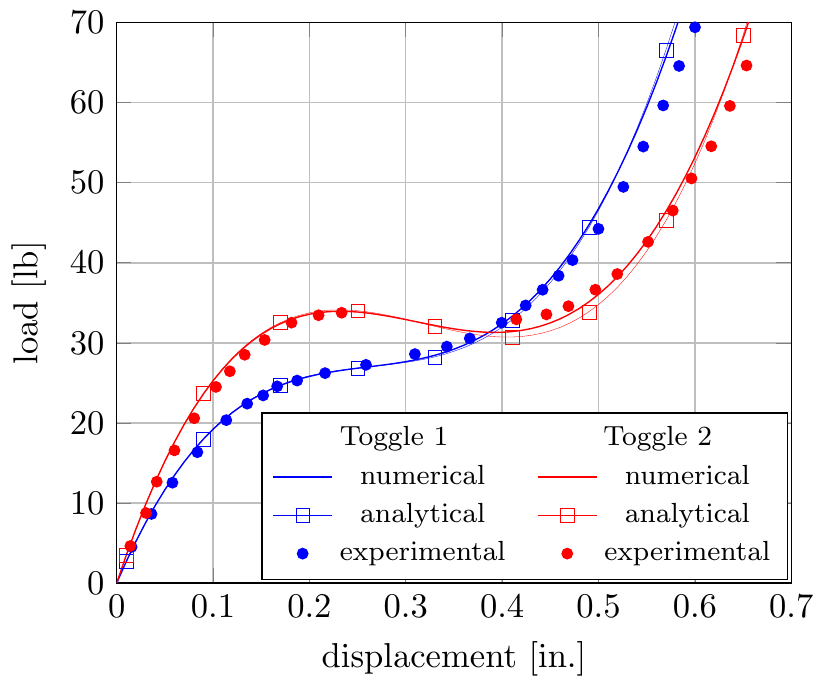}
    &  \includegraphics[width=0.49 \linewidth]{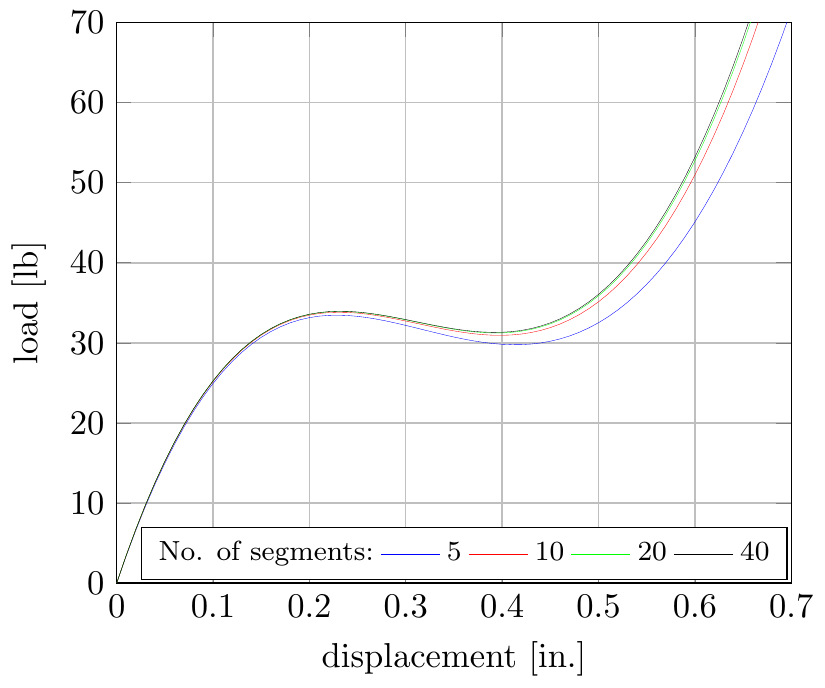}
    \end{tabular}
\caption{Williams toggle: (a) comparison of the present numerical results with the experimental data and approximate analytical solution by Williams, (b) effect of the number of segments used for numerical integration (for toggle 2).}
\label{f:will2}
\end{figure}

Williams tested two toggles with members of length $L=12.94$~in., made of aluminum alloy strips 
characterized by sectional stiffnesses  
$EA=1.885\cdot 10^6$~lb and $EI=9.27\cdot 10^3$~lb$\cdot$in.$^2$.
%The corresponding dimensionless parameter $EAL^2/EI$ is equal to 34,049.
The experimental results were reported for two geometries,
one with $\psi=0.0247$ and the other with $\psi=0.02985$.
The first case gives a monotonic load-displacement curve
while the second case leads to the snap-through behavior: 
the load-displacement curve exhibits a local maximum
followed by a local minimum, between which the equilibrium
state would be unstable under load control (but remains 
stable under displacement control). 
The experimental data are represented by individual points
(filled circles) in Fig.~\ref{f:will2}a while the 
approximate analytical solution derived by Williams
is shown as the continuous curve with hollow square symbols and the results of our numerical simulation as the continuous curves with no symbols. Blue color
refers to the first case ($\psi=0.0247$) and red color
to the second case ($\psi=0.02985$). Since Williams
performed his tests under load control, the descending
branch of the load-displacement diagram could not be
measured. Taking into account
that the measured values must be quite sensitive to
small changes in the initial geometry, the overall agreement between experimental and numerical results can be considered as very good. The simplified analytical solution derived by Williams is visually indiscernible from the present numerical
solution, except for a limited range of displacements between 0.4 and  0.6 in.\ in the second case  ($\psi=0.02985$, red curves in Fig.~\ref{f:will2}a).  In this range, the numerical solution is closer
to experimental results than the simplified analytical one.

The numerical results plotted in Fig.~\ref{f:will2}a
have been obtained with 40 integration segments,
to ensure high accuracy. The effect of the number of segments
is demonstrated in Fig.~\ref{f:will2}b. Already for 10 segments,
the numerical error is comparable to the experimental one,
and for 20 segments the complete computed curve is almost indiscernible from the curve obtained with 40 segments. 
In general, the errors are very small in the initial range up to the 
snap-through point (local maximum of the load-displacement curve), even for a simulation with just 5 integration segments.

\subsection{Buckling}
\label{sec:buckling}

The proposed beam element can efficiently handle highly nonlinear response, including potential loss of stability. Let us show a simple example that illustrates how instability phenomena can be treated.

Same as in Section~\ref{sec:3.1}, the example deals with a cantilever, but this time loaded by a concentrated force that induces axial compression. For a cantilever of length $L$, the buckling length is $L_b=2L$ and the corresponding Euler
critical load is evaluated using the well-known formula
\beq\label{eq:Pcr} 
P_{E} = \frac{EI\pi^2}{L_b^2} = \frac{EI\pi^2}{4L^2}
\eeq 
However, the derivation of this classical formula is based on the assumption of axial incompressibility. The adjusted derivation valid for axially compressible columns
is presented in detail in Appendix~\ref{appC}, and the resulting generalized version of formula (\ref{eq:Pcr}) is shown to be
\beq\label{eq:Pcr2} 
P_{cr} = \frac{EA}{2}\left(1-\sqrt{1-\frac{4EI\pi^2}{EAL_b^2}}\right) \approx
\frac{EI\pi^2}{L_b^2}\left(1+\frac{EI\pi^2}{EAL_b^2}\right)
\eeq 
The approximation is valid if $P_{cr}\ll EA$, which is always
the case here.

Numerically, the axially compressed cantilever can be described by a single element connecting two nodes.
Node 2 is fixed and the displacements and rotation 
of node 1 play the role of global degrees of freedom. 
If the beam is perfectly straight and the applied force $P$ is perfectly aligned with the beam axis, the numerically computed solution corresponds to axial compression and degrees of freedom $w_1$ and $\varphi_1$ remain equal to zero. The beam is uniformly compressed and, since we use here a model based on Biot strain, displacement $u_1$ is proportional to the applied force $P$. 
Of course, this
type of solution becomes unstable if the applied force
exceeds the critical one. 

The loss of stability can be
detected by checking the eigenvalues of the tangent structural stiffness matrix. Initially, all eigenvalues
are positive, which indicates that the tangent stiffness matrix is positive definite and the solution of the equilibrium equations corresponds to a minimum
of potential energy, i.e., to a stable state. 
Stability is lost when at least one eigenvalue becomes negative, and the onset of instability is characterized by
the smallest eigenvalue equal to zero.

\begin{figure}[h!]
    \centering
    \begin{tabular}{cc}
    (a) & (b) 
    \\
\includegraphics[width=0.49 \linewidth]{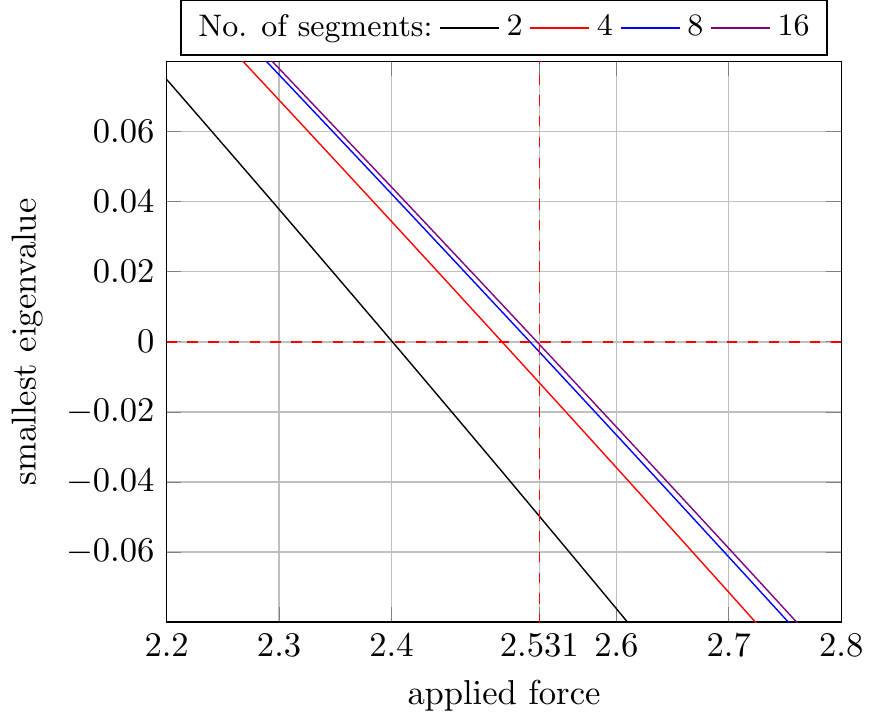}
    &  \includegraphics[width=0.49 \linewidth]{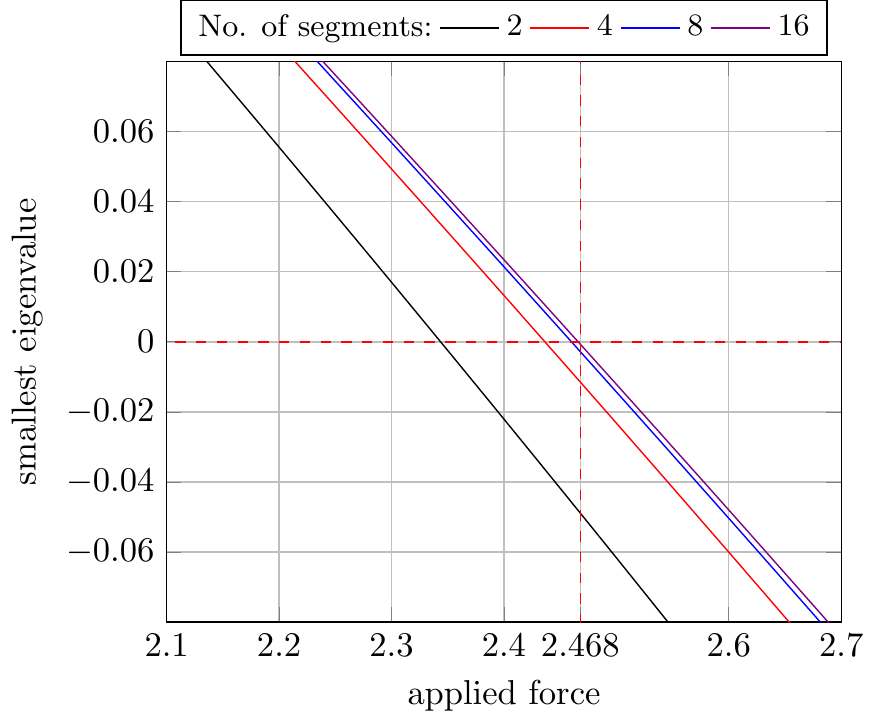}
    \end{tabular}
\caption{Smallest eigenvalue of the tangent stiffness matrix as function of the applied force: (a) $EA~=$ 100, (b)~$EA~=$~10,000.} 
\label{f:eigenval}
\end{figure}

The dependence of the smallest eigenvalue of the tangent stiffness matrix on the applied force is
plotted in Fig.~\ref{f:eigenval}. The problem is treated in the
dimensionless form---the beam length $L$ and the flexural
stiffness $EI$ are set to 1, which means that the
dimensionless value of the applied force in fact
corresponds to $PL^2/EI$ and the computed displacement
to $u_1/L$. 
The behavior of the model is affected by the 
axial sectional stiffness $EA$, which corresponds to the
dimensionless slenderness parameter $EAL^2/EI$.
This parameter is the square of the ratio
$L/i$ where $i=\sqrt{I/A}$ is the sectional radius of inertia. For instance, for the strip used by Williams
and described in Section~\ref{sec:3.2}, the span-to-depth
ratio is $L/h\approx 53$, which certainly represents
an extremely slender beam, and parameter $EAL^2/EI$
is in this case approximately equal to 34,000. In our simulations, we will
typically consider $EAL^2/EI=$10,000 or 100, the latter choice
representing a rather deep beam.

\begin{figure}[h!]
\centering
\includegraphics{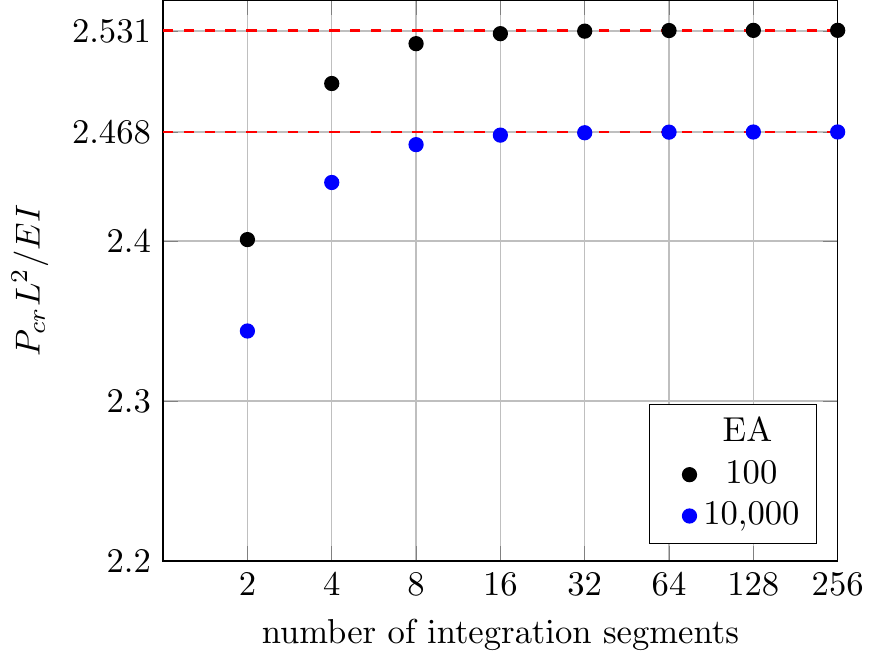}
\caption{Dimensionless critical force evaluated numerically from the criterion of zero minimum eigenvalue, plotted as function of the number of integration segments.}
\label{f:eigenval2}
\end{figure}

The numerical solution naturally depends
on the number of segments used for integration
of the governing equations on the element level.
As seen in Fig.~\ref{f:eigenval2}, the critical force
evaluated from the condition of zero minimum eigenvalue quickly converges as
the number of segments increases, but the limit value
is affected by the slenderness parameter.

%To facilitate the interpretation, the applied force
%is normalized by the Euler buckling load $P_{E}$ given by formula (\ref{eq:Pcr}).
In the dimensionless format (i.e., for $L$ and $EI$ set to 1),
the Euler critical load is $P_E=\pi^2/4\approx 2.4674$.
For highly accurate numerical simulations (a sufficiently high 
number of integration segments and very short incremental steps,
at least in the vicinity of the critical state),
the onset of instability occurs at $P_{cr}\approx 2.531$
for $EA=100$ and  at $P_{cr}\approx 2.468$
for $EA=$10,000. This is correct, because 
Euler formula (\ref{eq:Pcr}) is exact for the ideal
case of an axially incompressible beam. 
The generalized formula (\ref{eq:Pcr2}) gives 2.531485
 for $EA=100$ and 2.46801 for $EA=10,000$
 if the ``exact'' expression is used,  
 in perfect agreement  with the
loads for which the onset of instability has been detected
by highly accurate numerical evaluation of the tangent
stiffness matrix and its minimum eigenvalue.
The approximate formula (i.e., the last expression on the right-hand side
of (\ref{eq:Pcr2})) gives 2.5283
 for $EA=100$ and 2.46801 for $EA=10,000$.
 In the former case (deep beam), the approximation induces a difference of about 0.13 \%
 compared to the exact formula, while in the latter
 case (slender beam), the first 6 valid digits
 of the resulting value remain the same.

For loads exceeding the critical one, the straight-beam solution becomes unstable and thus physically irrelevant, and it would be desirable to compute the bifurcated stable solution that describes the actual shape
of the buckling beam. A rigorous approach would be
to find the eigenvector associated with the zero eigenvalue
of the stiffness matrix at the onset of buckling and
then search for a branch of the equilibrium diagram that
bifurcates from the main one in the direction given by this eigenvector. Sophisticated
techniques of this kind have been proposed and developed in the literature. 

As an alternative, one can simply perturb
the original problem and change the equilibrium
diagram with a bifurcation point into an equilibrium
diagram which closely follows one of the bifurcated stable branches but does not pass through a critical point. This is typically achieved by breaking symmetry
of the original problem. In our case, we can consider, e.g., the load as slightly eccentric, or the
beam as slightly curved.

\begin{figure}[h]
\centering
\begin{tabular}{cc}
(a) & (b) \\
\includegraphics[width=0.5 \linewidth]{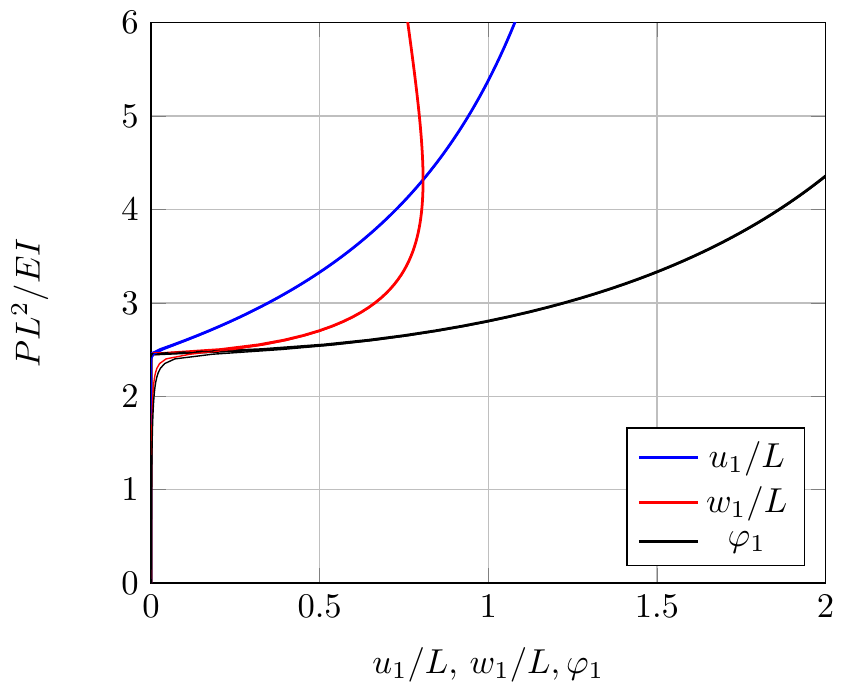}
&
\includegraphics[width=0.5 \linewidth]{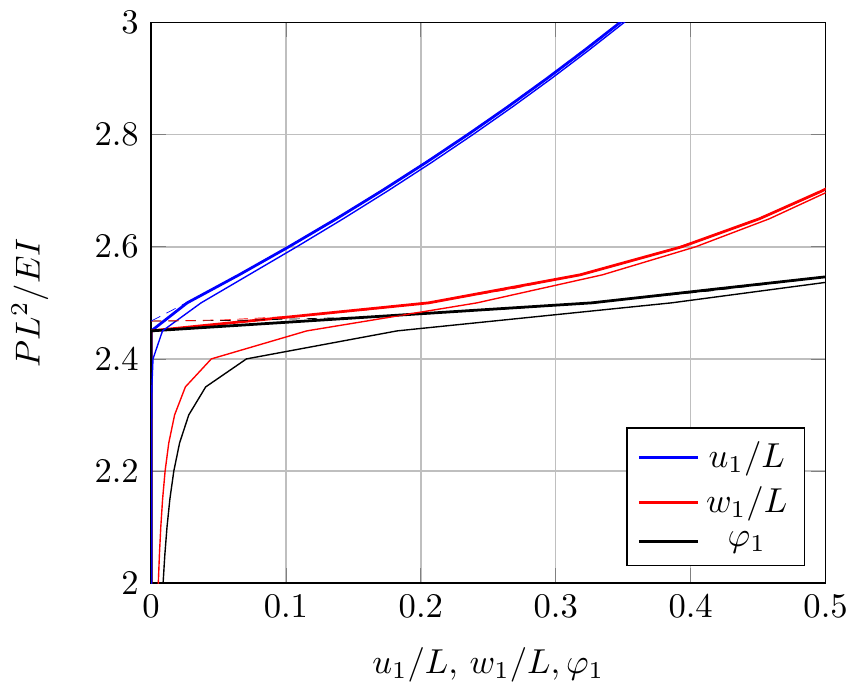}
\\
(c) \\
\includegraphics[width=0.5 \linewidth]{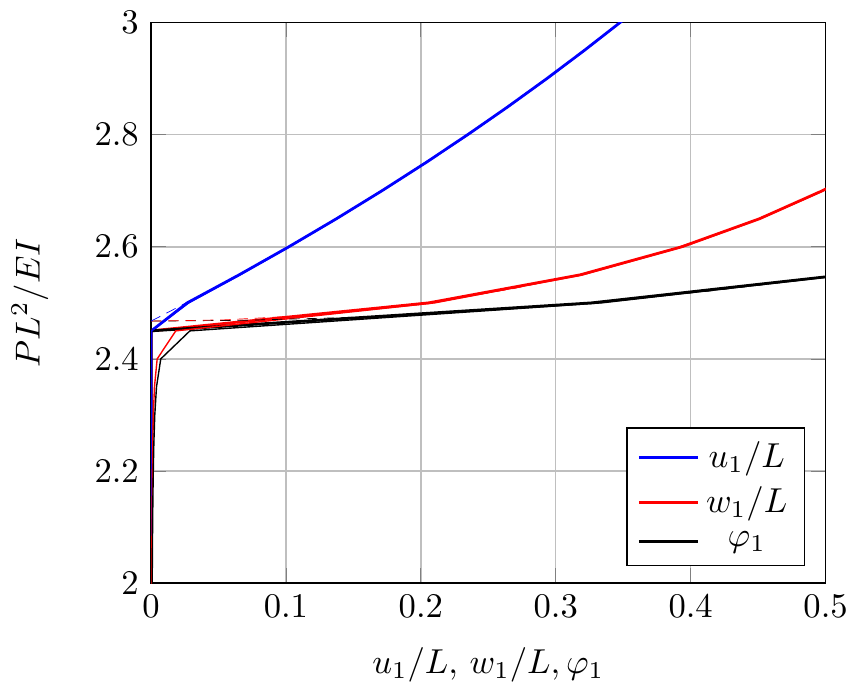}
\end{tabular}
\caption{Equilibrium diagrams for a compressed cantilever: (a) full range, (b)--(c) close-up around bifurcation point. Computed with $EAL^2/EI=$10,000 and step size $\Delta PL^2/EI=0.05$. Thick lines 
correspond to zero eccentricity, thin lines to eccentricity $e/L=10^{-3}$ in parts (a) and (b) and $e/L=10^{-4}$ in part (c), and dashed lines to the analytical solution in the axially incompressible case
($EA\to\infty$). 
%The load $P$ is normalized by $EI/L^2$ and the displacements $u_1$ and $w_1$ are normalized by $L$.
} 
\label{f:cantstab}
\end{figure}

Making use of the first option, we combine the
applied force $P$ with an applied moment $M=Pe$ 
where $e$ is a fixed small eccentricity. 
The obtained equilibrium diagrams are plotted in
Fig.~\ref{f:cantstab} by thin lines. The thick lines in
the same figure correspond to the original, unperturbed problem, i.e., to the load applied with zero eccentricity.
To follow the bifurcated stable branch instead of the 
main branch that becomes unstable for loads exceeding
the critical one, the equilibrium iteration after each increment of applied force is started from 
a perturbed trial state, obtained by solving an auxiliary equilibrium problem for loading by a small
applied moment added to the previously applied force. 
This moment is similar to the moment $Pe$ due to eccentricity but this time it is not considered as the 
actual part of applied loads---it is used  to 
generate a perturbed initial state for equilibrium 
iterations and then removed when the actual axial loading 
is increased. As a result, the final converged state 
corresponds to the original problem of an axially loaded
straight beam. If the load is below the critical level,
the iteration necessarily converges to the trivial solution
(i.e., the axially compressed beam remains straight), because this is the only solution of the
equilibrium equations. On the other hand, if the load is
above the critical level, there exist three equilibrium
states, one of which is unstable (straight beam) while
the other two are stable (buckling to one or the other side). An iterative process that starts from an unsymmetric
state is likely to end up on one of the two stable
bifurcated branches. This is indeed confirmed by numerical
simulations.

%for $EAL^2/EI=100$ and $e/L=10^{-3}$.

%A simpler, more heuristic approach 

The solid curves plotted in Fig.~\ref{f:cantstab} have been computed for
parameter $EAL^2/EI$ set to 10,000 using load increments
$\Delta P=0.05\,EI/L^2$. Up to $P=2.45\,EI/L^2$, the 
solution obtained when the load is axial and the iterations
start from a perturbed state remains on the main branch,
i.e., the lateral displacement $w_1$ and the rotation
$\varphi_1$ remain zero (up to the tolerated numerical
error) while the axial displacement $u_1$ increases
proportionally to the applied load (due to the high axial
stiffness, it also appears to be almost zero in the diagrams). On the other hand, the solutions obtained 
when the load is considered as eccentric gradually deviate
from the straight main branch. 

Fig.~\ref{f:cantstab}a
shows the full equilibrium diagrams for loads up to 
$P=6\,EI/L^2$. On this scale, the thick and thin solid curves almost coincide, except for the immediate vicinity of the bifurcation point. To better assess the difference, 
the diagrams in Fig.~\ref{f:cantstab}b,c are limited
to the range of normalized load $PL^2/EI$ between 2 and 3. 
For eccentricity $e=10^{-3}L$, the deviation from the
main branch becomes quite pronounced (Fig.~\ref{f:cantstab}b) while for a reduced eccentricity
$e=10^{-4}L$ it is much less important (Fig.~\ref{f:cantstab}c). 
For comparison, the dashed curves show the analytical solution derived for an axially incompressible beam in Appendix~\ref{sec:examplecantilever} and described by formulae (\ref{eq:ex298})--(\ref{eq:ex300}). 

\subsection{Frames}

\begin{figure}[h!]
    \centering
    \begin{tabular}{cc}
    (a) & (b) 
    \\
    \includegraphics[width=50mm]{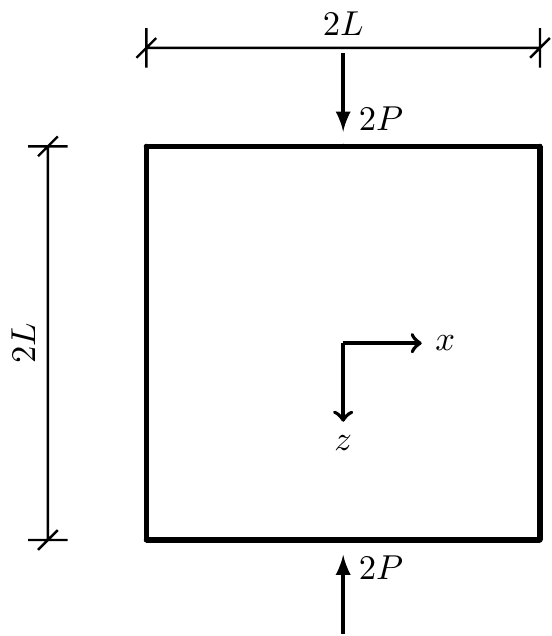}
    &  \includegraphics[width=40mm]{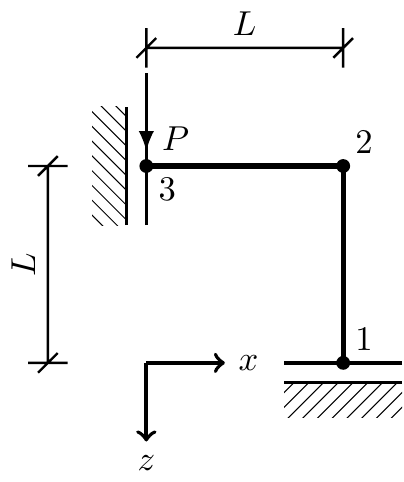}
    \end{tabular}
    \caption{Square frame: (a) Geometry and applied loads (case of compression),
    (b) one-quarter model that exploits symmetry, leading to five global unknowns.}
    \label{squareframe}
\end{figure}

As a more challenging example, let us consider two problems of large deflection of frames: a square frame loaded at the midpoints of a pair of opposite sides (Fig.~\ref{squareframe}a) and a square-diamond frame loaded at two vertices with hinges (Fig.~\ref{diamondframe}a).
Analytical solutions for the deflections and bending moments were presented
by Kerr \cite{kerr1964} for the square frame loaded at the midpoints of a pair of opposite sides and by Jenkins et al.~\cite{jenkins1966} for the diamond-shaped frame. The elliptic integrals were numerically evaluated and presented in a tabulated form by Mattiasson \cite{mattiasson1981} based on the procedure described by King \cite{king1924}, which has shown excellent convergence properties and highly accurate results. For this reason, other authors often consider Mattiasson's solutions as  analytical ones \cite{chorn2021,wood1977}.

Owing to symmetry, only a quarter of each frame needs to be analyzed. In our simulations, it is sufficient to use a mesh consisting of two elements for the square frame (Fig.~\ref{squareframe}b) and  a single-element mesh  for the diamond frame
(\ref{diamondframe}b). In the figures, the applied
force is oriented such that it induces compression, but the simulations cover the 
opposite orientation leading to tension, too. Mattiasson \cite{mattiasson1981} neglected  axial as well as shear deformations.
To be able to compare our numerical results
with his, we need  to set the
 axial stiffness  to a sufficiently large value.
 However, our simulations can also work with
 lower, more realistic values.
 
 The results are again presented in the dimensionless form, with all quantities normalized by suitable combinations of the flexural stiffness $EI$ and beam length $L$. 
 In numerical simulations, $EI$ and $L$ are set to unity and the input values of axial stiffness $EA$
 and applied force $P$ have the meaning of dimensionless parameters $EAL^2/EI$
 and $PL^2/EI$. The computed
 displacements
 then correspond to dimensionless fractions
 $u/L$ and $w/L$, and bending moments to $ML/EI$.
 
\begin{figure}[h!]
    \centering
    \begin{tabular}{cc}
    (a) & (b) 
    \\
    \includegraphics[width=75mm]{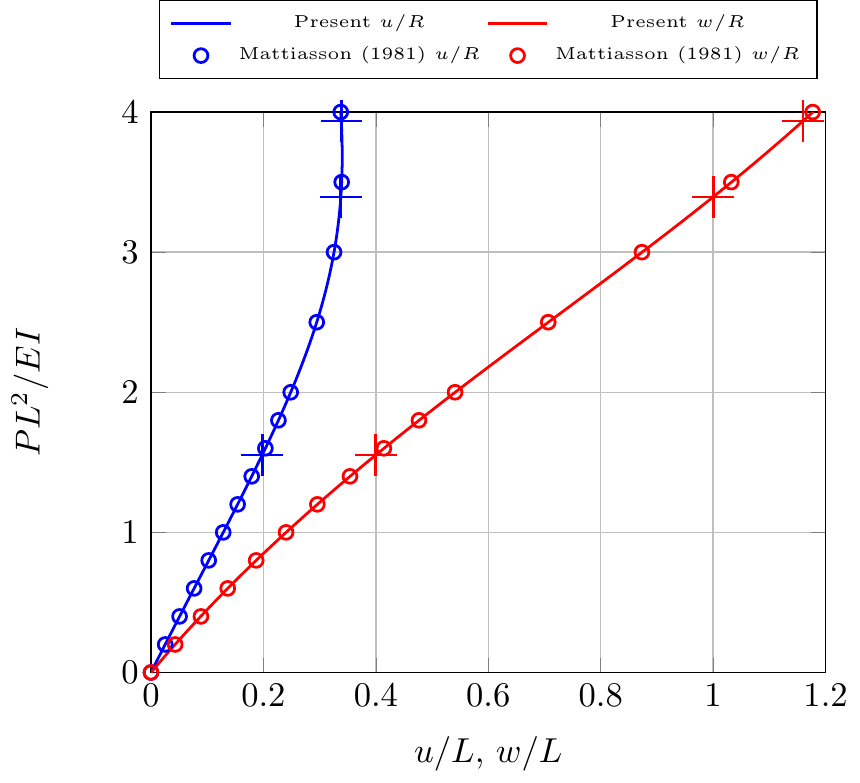}
    &  \includegraphics[width=75mm]{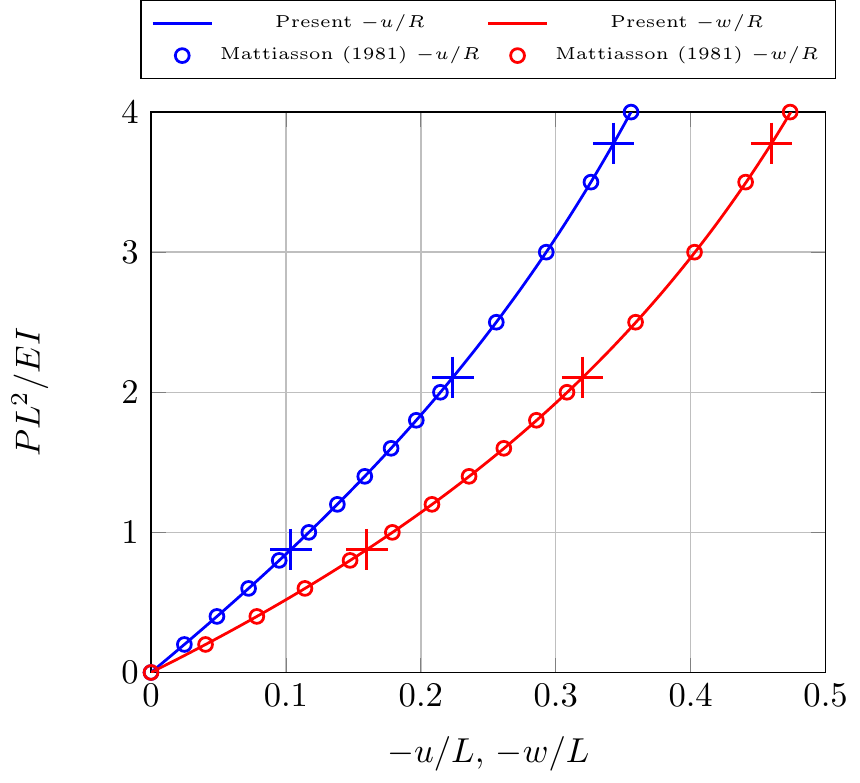}
    \end{tabular}
    \caption{Square frame: load-displacement curves for (a) compression loading, (b) tension loading. With reference to Fig.~\ref{squareframe}b, horizontal displacement $u$ is reported at node 1 while vertical displacement $w$ is reported at node 3. Dotted lines represent Mattiason's solution. Crosses indicate three states for which the deformed shapes are plotted in Fig.~\ref{deformedshape_SF}.}
    \label{load_dispSF}
\end{figure}

\begin{figure}[h!]
    \centering
    \begin{tabular}{cc}
    (a) & (b) 
    \\
    \includegraphics[width=60mm]{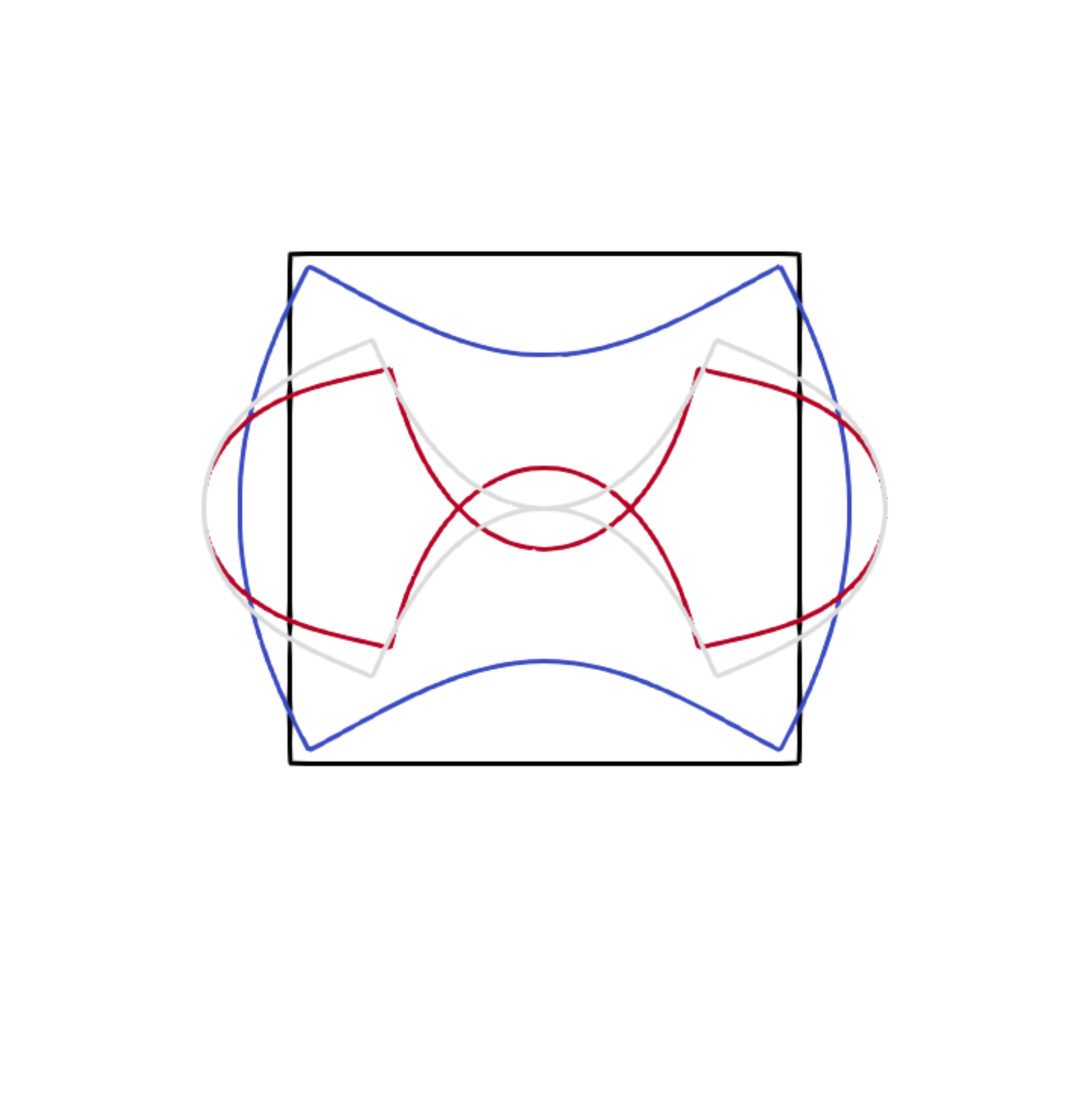}
    &  \includegraphics[width=60mm]{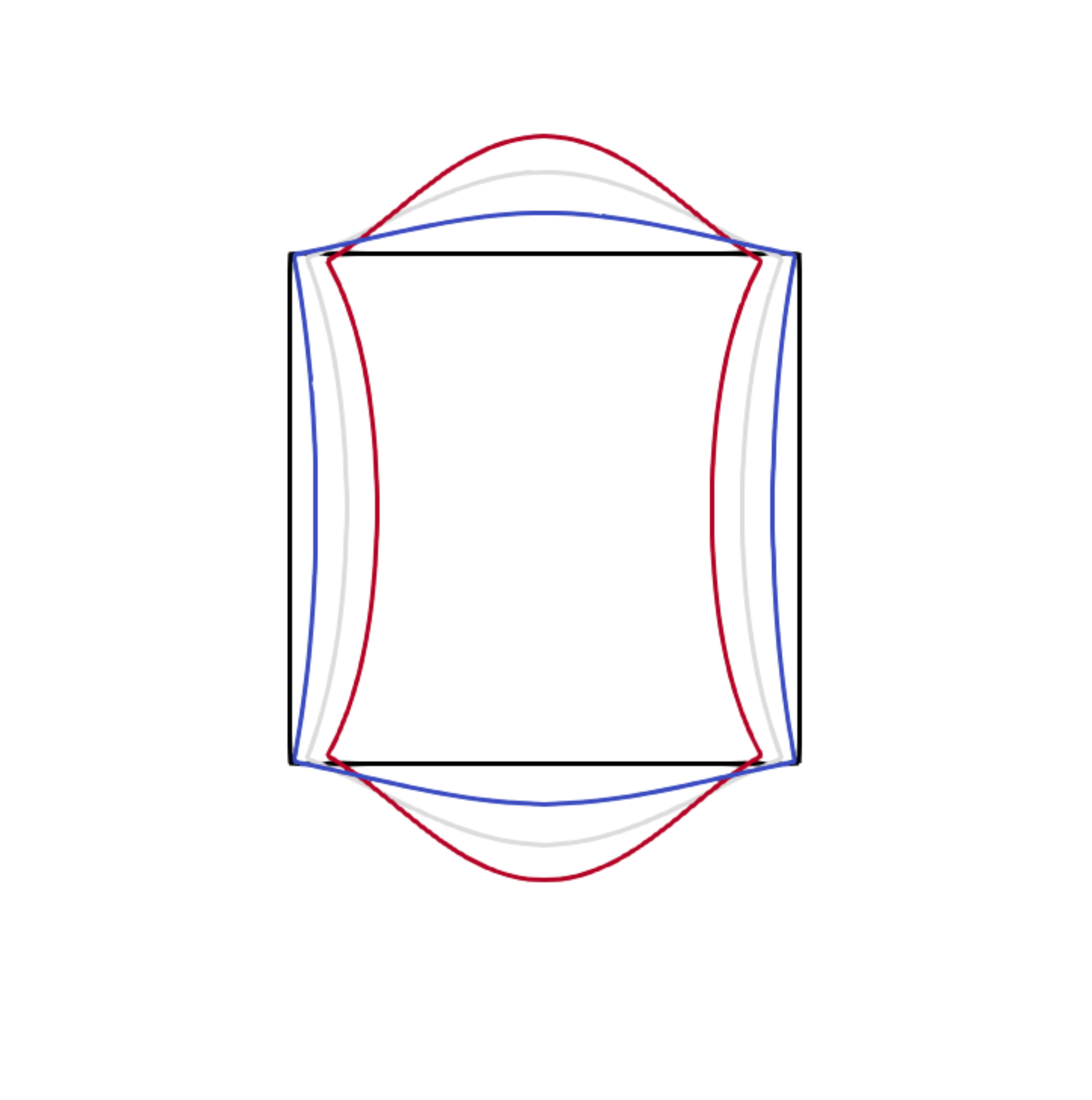}
    \end{tabular}
    \caption{Square frame: deformed shapes (scale factor equal to 1) for (a) compression loading (b) tensile loading. The corresponding points in the load-displacement diagrams are indicated in Fig. \ref{load_dispSF} by '+' symbols. }
    \label{deformedshape_SF}
\end{figure}

\begin{figure}[h!]
    \centering
    \begin{tabular}{cc}
    (a) & (b) 
    \\
    \includegraphics[width=50mm]{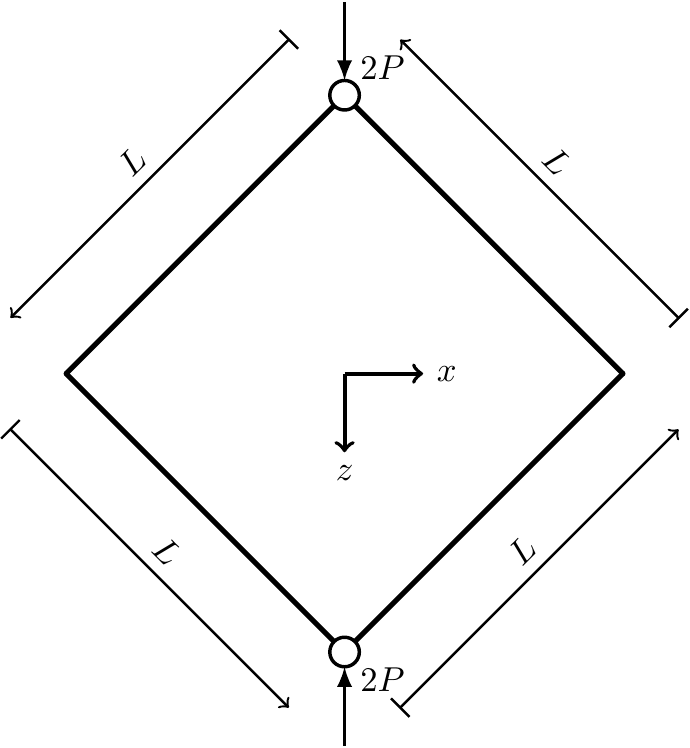}
    &  \includegraphics[width=40mm]{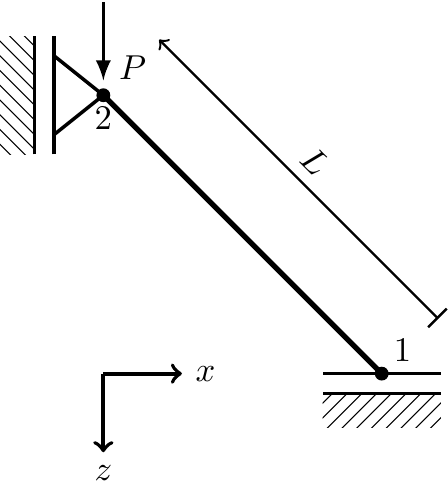}
    \end{tabular}
    \caption{Diamond frame: (a) Geometry and applied loads (case of compression),
    (b) one-quarter model that exploits symmetry, leading to three global unknowns.}
    \label{diamondframe}
\end{figure}

 Consider first the square frame shown in Fig.~\ref{squareframe}. Its response under compressive loading is characterized by the 
 load-displacement diagrams in Fig.~\ref{load_dispSF}a, with the red curve
 corresponding to the vertical deflection $w_3$
 and the blue curve to the horizontal displacement $u_1$ (both normalized by $L$). 
 Empty circles represent Mattiasson's data
 and the crosses indicate three states for 
 which the deformed shapes are plotted in Fig.~\ref{deformedshape_SF}a. Analogous results
 for the case of tensile loading are presented
 in Fig.~\ref{load_dispSF}b in terms of the load-displacement diagrams and in Fig.~\ref{deformedshape_SF}b in terms of the deformed shapes at three selected states.

The agreement of our results with
 Mattiasson's solution is seen to be excellent. For compression, the response after reaching the load level $PL^2/EI$=3.3942 (second cross in Fig.~\ref{load_dispSF}a) loses its physical meaning because of non-physical penetration of node 3 into its mirrored counterpart (see the deformed shape red colored in Fig.~\ref{deformedshape_SF}a). The present
 numerical results have been computed using
30 integration segments per element and with
 the axial stiffness parameter $EAL^2/EI$
 set to $10^6$. This value is fully sufficient to get very close to the inextensible limit.
 If the parameter is increased to $10^7$,
 the relative change of the vertical displacement
 at the end of the simulation (i.e., at load level 
 $PL^2/EI=4$)  is only $3\cdot 10^{-6}$.
 Even for $EAL^2/EI=10^4$, the relative change
 with respect to the inextensible case would be 
 $3\cdot 10^{-4}$, which is still negligible.
 On the other hand, $EAL^2/EI=10^2$ would lead
 to a relative change of 2.9 \% in the vertical displacement and 1.3 \% in the horizontal displacement, which may already play some role. 
 Still lower values of the axial stiffness parameter are not relevant because
 they correspond to deep beams for which 
 the beam theory with neglected shear distortion
 would be inappropriate. 
 
 The choice of the axial stiffness parameter reflects the geometry of the frame
 (shape of the cross section and span-to-depth ratio). Let us now explore the effect of 
 a numerical parameter---the number of integration segments per element. The ``high-precision''
 value (i.e., the value computed with an extremely high number 
 of integration segments) of normalized displacement $w_3/L$
 computed at load level $PL^2/EI=4$
 is 1.177368 for $EAL^2/EI=$10,000 and 1.211001
 for $EAL^2/EI=100$. The values obtained for
 various numbers of integration segments per element and the corresponding relative errors are summarized in Table~\ref{tab:2}.  Already for 8 segments, the discretization error is below 
  1 \%, and it decreases proportionally to the square of the segment size.
 
 \begin{table}[htb]
     \centering
    \caption{Square frame: evaluation of errors in displacement caused by numerical integration along the beam elements.}
     \begin{tabular}{rllll}
     \toprule
     & \multicolumn{2}{c}{$EAL^2/EI=$ 10,000} & \multicolumn{2}{c}{$EAL^2/EI=$ 100} \\
         \textbf{number of segments} & \textbf{displacement} & \textbf{error $[\%]$} & \textbf{displacement} & \textbf{error} $[\%]$  \\
         \midrule
         8 & 1.221452 & 0.863 & 1.187822 & 0.888\\
         10 & 1.217680 & 0.551 & 1.184051 & 0.568\\
         20 & 1.212668 & 0.138 & 1.179036 & 0.142\\
         40 & 1.211418 & 0.0344 &  1.177785 & 0.0354\\
         80 & 1.211105 & 0.0086 & 1.177472 & 0.0089\\
         160 & 1.211027 & 0.0022 & 1.177394 & 0.0022\\
         320 & 1.211008 & 0.0006 & 1.177374 & 0.0005\\
         %640 & 1.211003 &  & 1.177369 & \\
         $\to\infty$ & 1.211001 & 0 & 1.177368 & 0\\
         \bottomrule
     \end{tabular}
     \label{tab:2}
 \end{table}

The example of the square frame loaded in compression can further be  exploited for illustration of convergence properties of the proposed numerical method.
The objective is to explore how the computational procedure converges if a very large load increment is applied.
The model deals with five
global unknowns ($u_1$, $u_2$, $w_2$, $\varphi_2$, and $w_3$), which are found iteratively by Newton-Raphson equilibrium iterations. 
In each iteration, the end forces and the tangent stiffness need to be evaluated for given values
of the global unknowns, and this evaluation is
also performed iteratively, using the technique
described in Section~\ref{sec:endforces}. Each iterative process
uses a certain error tolerance, which can influence the
number of iterations needed to satisfy the underlying
equations with sufficient accuracy. 

It turns out that the numerical scheme is more robust for lower values of the axial stiffness parameter.
For $EAL^2/EI=100$ (considered as low), it is possible to apply the total load $PL^2/EI=4$ in one single step, starting from the undeformed configuration. 
The error (defined as the norm of the unbalanced forces normalized by the same factor $EI/L^2$ as the actual load) first increases
from 4.0 to 43.9,
but after 6 iterations it is below $10^{-3}$
and after 8 iterations below $10^{-9}$. 
For $EAL^2/EI=1000$, the maximum step size is $\Delta PL^2/EI=1.99$,
and for $EAL^2/EI=$10,000, it is $\Delta PL^2/EI=0.6$, which means that the whole diagram 
depicted in Fig.~\ref{load_dispSF}a can be covered respectively in 3 or 7 incremental steps. Of course, larger steps require more global
Newton-Raphson iterations. 

For comparison, Table~\ref{tab:1} shows the average
numbers of global iterations per step needed
to increase the load to $PL^2/EI=4$, depending
on the step size, axial stiffness parameter and
relative tolerance (maximum allowed norm of
unbalanced forces normalized by $EI/L^2$).
In each row, the step size is indicated in the first column and the other columns contain the
average numbers of iterations per step for various
combinations of parameters $(EAL^2/EI,\varepsilon_{\rm tol})$, in each case specified in the column heading. For sufficiently short steps, convergence is very regular. For instance,
for $EAL^2/EI=100$ and step size $\Delta PL^2/EI=0.25$, the whole curve is covered in 16 steps and, in each step, 2 iterations are sufficient to bring the error below $10^{-3}$ and
2 additional iterations bring the error below $10^{-9}$. On the other hand, for larger steps or higher axial stiffness, more iterations
are needed in the initial part of the iterative process, during which the evolution of error is
typically less regular. Once the computed
approximation gets close to the exact solution, quadratic convergence is observed 
and the error is easily reduced from $10^{-3}$
to $10^{-9}$ in at most two iterations.

\begin{table}[htb]
    \centering
        \caption{Square frame loaded by compression: average  numbers of global iterations per step, depending on $(EAL^2/EI,\varepsilon_{\rm tol})$ where $EAL^2/EI=100$ or 10,000 is the axial stiffness parameter and $\varepsilon_{\rm tol}=10^{-3}$ or $10^{-9}$ is the tolerance.}
    \begin{tabular}{lllll}
    \toprule
       $\mathbf{\Delta PL^2/EI}$  &  $\mathbf{(100}$\textbf{,} $\mathbf{10^{-3})}$ &  $\mathbf{(100}$\textbf{,} $\mathbf{10^{-9})}$ &  $\mathbf{(10000}$\textbf{,} $\mathbf{10^{-3})}$ &  $\mathbf{(1000}$\textbf{,} $\mathbf{10^{-9})}$\\
       \midrule
    4 & 6 & 8 \\ 
    2 & 4 & 5.5 \\
    1 & 4 & 5 \\
    0.5 & 3 & 4 & 5.125 & 6.75 \\
    0.25 & 2 & 4 & 4 & 5.438\\
    \bottomrule
    \end{tabular}
    \label{tab:1}
\end{table}

\begin{figure}[h!]
    \centering
    \begin{tabular}{cc}
    (a) & (b) 
    \\
    \includegraphics[width=75mm]{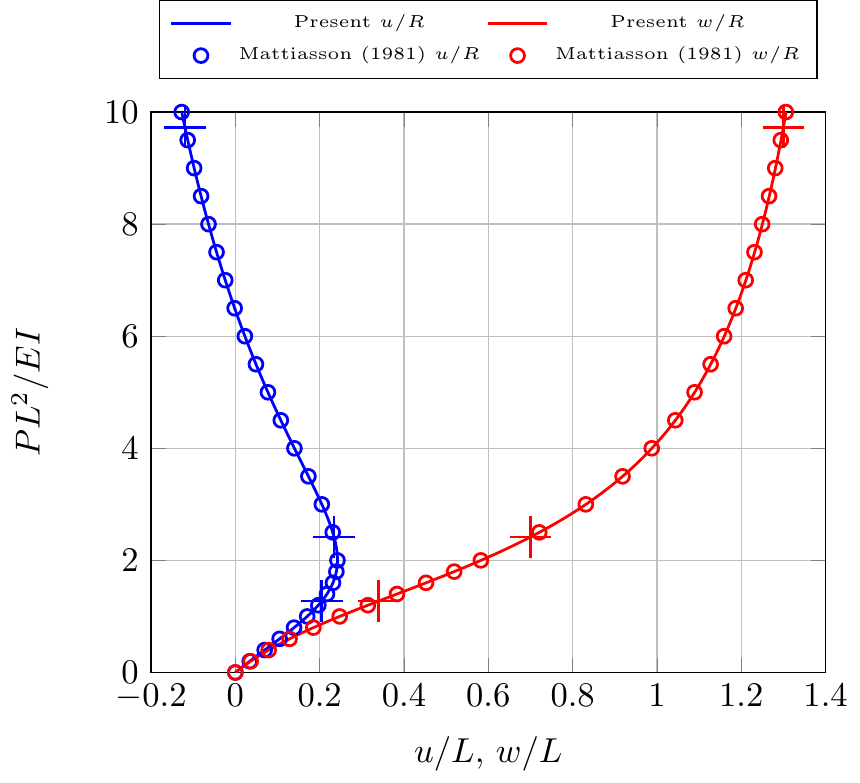}
    &  \includegraphics[width=75mm]{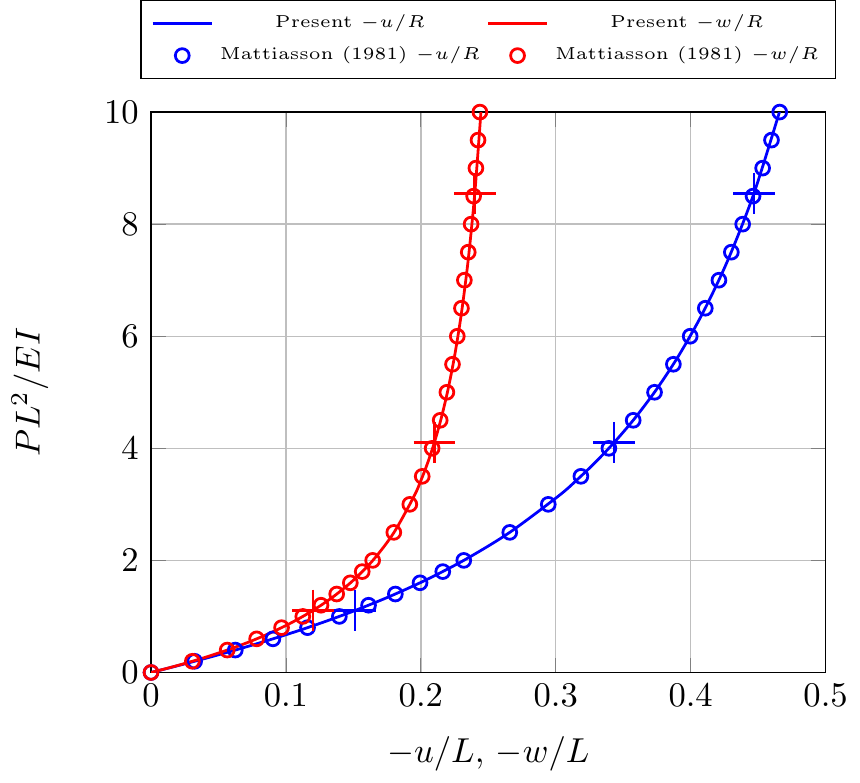}
    \end{tabular}
    \caption{Diamond frame: load-displacement curves for (a) compression loading, (b) tension loading. With reference to Fig. \ref{diamondframe}b, horizontal displacement $u$ is reported at node 1 while vertical displacement $w$ is reported at node 2. Dotted lines represent Mattiason's solution. Crosses indicate three states for which the deformed shapes are plotted in Fig.~\ref{deformedshape_DF}.}
    \label{load_dispDF}
\end{figure}

\begin{figure}[h!]
    \centering
    \begin{tabular}{cc}
    (a) & (b) 
    \\
    \includegraphics[width=60mm]{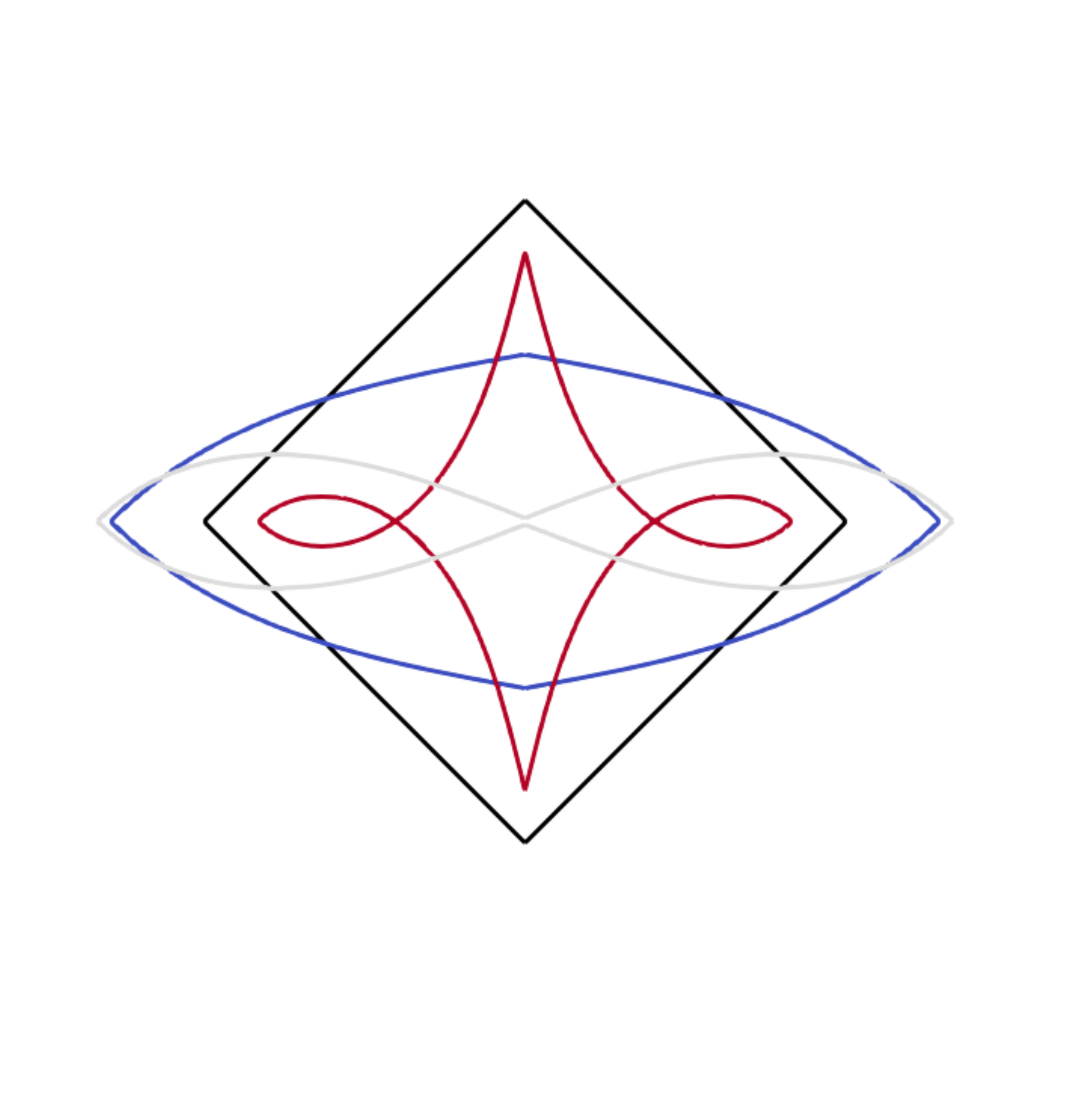}
    &  \includegraphics[width=60mm]{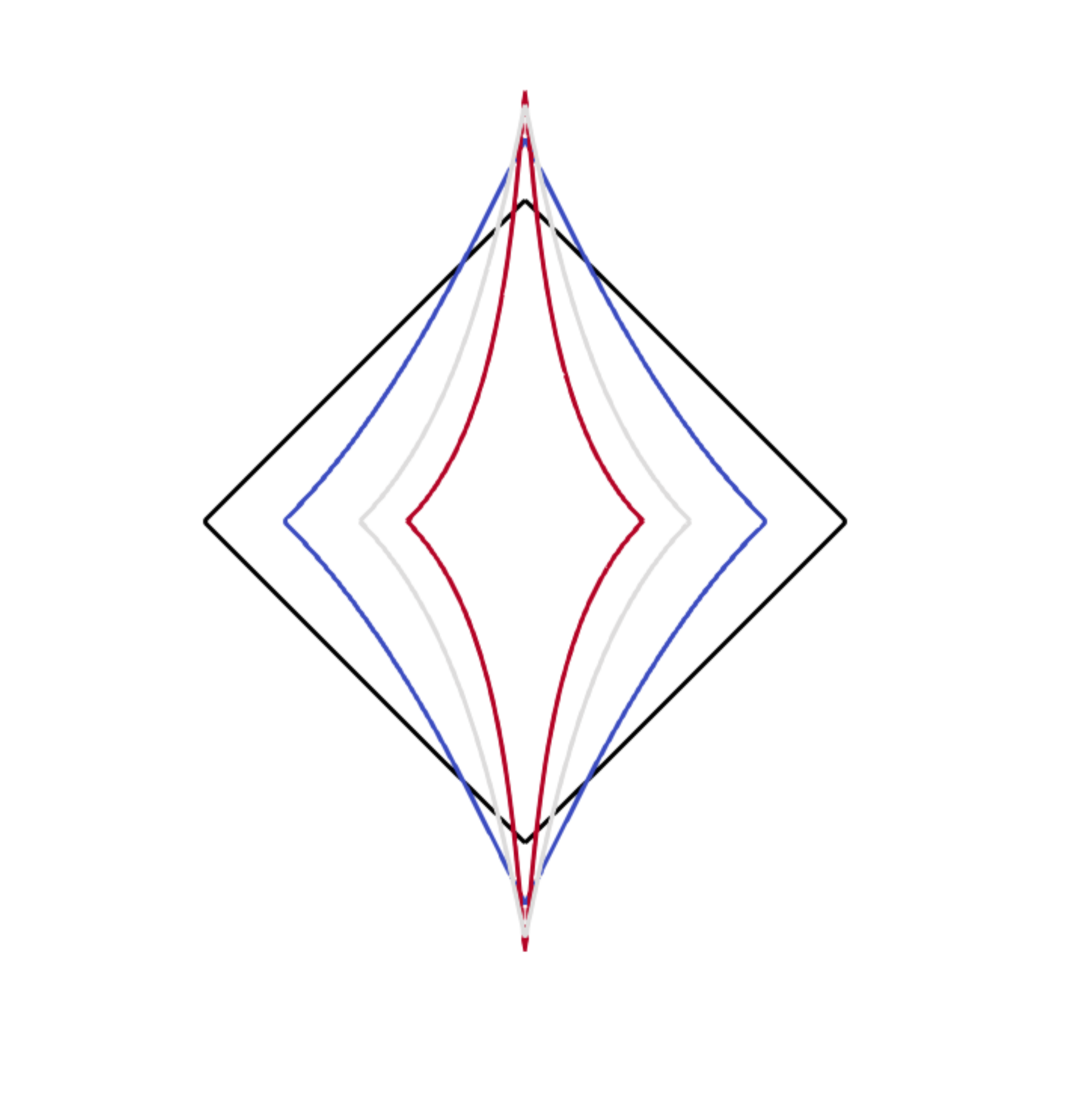}
    \end{tabular}
    \caption{Diamond frame: deformed shapes (scale factor equal to 1) for the (a) compression loading, (b) tension loading. The corresponding points in the load-displacement diagrams are indicated in  Fig. \ref{load_dispDF} by '+' symbols.}
    \label{deformedshape_DF}
\end{figure}

Let us now proceed to the diamond frame shown in 
Fig.~\ref{diamondframe}. The corresponding load-displacement diagrams are plotted in
Fig.~\ref{load_dispDF} and the deformed shapes
are shown in Fig.~\ref{deformedshape_DF}.
Red color in Fig.~\ref{diamondframe} corresponds to
the deflection $w_2$ and blue color to the horizontal displacement $u_1$.
The agreement of the present numerical results
with Mattiason's solution is again excellent.

\newpage
\subsection{Honeycomb lattice}
\label{sec:honeycomb}
%%%%%%%%%%%%%%%%%%%%%%%%%%%%%%

\subsubsection{Problem description}

In the last example we consider a material with internal
microstructure that corresponds to a two-dimensional elastic honeycomb lattice. Samples of such material will be subjected to prescribed displacements that would induce
uniaxial tension or compression if the material behaved as a homogeneous continuum. The objective
is to study the size effect, i.e., to investigate how
the macroscopic response of a finite sample deviates
from the limit behavior of an infinite lattice, which
can be under certain conditions described using a
hexagonal unit cell with imposed periodicity conditions.

\begin{figure}[h]
\centering
\begin{tabular}{ccc}
(a) & (b) & (c)\\
\\
\includegraphics[width=0.25 \linewidth]{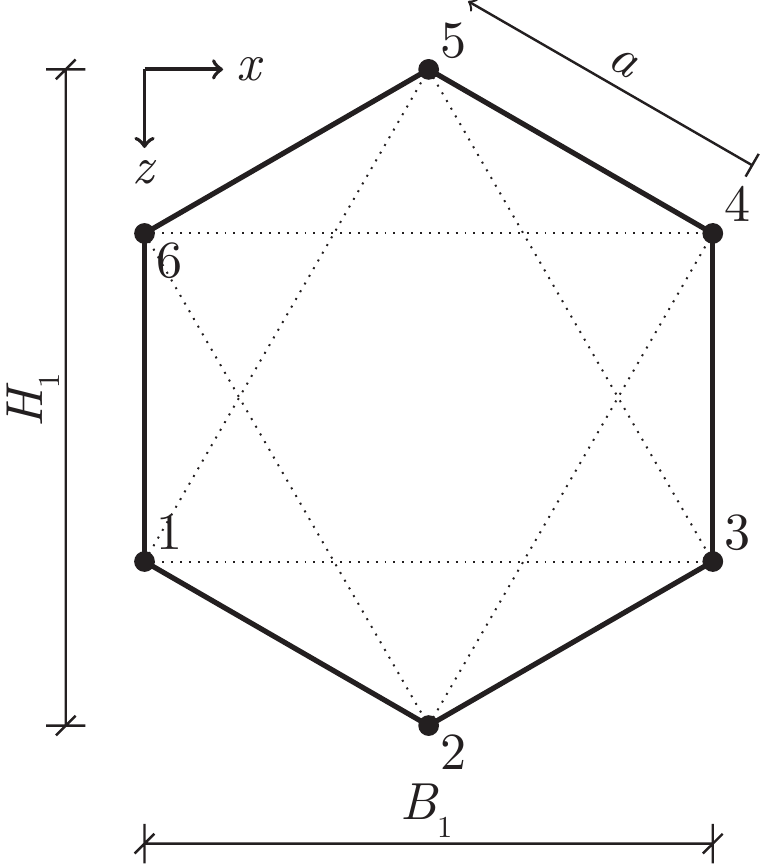} & \includegraphics[width=0.275 \linewidth]{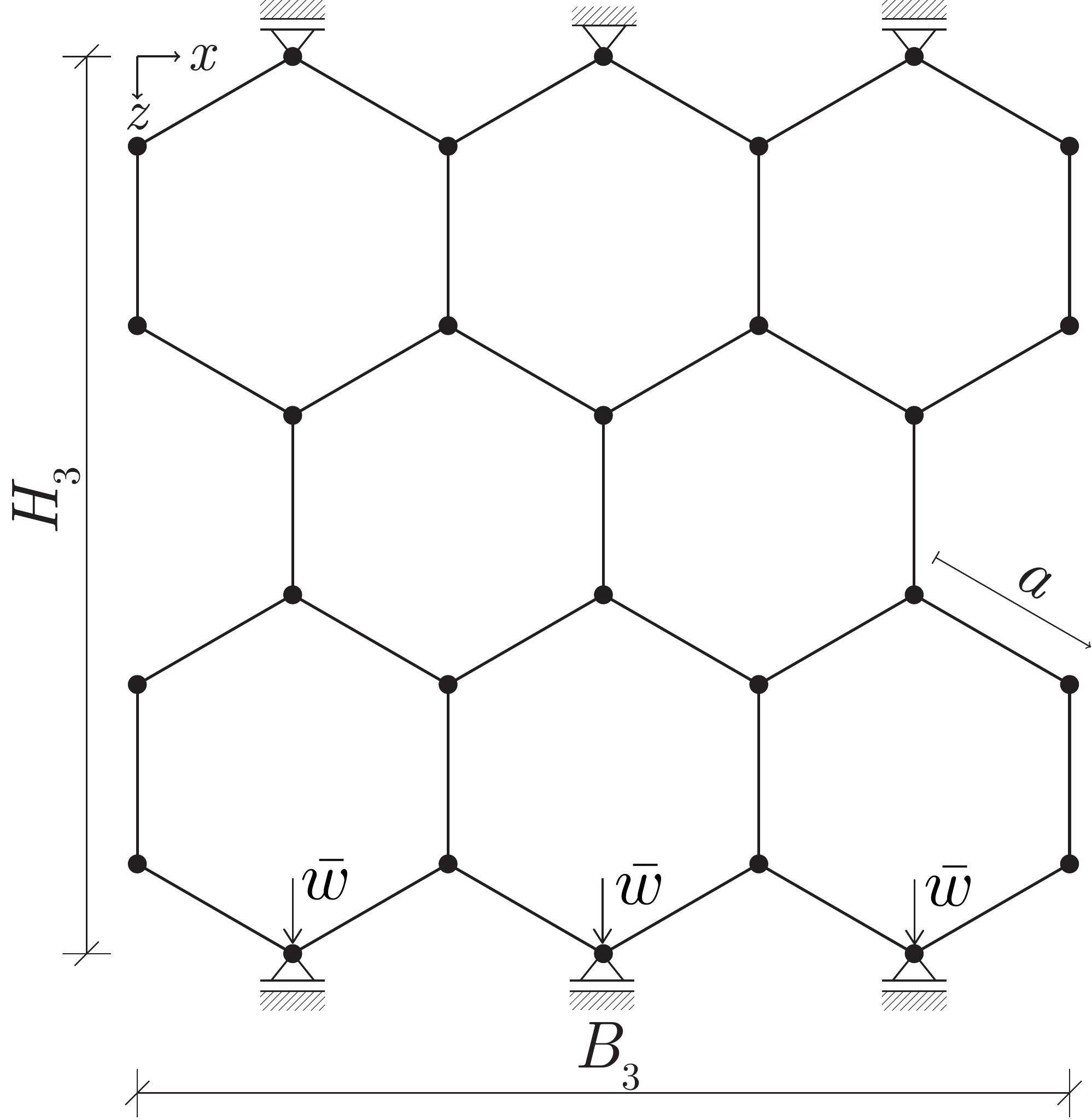} &
\includegraphics[width=0.29 \linewidth]{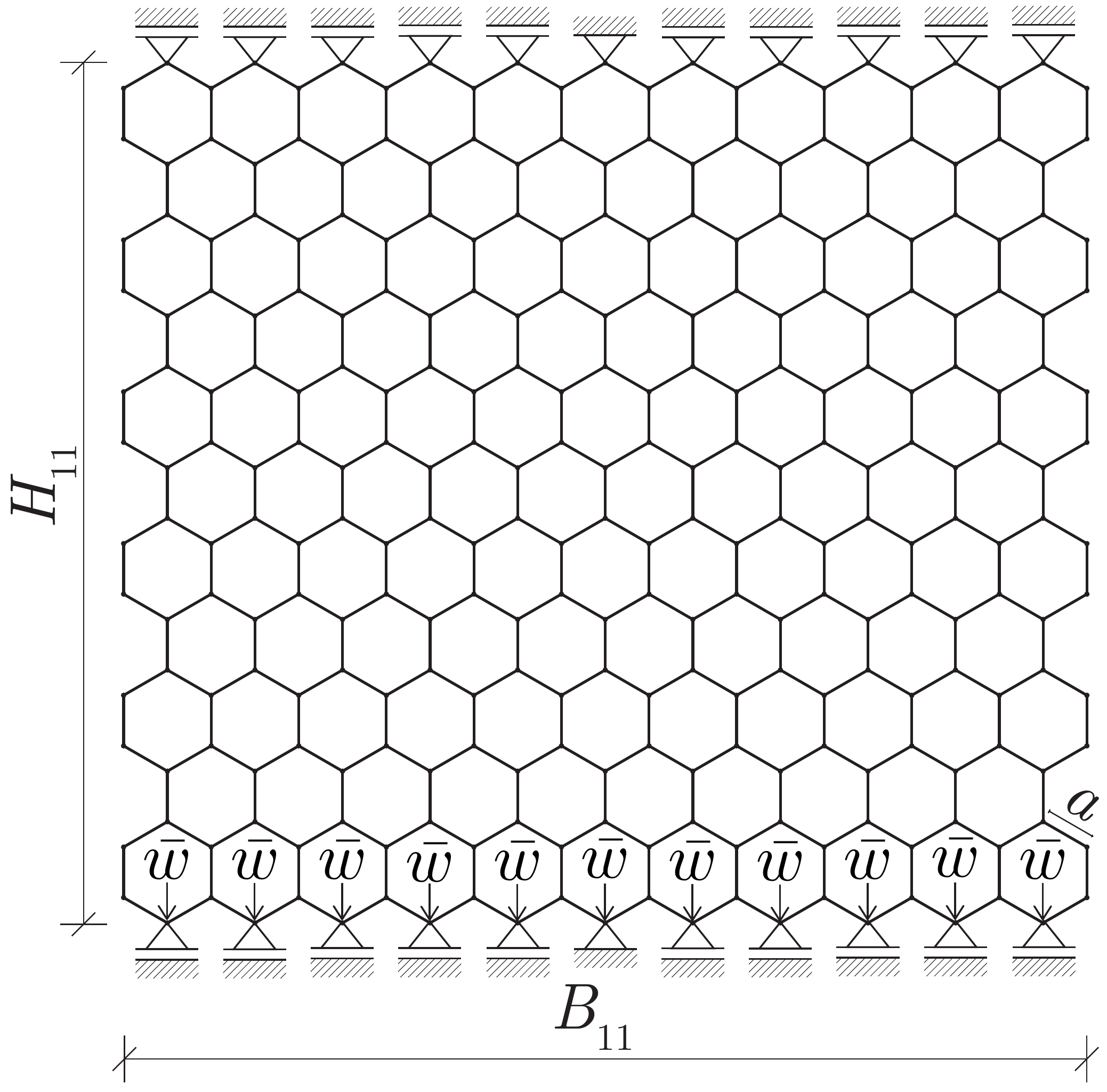}
\end{tabular}
\caption{Honeycomb lattices: (a) regular hexagonal cell with numbering of nodes and global coordinate system,
%Periodic nodes are symbolically connected through dashed lines. 
(b) $3 \times 3$ lattice, (c) $11 \times 11$ lattice.} 
\label{f:hex}
\end{figure}

The unit cell is a regular hexagon  consisting of six beam elements of length $a$, as shown in Fig.~\ref{f:hex}a.
Larger assemblies are obtained by stacking the cells 
horizontally and vertically in a honeycomb pattern.
For instance, the assembly in Fig.~\ref{f:hex}b will
be referred to as the $3\times 3$ pattern. It can be considered
as consisting of 3 layers; the odd layers contain
3 cells each while the even layer contains 2 full cells
and 2 half-cells. Similarly, the  assembly in Fig.~\ref{f:hex}c represents the $11\times 11$ pattern.

The projected dimensions of the unit cell are
$B_1=a\sqrt{3}$ horizontally and $H_1=2a$ vertically.
In general, an $n\times n$ lattice has width $B_n=na\sqrt{3}$
and height $H_n=(3n+1)a/2$. 
The lattices are loaded by tension or compression in the vertical direction. The loading is applied by prescribed
monotonically increasing 
vertical displacements $\bar{w}$ at $n$ bottom nodes with
coordinate $z=H_n$, while the $n$ top nodes with
coordinate $z=0$ are constrained vertically
(prescribed displacements $w=0$). All nodes are free
to move horizontally, only one node (an arbitrary one) 
is fixed horizontally in order to suppress rigid body translations.

%We present the numerical solution of a periodic honeycomb unit cell which is compared with the ones associated to larger lattices not considered as periodic media. 

%On the other hand, we consider the lattice system shown in Fig. \ref{f:hex} (b) constituted by three-by-three ($3\_3$) unit cells with 35 beam elements and 28 nodes. For both the unit cell and the larger lattice system we considered uniaxial tension and compression while periodicity was applied just to the unit cell. Similar to the honeycomb lattice (HL) showed in Fig. (\ref{f:hex}) (b), we considered larger representative volume elements named $5\_5$ HL, $7\_7$ HL, $9\_9$ HL, $11\_11$ HL whose details in terms of number of elements, nodes and characteristic dimensions are given in Table \ref{tab:hex}. Vertical displacements were imposed at nodes with zero-$z$ coordinate while the same degree of freedom was fixed for nodes at $z$-coordinate equal to $H$ (to prevent rigid body motion the displacement of one node is prevented in both directions).

%Figure \ref{f2:hex} show the stress strain curves obtained for the periodic unit cell (Fig. \ref{f2:hex} (a)) and for the larger lattice (Fig. \ref{f2:hex} (b)). 

\subsubsection{Hexagonal unit cell with periodic conditions}

\begin{figure}[h]
\centering
\begin{tabular}{cccc}
(a) & (b) & (c) & (d)\\
\\
\includegraphics[width=0.2 \linewidth]{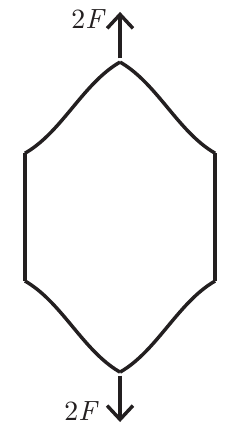} & \includegraphics[width=0.2 \linewidth]{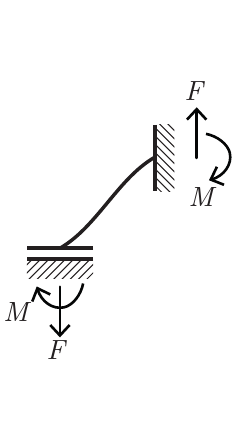} &
\includegraphics[width=0.2 \linewidth]{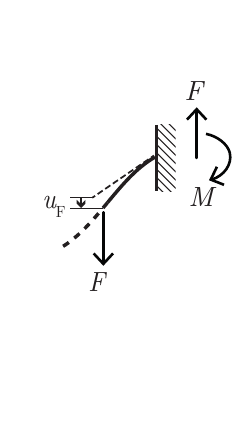} &
\includegraphics[width=0.2 \linewidth]{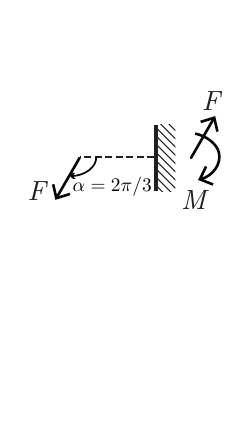}
\end{tabular}
\caption{Hexagonal unit cell with periodic conditions: (a) deformed cell under tension, (b) upper left inclined beam
with boundary conditions that follow from periodicity and double symmetry, (c)  upper half of the inclined beam, (d) rotated cantilever (the straight dashed line is the undeformed shape).} 
\label{f:periodic}
\end{figure}

For uniaxial tension or compression of an infinite lattice in the vertical direction, the solution can be expected
to exhibit periodicity and symmetry with respect to the 
horizontal and vertical axes of each hexagonal cell
(unless it is disturbed by instabilities, which are 
currently disregarded). 
Periodicity conditions combined with symmetry 
lead to zero rotations of all joints. 
Moreover, due to symmetry, vertical beams 1-6 and 3-4
deform only axially, and the deformation of
the four inclined beams can be obtained by mirroring
the deformation pattern of one of them (Fig.~\ref{f:periodic}a). Consequently,
it is sufficient to analyze only one selected inclined 
beam with one end clamped and the other forced to displace vertically
and allowed to move horizontally (Fig.~\ref{f:periodic}b). Once the forces in this
beam are computed, the axial force in vertical beams
is obtained from equilibrium and the macroscopic strain
and stress can be evaluated. 

In the incompressible limit, the behavior of
the inclined beam  in Fig.~\ref{f:periodic}b can be described analytically. One can even further reduce the problem to a cantilever of length $a/2$ shown in Fig.~\ref{f:periodic}c, because the inflexion point must be located at midspan of the original beam of length $a$. 
When the cantilever is analyzed in its local coordinate system, the problem is equivalent
to the one solved  in Appendix~\ref{sec:examplecantilever} and depicted in Fig.~\ref{f:cantilever},
with the direction of the applied force
inclined by $\alpha=2\pi/3$ with respect to the beam axis in the undeformed state (Fig.~\ref{f:periodic}d).
Formulae (\ref{eq:F})--(\ref{eq:a68})  derived in the appendix are directly applicable, with cantilever length $L$ set to $a/2$.

The end-point displacements $u_a$ and $w_a$ given by
(\ref{eq:ua})--(\ref{eq:wa})
refer to the local beam coordinates, 
but their projection onto the direction of the 
applied force, $u_F$, which can be evaluated from (\ref{eq:a68}), corresponds to one half of the difference
between global vertical displacements of nodes  5 and 6 of the unit cell.
Since the vertical displacements of nodes 6 and 1
are the same (due to inextensibility of the vertical beam that connects these nodes), the difference between 
global vertical displacements of nodes  5 and 1 is $2\,u_F$, and the corresponding macroscopic normal strain in the periodic lattice in the vertical direction is
\beq\label{eq:epsz}
\eps_z = \frac{2u_F}{3H_1/4}=\frac{4u_F}{3\,a}
\eeq
The force $F$ evaluated from (\ref{eq:F}) can be  converted into the macroscopic normal stress 
\beq\label{eq:sigz} 
\sigma_z = \frac{F}{tB_1/2} = \frac{2F}{ta\sqrt{3}} 
\eeq
Here, $t$ denotes the out-of-plane thickness.

Based on formulae (\ref{eq:F}) and (\ref{eq:a68}) and on the transformation 
of displacement and force into strain and stress given by (\ref{eq:epsz})--(\ref{eq:sigz}), it is possible to construct the stress-strain diagram. In (\ref{eq:F}) and (\ref{eq:a68}),
the force and the displacement are expressed
as functions of the rotation of the free end,
which plays the role of a parameter. 
In the analysis presented in Appendix~\ref{sec:examplecantilever},
 $F$ is considered as a positive quantity that represents the magnitude of the force, and the oriented direction is taken into account by an appropriate choice
of angle $\alpha$. For tension, $\alpha$ needs to be set to $2\,\pi/3$. Since the analytical solution derived in Appendix~\ref{sec:examplecantilever} is valid for 
$\alpha$ between 0 and $\pi$, vertical compression of the honeycomb lattice needs to be handled by setting 
$\alpha=\pi/3$ and adding negative signs in front
of the fractions in (\ref{eq:epsz})--(\ref{eq:sigz}).

%At very low loading levels, $\kappa_a$ must be almost proportional to $F$ and the fraction $\kappa_a^2/F$ becomes negligible. The initial value of parameter $\tilde k$ is thus $\sqrt{(1-\cos\beta)/2}=1/2$ and, since $\sin(\beta/2)=1/2$, we get $(1/\tilde{k})\sin(\beta/2)=1$
%and $\tilde a=F_J(\arcsin(1),1/2)=K(1/2)$
%then $\tilde k$ grows.

\begin{figure}[htb]
\centering
\begin{tabular}{cc}
(a) & (b) \\
\includegraphics[width=0.5 \linewidth]{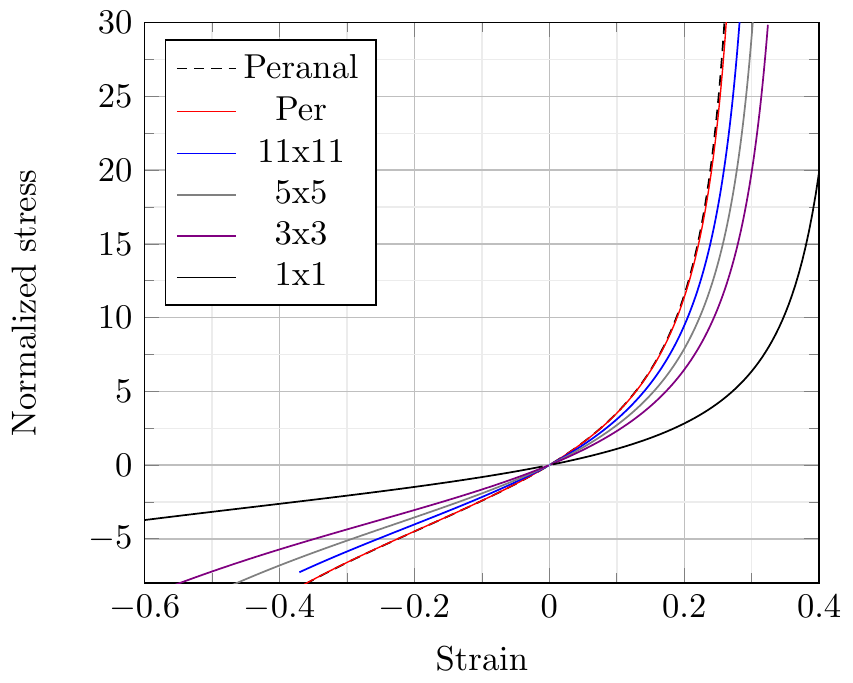}
&
\includegraphics[width=0.5 \linewidth]{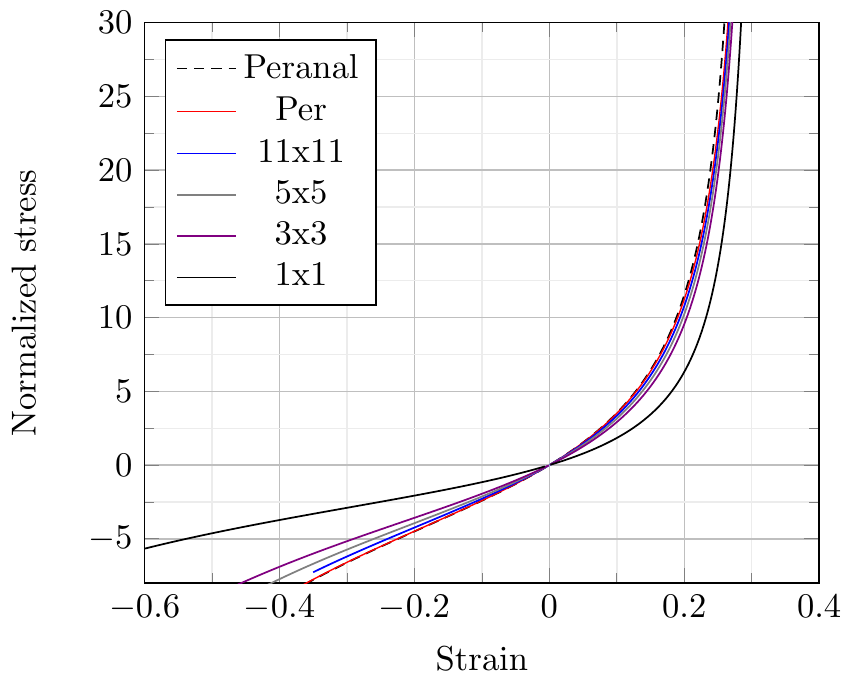}
\end{tabular}
\caption{Honeycomb lattice, complete stress-strain curves: (a) raw  results, (b) results with compensation for missing stiff layer.} 
\label{f:hex-all}
\end{figure}

\begin{figure}[h!]
\centering
\begin{tabular}{cc}
(a) & (b) \\
\includegraphics[width=0.5 \linewidth]{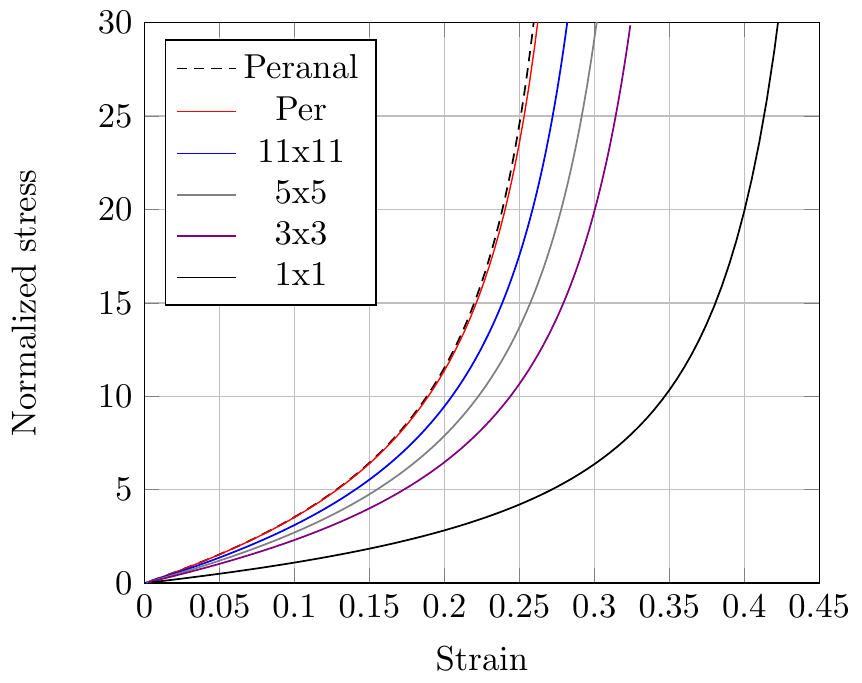}
&
\includegraphics[width=0.5 \linewidth]{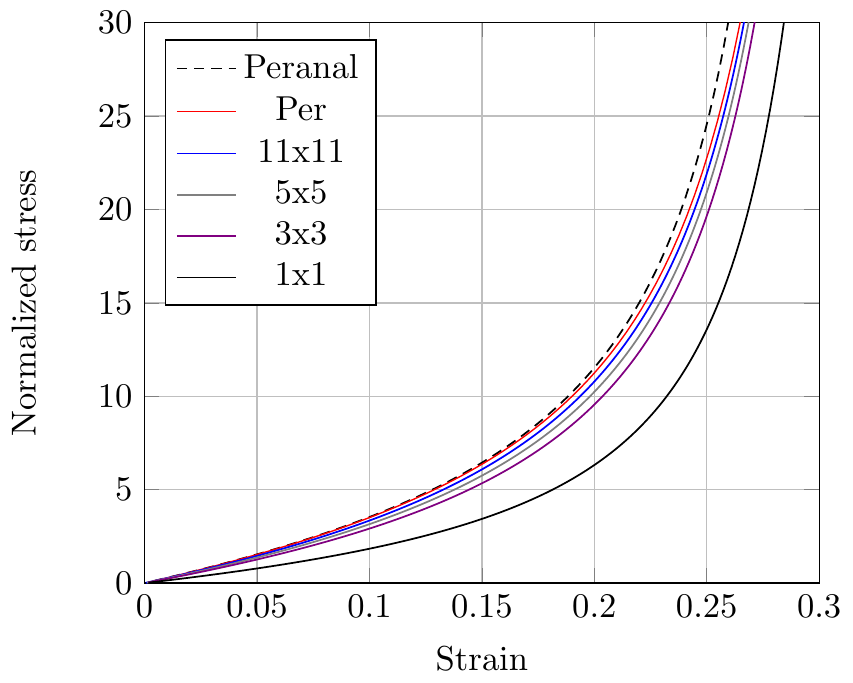}
\end{tabular}
\caption{Honeycomb lattice,  stress-strain curves for tension: (a) raw  results, (b) results with compensation for missing stiff layer.} 
\label{f:hex-ten}
\end{figure}

The macroscopic stress-strain curves extracted from
the results obtained for one hexagonal cell with imposed periodicity are plotted in both parts
of Fig.~\ref{f:hex-all}. The dashed black curve represents
the analytical solution derived for axially inextensible beams and the red curve has been obtained by
a numerical simulation of the inclined beam 5-6 considered
as extensible, with the contribution of the vertical
(also extensible) beam 6-1 added in closed form.
The dimensionless stress plotted in  Fig.~\ref{f:hex-all} is 
the actual macroscopic stress divided by the
normalizing factor $EI/(ta^3)$. The simulations have been
performed with a relatively high normal stiffness characterized by the
dimensionless parameter $EAa^2/EI=10,000$.
The numerical and analytical results are in very good
agreement and the red and dashed black curves slightly differ only in the 
regime of high tension, as seen more clearly in Fig.~\ref{f:hex-ten}.

\subsubsection{Effect of lattice size}

\begin{figure}[h!]
\hspace{-2.8cm}
\begin{tabular}{cc}
(a) & (b) \\
\includegraphics[width=0.55 \linewidth]{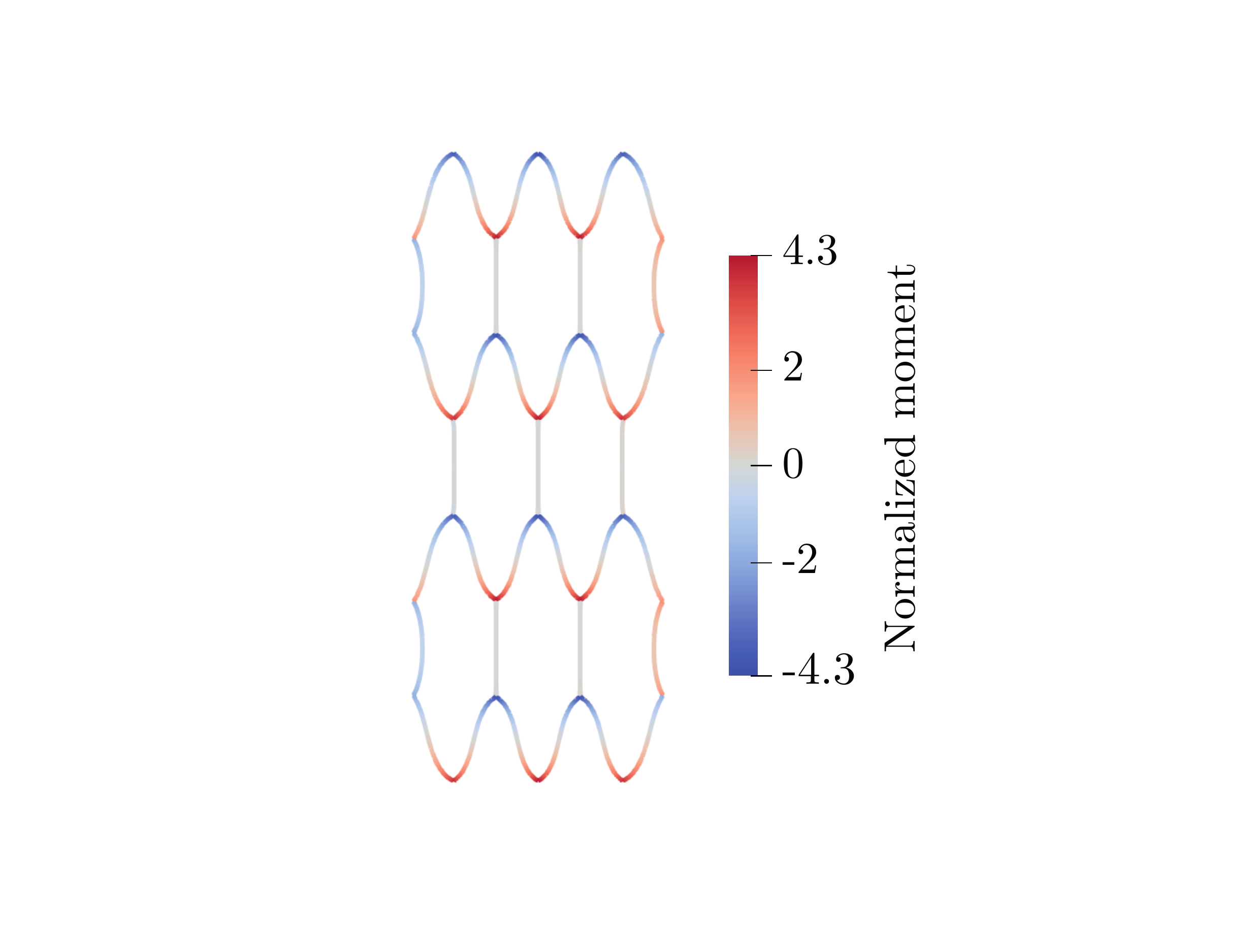}
&
\includegraphics[width=0.55 \linewidth]{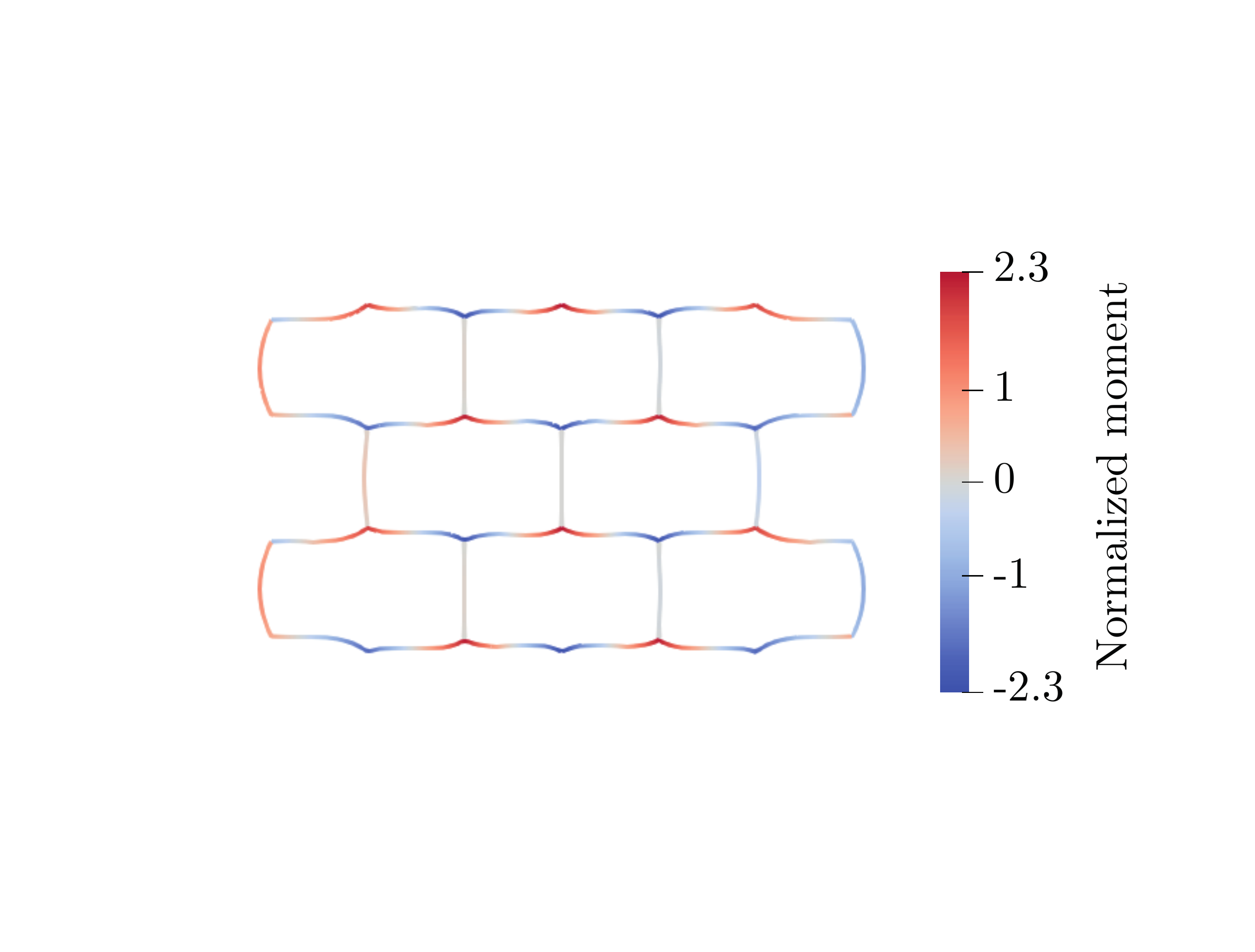}\\
(c) & (d) \\
\includegraphics[width=0.55 \linewidth]{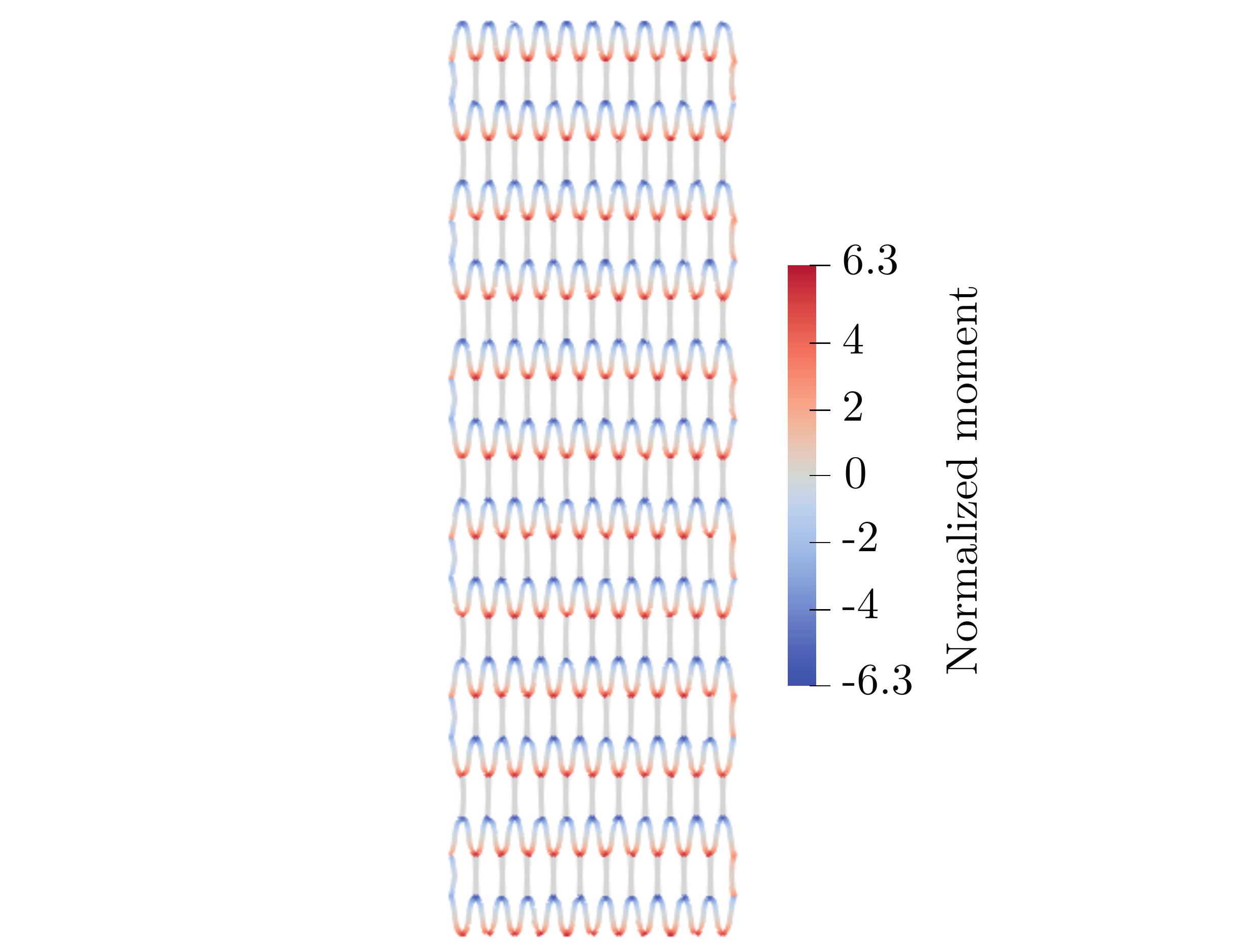}
&
\includegraphics[width=0.55 \linewidth]{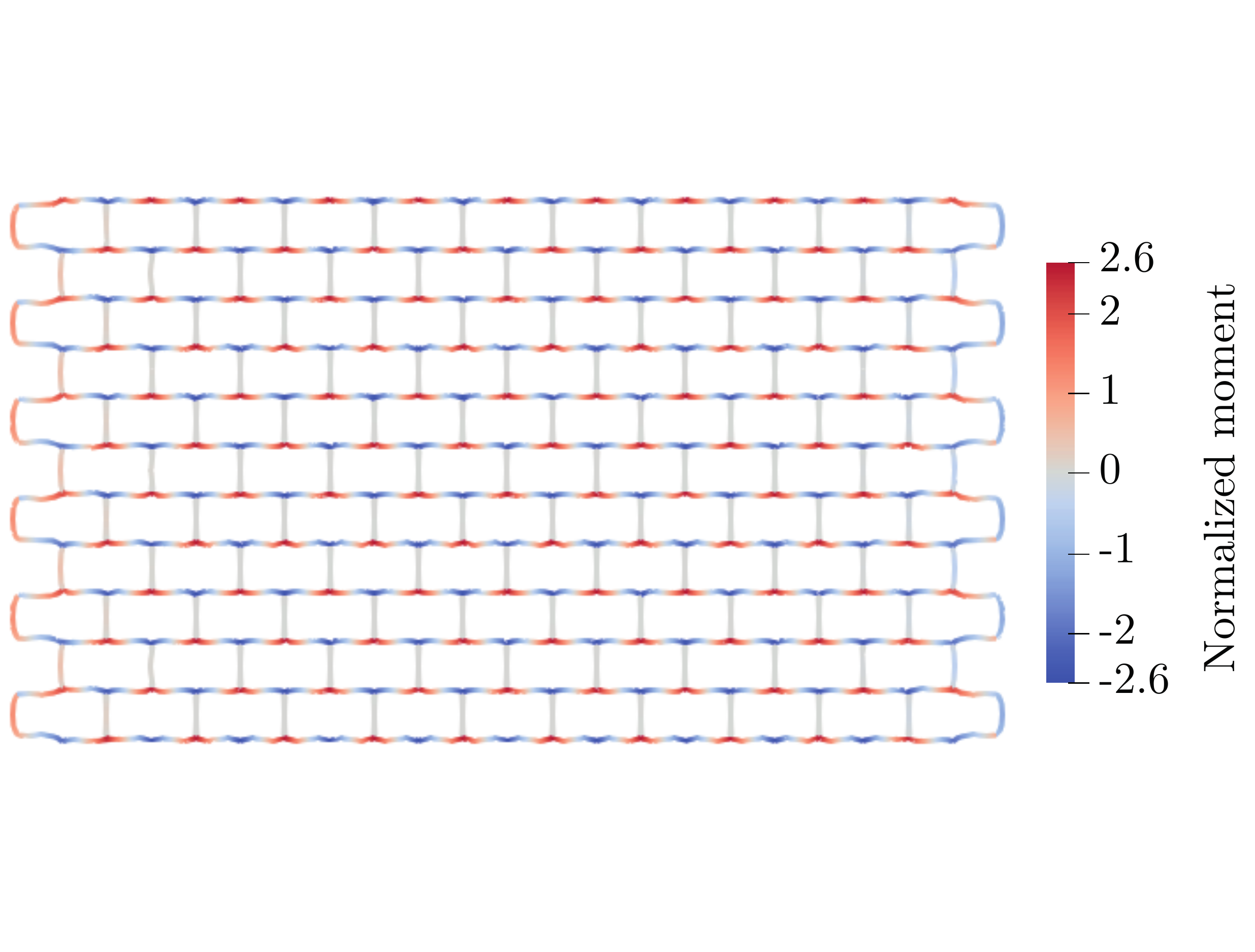}
\end{tabular}
\caption{Honeycomb lattice, deformed shapes at 30$\%$ of strain level: (a) $3 \times 3$ lattice loaded in tension, (b) $3 \times 3$ lattice loaded in compression, (a) $11 \times 11$ lattice loaded in tension, (b) $11 \times 11$ lattice loaded in compression. Colour maps refer to the normalized moment, i.e. moment divided by the factor $EI/a$.} 
\label{f:hex-ds}
\end{figure}

The solutions that satisfy periodicity conditions
correspond to the theoretical limit of an infinite
lattice filling the whole plane.
To assess the effect of finite size,
numerical simulations have been performed on finite 
$n\times n$ lattices exemplified in Fig.~\ref{f:hex}b,c,
with $n=1,3,5$ and 11. The imposed vertical displacement
$\bar w$ and the computed reactions at the 
bottom nodes, $R_i$, $i=1,2,\ldots n$, have been 
converted into the corresponding average
stress and strain,   
\bea\label{honey:sigma} 
\sigma &=& \frac{\sum_{i=1}^n R_i}{tB_n} \\
\label{honey:eps} 
\eps &=& \frac{\bar w}{H_n}
\eea
Stress-strain curves obtained in this way are shown
in Figs.~\ref{f:hex-all}a and \ref{f:hex-ten}a.
A strong size effect is observed---smaller samples 
lead to a more compliant response, both in tension and
in compression. The overall shapes of all stress-strain
curves are similar. Only a slight nonlinearity is observed
in compression while the tensile response 
is highly nonlinear, with a dramatic increase of
tangent stiffness at average strains exceeding 10 \%.
This stiffening is caused by the fact that the inclined
beams initially deform by bending but this relatively
soft deformation mode has a limited capacity and 
as the inclined beams get aligned with the direction
of applied tensile loading, their high axial stiffness
is activated. This is nicely illustrated by the deformed shapes in parts (a) and (c) of
Fig.~\ref{f:hex-ds}.

%\begin{figure}[h!]
%\centering
%\begin{tabular}{cc}
%(a) & (b) \\
%\includegraphics[width=0.5 %\linewidth]{figures/honeycomblattices5.pdf}
%&
%\includegraphics[width=0.5 \linewidth]{figures/honeycomblattices6.pdf}
%\end{tabular}
%\caption{Honeycomb lattice,  stress-strain curves for compression: (a) raw  results, (b) results with compensation for missing stiff layer.} 
%\label{f:hex-com}
%\end{figure}

The size effect exhibited by the stress-strain diagrams in Fig.~\ref{f:hex-all}a is to some extent caused by
the fact that the finite lattices of the kind
depicted in Fig.~\ref{f:hex}b,c do not represent a decomposition of the infinite lattice into periodically
repeatable units. Indeed, when we stack two such lattices
vertically, an additional layer of vertical beams 
needs to be inserted in between, and when we stack them
horizontally, the vertical beams on the boundaries
that are now glued together would be doubled. 

The first effect turns out to be stronger than the
second one. It can be eliminated by adding an extra
layer of vertical beams to the nodes at the bottom
of each lattice. The resulting modified
lattice would  not be practical for testing but
its numerical treatment is straightforward. 
In fact, it is not even necessary to perform 
additional numerical simulations of the modified
lattices because the effect of the added layer on
the stress-strain diagram can easily be estimated.
At a given stress $\sigma$ evaluated from
(\ref{honey:sigma}), the average axial force in the
vertical beams is equal to the average reaction, 
$\bar{R}=\sigma tB_n/n$, and the contribution of the
added layer of axially deformed beams to the vertical displacement on the boundary is $\Delta\bar w=\bar{R}a/EA$. At the same time, the added layer
increases the height of the sample by $\Delta H=a$.
The corrected average strain is then estimated as
\beq
\label{honey:epscor} 
\eps+\Delta\eps = \frac{\bar w+\Delta\bar w}{H_n+\Delta H} = \frac{\eps(3n+1)a/2+\bar{R}a/EA}{(3n+1)a/2+a}
=\frac{\eps(3n+1)+2\sigma ta\sqrt{3}/EA}{3n+3}
= \frac{3n+1}{3(n+1)}\,\eps + \frac{2 ta}{\sqrt{3}\,(n+1)EA}\,\sigma
\eeq
from which
\beq
\label{honey:epscor2} 
\Delta\eps =-\frac{2}{3(n+1)}\,\eps + \frac{2 }{\sqrt{3}\,(n+1)}\,\frac{EI }{EAa^2}\,\frac{\sigma ta^3}{EI}
\eeq
For beams with high axial stiffness, the second term on the right-hand side is usually
negligible compared to the first term. It is expressed
as a product of three fractions, the second of which is the reciprocal value of the dimensionless stiffness
coefficient $EAa^2/EI$, in our example
equal to 10,000. In Fig.~\ref{f:hex-all}, normalized stresses $\sigma ta^3/EI$ do not exceed 30, and so even for
$n=1$ the  strain correction described by the second term
 on the right-hand side of (\ref{honey:epscor2})
is $30/(10,000\cdot\sqrt{3})\approx 0.0017$.
On the other hand, the first term on the right-hand side of (\ref{honey:epscor2}) is important, especially
for small $n$. The strain correction $\Delta\eps$ has
the opposite sign than the originally evaluated strain
$\eps$,
which means that the correction leads to stiffer
response.

Stress-strain diagrams for the modified lattices
with an added layer of vertical beams are shown
in . The size effect is reduced, but it is still
present, especially for small lattices. 
The residual size effect originates from 
softer response of cells located near the lateral (vertical)
boundaries of the sample. For a large lattice,
cells that are far from the lateral boundaries
deform in a pattern similar to the periodic cell,
i.e., with negligible rotations. On the other hand,
cells located in boundary layers are less constrained
and deform more easily. Images of deformed lattices
in Fig.~\ref{f:hex-ds} indicate that nodes on the
lateral boundaries rotate and the deformation is
more equally distributed among the vertical and
inclined beams, which reduces the apparent macroscopic stiffness of the sample. Since the affected boundary
layers occupy a relatively larger area fraction 
in a small sample than in a large one, smaller
samples behave as if the material were softer.

\section{Concluding remarks}
\label{sec:conclusions}

The %two-dimensional 
geometrically nonlinear formulation of a beam element
developed in this paper is based on kinematic equations
valid for arbitrarily large rotations of cross sections. 
These equations are combined with 
the integrated form of equilibrium equations and with
generalized material equations that link the internal forces to
the deformation variables (curvature and centerline stretch).
The resulting fundamental equations have the form of
a set of three first-order ordinary differential equations
for two displacements and one rotation as basic unknowns.
Instead of approximating these functions by apriori
selected shape functions, it is proposed to convert
the boundary value problem into an initial value problem 
using a technique inspired by the
shooting method, and then to 
adopt a finite difference scheme. 
The advantage is that accuracy of the numerical approximation 
can be conveniently increased by 
refining the integration scheme on the element level 
while the number of global degrees of freedom is kept constant.

A number of examples have been presented to illustrate
the flexibility and efficiency of the proposed approach and to assess its
accuracy and robustness.
In all examples, the number of global degrees of freedom has been kept as low
as possible. Spatial refinement has been taken care of by
reducing the spacing between auxiliary grid points used for
finite-difference approximation of the fundamental equations.
The results have been shown to be in excellent agreement with
previously published numerical results or analytical solutions.
For the Williams toggle frames, a good agreement with 
experimental data has been observed. 

It has been shown that when the integration grid is refined,
the error decreases in proportion to the square of the grid spacing. The accuracy is comparable to results obtained
when the beam is discretized by standard finite elements
of the same size as the grid spacing. However, the present approach
has the big advantage that no additional global
degrees of freedom need to be introduced when the grid is refined.
In the context of an incremental iterative analysis on the
global (structural) level, relatively large increments of
nodal displacements and rotations can be handled. 
The shooting method adopted on the element level at the same time
permits to set up the element tangent stiffness matrix,
which is then processed by a standard assembly procedure to
construct the structural tangent stiffness. Linearization
of the global equations is consistent, as confirmed by
quadratic convergence of the global Newton-Raphson iteration
procedure. This process is more robust for low values of
the axial stiffness (reflected by the dimensionless stiffness
parameter $EAL^2/EI$), and for higher values it works fine
provided that the increment size is not extreme.

To keep the paper focused and its size limited, 
we have presented a basic version of the proposed
approach, with a number of simplifying assumptions.
For instance, the stress-strain law has been assumed
to be linear elastic, and the beam element has been
considered as two-dimensional and initially straight.
This relatively simple setting permits to explain
the main ideas and highlight the essence of the
numerical procedure. However, extensions and
generalizations are possible and some of them
are currently being explored. In particular,
an extension to beam elements with initial curvature
seems to be potentially attractive and powerful,
and it will be reported on in a follow-up publication.

%\backmatter

\section*{Acknowledgments}
The authors acknowledge the support of the Czech Science Foundation (project No. 19-26143X).

%\subsection*{Author contributions}
%This is an author contribution text. This is an author contribution text. This is an author contribution text. This is an author contribution text. This is an author contribution text. 

%\subsection*{Financial disclosure}
%None reported.

\subsection*{Conflict of interest}
The authors declare no potential conflict of interests.

%\section*{Supporting information}

%The following supporting information is available as part of the online article:

%\noindent
%\textbf{Figure S1.}
%{500{\uns}hPa geopotential anomalies for GC2C calculated against the ERA Interim reanalysis. The period is 1989--2008.}

%\noindent
%\textbf{Figure S2.}
%{The SST anomalies for GC2C calculated against the observations (OIsst).}

\newpage
\appendix

%%%%%%%%%%%%%%%%%%%%%%%%%%%%%%%%%%%%%%%%%%%%%%%
\section{Nonlinear elastic beam -- analytical treatment}\label{app:analytical}
%%%%%%%%%%%%%%%%%%%%%%%%%%%%%%%%%%%%%%%%%%%%%%%

\subsection{Solution of second-order equation for sectional rotation}\label{app:anal1}

The problem described by equations (\ref{e212})--(\ref{e214}), which were derived in Section~\ref{sec:fund1} and form the basis of the 
algorithms used in this paper, can be handled 
analytically.
It is convenient to start from the transformed version
of these equations given in (\ref{e221}), because it contains the rotation function $\varphi$ as the only unknown.
Equation (\ref{e221}) is a nonlinear second-order differential equation
in a special form 
\beq\label{eq:varphif} 
\varphi'' + g(\varphi) = 0
\eeq 
and it always allows for a formal ``analytical'' solution.
The usefulness of the
result depends on the existence of closed-form expressions
for integrals that are arise in the solution process.

To derive a formal analytical solution, let us
first multiply (\ref{eq:varphif}) by $2\varphi'$ and
then express the resulting equation as
\beq\label{eq:varphiF} 
\left(\varphi'^2 + 2G(\varphi)\right)' = 0
\eeq 
where $G$ is the indefinite integral (antiderivative) of $g$. The resulting first-order equation
\beq\label{eq:varphiF2} 
\varphi'^2 + 2G(\varphi) = C
\eeq 
can be handled by separation of variables, which leads to
\beq\label{eq:varphiF3} 
\pm\frac{\dif\varphi}{\sqrt{C-2G(\varphi)}} = \dif x
\eeq 
and finally to the implicit formula for the solution,
\beq\label{eq:varphiF4} 
\pm\int_{\varphi_0}^{\varphi(x)}\frac{\dif\bar\varphi}{\sqrt{C-2G(\bar\varphi) }} = x-x_0
\eeq 
Here, $C$ is an arbitrary integration constant,
and $\varphi_0$ is the value of $\varphi$ at $x=x_0$.
%in our case typically equal to 0.
The sign to be used in front of the integral
corresponds to the sign of the first derivative
of $\varphi$ in the interval of interest. 

When this technique is applied to equation  (\ref{e221}), 
the corresponding function $g$ is given by
\bea\nonumber
g(\varphi) &=& \frac{1}{EI}\left(1+\frac{-X_{ab}\cos\varphi + Z_{ab}\sin\varphi}{EA}\right)
\left(X_{ab}\sin\varphi + Z_{ab}\cos\varphi\right) =
\\
&=& \frac{X_{ab}}{EI}\sin\varphi + \frac{Z_{ab}}{EI}\cos\varphi +\frac{Z_{ab}^2-X_{ab}^2}{2EIEA}\sin 2\varphi - \frac{Z_{ab}X_{ab}}{EIEA}\cos 2\varphi
\label{ee250}
\eea
Integration with respect to $\varphi$ is easy and the 
antiderivative of $g$ can be selected as
\beq \label{ee251}
G(\varphi) = \frac{X_{ab}}{EI}\left(1-\cos\varphi\right) + \frac{Z_{ab}}{EI}\sin\varphi +\frac{Z_{ab}^2-X_{ab}^2}{4EIEA}\left(1-\cos 2\varphi\right) - \frac{Z_{ab}X_{ab}}{2EIEA}\sin 2\varphi
\eeq 
For initial conditions
\bea 
\varphi(0) &=& \varphi_a \\
\varphi'(0) &=& \kappa_a
\eea 
we obtain integration constant
\beq\label{ee254} 
C = \kappa_a^2 + 2G(\varphi_a)
\eeq 
expressed in terms of the left-end rotation, $\varphi_a$, and left-end curvature, $\kappa_a$.

Up to here, everything has been quite straightforward.
However, in the next step we need to evaluate the integral on the left-hand side of (\ref{eq:varphiF4}),
which is given by
\beq 
\int \frac{{\rm d}\bar\varphi}{\sqrt{\kappa_a^2+2G(\varphi_a)-2G(\bar\varphi)}} = \int \frac{{\rm d}\bar\varphi}{\sqrt{c_1+c_2\cos\bar\varphi+c_3\sin\bar\varphi+c_4\cos2\bar\varphi+c_5\sin2\bar\varphi}}
\eeq 
In general, this would be very difficult. In the special case
of an axially inextensible beam, the problem is simplified
because coefficients $c_4$ and $c_5$, which multiply
the terms with $\cos 2\bar\varphi$ and $\sin 2\bar\varphi$,
vanish (they contain $EA$ in the denominator). The remaining coefficients are given by
\bea \label{ee1}
c_1 &=& \kappa_a^2 -\frac{2X_{ab}}{EI}\cos\varphi_a+\frac{2Z_{ab}}{EI}\sin\varphi_a
= \kappa_a^2 +\frac{2N_{ab}}{EI}
\\
c_2 &=& \frac{2X_{ab}}{EI} \label{ee2} \\
c_3 &=& -\frac{2Z_{ab}}{EI} \label{ee3}
\eea 
in which
\beq 
N_{ab} = -X_{ab}\cos\varphi_a +Z_{ab}\sin\varphi_a
\eeq 
is the normal force at the left end section.
Let us introduce auxiliary constants,
\bea \label{eq:Fab}
F_{ab} &=& \sqrt{X_{ab}^2+Z_{ab}^2}
\\
A &=& \sqrt{c_2^2+c_3^2}  =\frac{2}{EI}\sqrt{X_{ab}^2+Z_{ab}^2}=\frac{2F_{ab}}{EI} \label{ee4}\\
\alpha &=& -\arctan\frac{c_3}{c_2} = \arctan\frac{Z_{ab}}{X_{ab}} \label{ee5}
\eea
and a transformed variable,
\beq 
\tilde\varphi = \frac{\bar\varphi+\alpha}{2}
\eeq 
Their purpose is to replace $c_2\cos\bar\varphi+c_3\sin\bar\varphi$ by
$A\cos 2\tilde\varphi = A(1-2\sin^2\tilde\varphi)$.
Then we can proceed to the integral
\bea\nonumber
\int_{\varphi_a}^{\varphi} \frac{{\rm d}\bar\varphi}{\sqrt{c_1+c_2\cos\bar\varphi+c_3\sin\bar\varphi}} &=& \int_{(\varphi_a+\alpha)/2}^{(\varphi+\alpha)/2} \frac{2\,{\rm d}\tilde\varphi}{\sqrt{c_1+A(1-2\sin^2\tilde\varphi)}} =  \frac{2}{\sqrt{c_1+A}}\int_{(\varphi_a+\alpha)/2}^{(\varphi+\alpha)/2} \frac{{\rm d}\tilde\varphi}{\sqrt{1-\displaystyle\frac{2A}{c_1+A}\sin^2\tilde\varphi}}=
\\
\label{eq267}
&=& \frac{2}{\sqrt{c_1+A}}\int_{(\varphi_a+\alpha)/2}^{(\varphi+\alpha)/2} \frac{{\rm d}\tilde\varphi}{\sqrt{1-k^2\sin^2\tilde\varphi}}
\eea
where
\beq 
k=\sqrt{\frac{2A}{c_1+A}} = %\frac{2}{\sqrt{4-2\cos(\varphi_a+\alpha)+EI\kappa_a^2/\sqrt{X_{ab}^2+Z_{ab}^2}}}
\sqrt{\frac{4F_{ab}}{2(F_{ab}+N_{ab})+EI\kappa_a^2}} \label{eee264}
\eeq 
Recall that $F_{ab}$ given by (\ref{eq:Fab}) is the magnitude of the end force, which is by definition non-negative.
It is worth noting that the denominator
of the fraction under the square root in (\ref{eee264}), given by $2(F_{ab}+N_{ab})+EI\kappa_a^2$, is also non-negative because $N_{ab}\ge -F_{ab}$. This denominator could be zero only if
$N_{ab}=-F_{ab}$ and $\kappa_a=0$ (equivalent to $M_{ab}=0$), which is the case of uniaxial tension, leading to the trivial solution  $\varphi(x)=0$.

Now we can substitute the right-hand side of (\ref{eq267}) into (\ref{eq:varphiF4}), setting $x_0=0$
and $\varphi_0=\varphi_a$ and selecting the sign
in front of the integral as ${\rm sgn}\,\kappa_a$,
so that it agrees with the sign of $\varphi'$
at $x=0$. If the curvature at the left end, $\kappa_a$, happens
to be zero (in cases when the end moment $M_{ab}$ vanishes), the sign should correspond to the expected sign of the curvature for small positive values
of $x$, which can be deduced from $\alpha$ and $\varphi_a$.

The integral on the left-hand
side of (\ref{eq:varphiF4}) is 
evaluated for function 
\beq \label{ee251x}
G(\varphi) = \frac{X_{ab}}{EI}\left(1-\cos\varphi\right) + \frac{Z_{ab}}{EI}\sin\varphi 
\eeq 
which is the reduced version of (\ref{ee251})
valid for an axially inextensible beam ($EA\to\infty$). Constant $C$ is substituted
from (\ref{ee254}) and the integral is expressed based on (\ref{eq267}).
The resulting equation
\beq \label{ee263}
\frac{2\,{\rm sgn}\,\kappa_a}{\sqrt{c_1+A}}\int_{(\varphi_a+\alpha)/2}^{(\varphi(x)+\alpha)/2} \frac{{\rm d}\tilde\varphi}{\sqrt{1-k^2\sin^2\tilde\varphi}} = x
\eeq 
implicitly defines function $\varphi(x)$ that
describes the sectional rotation.

%%%%%%%%%%%%%%%%%%%%%%%%%%%%%%%%%%%%%%%%%%%%%%
\subsection{Mathematical tools -- elliptic integrals and elliptic functions}

The integral on the right-hand side of (\ref{ee263}) is recognized as one of the elliptic integrals.
Before we proceed with the solution, let us recall the definitions of
elliptic integrals and elliptic functions,
which  will later be used in analytical expressions
describing the rotation and displacement functions.
A systematic overview can be found in standard
mathematical literature, e.g., in \cite{ByrdFriedman54}.

The {\bf incomplete elliptic integral of the first kind} is  given by 
\beq \label{ee264}
{\rm F_J}(\varphi,k) = \int_0^\varphi \frac{{\rm d}t}{\sqrt{1-k^2\sin^2t}}
\eeq
If the upper bound in  the integral  is set
to $\pi/2$, we obtain the {\bf complete elliptic integral of the first kind},
\beq \label{ee264c}
{\rm K}(k)={\rm F_J}(\pi/2,k) = \int_0^{\pi/2} \frac{{\rm d}t}{\sqrt{1-k^2\sin^2t}}
\eeq

The {\bf Jacobi amplitude function} ${\rm am}(x,k)$ is the inverse of ${\rm F_J}$ with respect
to $\varphi$, with $k$ considered as a fixed parameter (again, usually in the range between 0 and 1).
This means that 
\beq 
{\rm am}(x,k) = \varphi
\eeq
is equivalent with
\beq 
x = {\rm F_J}(\varphi,k)
\eeq

The {\bf elliptic sine} and {\bf elliptic cosine} are defined as
\bea\label{eq:sn} 
{\rm sn}(x,k) &=& \sin({\rm am}(x,k))
\\ \label{eq:cn} 
{\rm cn}(x,k) &=& \cos({\rm am}(x,k))
\eea
and belong to the family of {\bf Jacobi elliptic functions}. Another useful member of this family
is the so-called {\bf delta amplitude}
\beq \label{eq:dn} 
{\rm dn}(x,k) = \sqrt{1-k^2\,{\rm sn}^2(x,k)}
\eeq 

Finally, the {\bf incomplete and complete elliptic integrals of the second kind} are defined as
\bea\label{eq:EJ}
{\rm E_J}(\varphi,k) &=& \int_0^\varphi \sqrt{1-k^2\sin^2t}\,{\rm d}t
\\
{\rm E}(k) &=& {\rm E_J}(\pi/2,k) = \int_0^{\pi/2} \sqrt{1-k^2\sin^2t}\,{\rm d}t
\eea

In the above expressions, parameter $k$ is usually considered in the range between 0 and 1. The integrals in (\ref{ee264}) and (\ref{eq:EJ}) are well defined even for $k>1$ as long as $\varphi$ remains below $\arcsin(1/k)$. However, some implementations of incomplete elliptic integrals and elliptic functions consider the case of $k>1$ as inadmissible, for any $\varphi$. One can then exploit the transformation
\beq \label{ee264a}
{\rm F_J}(\varphi,k) = \int_0^\varphi \frac{{\rm d}t}{\sqrt{1-k^2\sin^2t}} = \frac{1}{k}\int_0^{\arcsin (k\sin\varphi)} \frac{{\rm d}s}{\sqrt{1-k^{-2}\sin^2s}} = \frac{1}{k}{\rm F_J}(\arcsin (k\sin\varphi),1/k)
\eeq
and use function ${\rm F_J}$ with parameter $k>1$ replaced by its reciprocal value, $\tilde{k}=1/k<1$. 
Based on (\ref{ee264a}), we can also write
\beq\label{eq:amtrans} 
{\rm am}(x,k) = \arcsin\frac{\sin{\rm am}(kx,1/k)}{k}
\eeq 
In terms of the elliptic sine function, formula
(\ref{eq:amtrans}) can be rewritten in an elegant form
\beq 
{\rm sn}(x,k) = \frac{1}{k}\,{\rm sn}\left(kx,\frac{1}{k}\right)
\eeq 
Analogous expressions can be derived for the other
elliptic functions and for the incomplete elliptic integral
of the second kind. In summary, the expressions useful
for conversion of $k$ into its reciprocal value read
\bea 
 \label{eq:amktilde}
{\rm am}(x,k) &=& \arcsin(\tilde{k}\,{\rm sn}(x/\tilde{k},\tilde{k}))
\\
{\rm sn}(x,k)&=&\tilde{k}\,{\rm sn}(x/\tilde{k},\tilde{k})
\\
{\rm cn}(x,k)&=&{\rm dn}(x/\tilde{k},\tilde{k})
\\ \label{eq:dnktilde}
{\rm dn}(x,k)&=&{\rm cn}(x/\tilde{k},\tilde{k})
\\ \label{eq:fjktilde}
{\rm F_J}(\varphi,k) &=& \tilde{k}\, {\rm F_J}(\tilde\varphi,\tilde{k})
\\ \label{eq:ejktilde}
{\rm E_J}(\varphi,k) &=& {\rm E_J}(\tilde\varphi,\tilde k)/\tilde k+(\tilde k -1/\tilde k)\,{\rm F_J}(\tilde\varphi,\tilde k)
\eea
in which $k=1/\tilde k$ and   $k\sin\varphi=\sin\tilde{\varphi}$.

%%%%%%%%%%%%%%%%%%%%%%%%%%%%%%%%%%%%%%%%%%%%%%
\subsection{Expressions for rotation and displacement functions}

Let us get back to the beam deformation problem.
Making use of the definition of the incomplete elliptic integral of the first kind, ${\rm F_J}$, equation (\ref{ee263}) can be written as
\beq \label{ee270}
\frac{2\,{\rm sgn}\,\kappa_a}{\sqrt{c_1+A}}
\left({\rm F_J}\left(\frac{\varphi(x)+\alpha}{2},k\right)-{\rm F_J}\left(\frac{\varphi_a+\alpha}{2},k\right)\right)= x
\eeq 
from which
%\beq 
%{\rm F_J}\left(\frac{\varphi+\alpha}{2},k\right) = {\rm F_J}\left(\frac{\varphi_a+\alpha}{2},k\right) + \,{\rm sgn}\,\kappa_a\frac{x\sqrt{c_1+A}}{2}
%\eeq 
%and finally
\beq \label{ee267}
\varphi(x) = 2\,{\rm am} \left( a +bx,k\right)-\alpha
\eeq 
where
\bea \label{ee267b}
a &=& {\rm F_J}\left(\frac{\varphi_a+\alpha}{2},k\right)
\\
\label{ee267c}
b &=& \frac{\sqrt{c_1+A}}{2}\,{\rm sgn}\,\kappa_a
\eea 
This is the analytical solution for the sectional
rotation $\varphi$ as function of the local coordinate $x$. 
Auxiliary constants $\alpha$, $k$, $c_1$ and $A$ depend
on the end forces and moment at the left end of the beam
and on the rotation of the left end section
(note that $\kappa_a=-M_{ab}/EI$). Of course, not all
of these variables are prescribed in advance and they
need to be determined from appropriate boundary conditions.

Based on expression (\ref{ee267}) for the sectional rotation, it is possible to evaluate
\bea
\sin\varphi(x) &=& 2\cos\alpha\,{\rm sn}(a+bx,k)\,{\rm cn}(a+bx,k)-\sin\alpha\,(1-2\,{\rm sn}^2(a+bx,k))
\\
\cos\varphi(x) &=& 2\sin\alpha\,{\rm sn}(a+bx,k)\,{\rm cn}(a+bx,k)+\cos\alpha\,(1-2\,{\rm sn}^2(a+bx,k))
\eea 
When this is substituted into the right-hand sides of  (\ref{eq35})--(\ref{eq36}) with $\lambda_s$ set to 1 (in accordance with the assumption of axial inextensibility), integration of the resulting equations leads to
analytical expressions for displacement functions,
\bea 
u_s(x) &=& C_u-x-\frac{2}{bk^2}\,{\rm dn}(a+bx,k)\sin\alpha+\left(\frac{2}{bk^2}{\rm E_J}({\rm am}(a+bx,k),k)+x-\frac{2x}{k^2}\right)\cos\alpha
\\
w_s(x) &=& C_w+\frac{2}{bk^2}\,{\rm dn}(a+bx,k)\cos\alpha+\left(\frac{2}{bk^2}{\rm E_J}({\rm am}(a+bx,k),k)+x-\frac{2x}{k^2}\right)\sin\alpha
\eea 
in which $C_u$ and $C_w$ are integration constants that need to be determined from the boundary conditions.

In cases when $k>1$, the derived formulae for the sectional rotation
and centerline displacements can be transformed into
expressions that use parameter $\tilde k=1/k$.
It is worth noting that the case of $k>1$ arises if $A>c_1$, which is equivalent
to $2F_{ab}(1+\cos(\varphi_a+\alpha))>M_{ab}^2/EI$.
Making use of identities (\ref{eq:amktilde}), (\ref{eq:dnktilde})  and (\ref{eq:ejktilde}), we obtain
\bea\label{eq:284} 
\varphi(x) &=& 2\,\arcsin(\tilde{k}\,{\rm sn}(\tilde{a}+\tilde{b}x,\tilde{k}))-\alpha 
\\
\label{eq288}
u_s(x) &=& C_u-x-\frac{2\tilde{k}}{\tilde{b}}\,{\rm cn}(\tilde{a}+\tilde{b}x,\tilde{k})\sin\alpha+\left(\frac{2}{\tilde{b}}\,{\rm E_J}({\rm am}(\tilde{a}+\tilde{b}x,\tilde{k}),\tilde{k})-x\right)\cos\alpha
\\
\label{eq287}
w_s(x) &=& C_w+\frac{2\tilde{k}}{\tilde{b}}\,{\rm cn}(\tilde{a}+\tilde{b}x,\tilde{k})\cos\alpha+\left(\frac{2}{\tilde{b}}\,{\rm E_J}({\rm am}(\tilde{a}+\tilde{b}x,\tilde{k}),\tilde{k})-x\right)\sin\alpha
\eea 
where
\bea \label{eq:tilk}
\tilde{k}&=&\frac{1}{k}=\sqrt{\frac{c_1+A}{2A}} =
\sqrt{\frac{2(F_{ab}+N_{ab})+EI\kappa_a^2}{4F_{ab}}}
=\sqrt{\sin^2\frac{\varphi_a+\alpha}{2}+\frac{M_{ab}^2}{4EIF_{ab}}}
\\ \label{eq:tila}
\tilde a &=& ka =\sqrt{\frac{2A}{c_1+A}} {\rm F_J}\left(\frac{\varphi_a+\alpha}{2},k\right)
={\rm F_J}\left(\arcsin\left(\frac{1}{\tilde k}\sin\frac{\varphi_a+\alpha}{2}
\right),\tilde k\right) \\
\label{eq:tilb}
\tilde b &=& kb =\sqrt{\frac{2A}{c_1+A}} \frac{{\rm sgn}\,\kappa_a}{2}\sqrt{c_1+A} =
\sqrt{\frac{A}{2}} \,{\rm sgn}\,\kappa_a = \sqrt{\frac{F_{ab}}{EI}}\,{\rm sgn}\,\kappa_a
\eea 

It is worth noting that the solution described by (\ref{ee267}) or (\ref{eq:284}) is valid only 
as long as the sign of the curvature
does not change. 
If the sign changes inside the interval of interest $(0,L)$ that represents the analyzed beam of length $L$,
formula (\ref{ee267}) can be used only up to the
inflexion point of the deformed centerline.
A systematic treatment of deformed shapes with inflexion points is presented in Appendix~\ref{app:inflex} and leads to formula (\ref{ee270w}).

\subsection{Special cases -- straight beam and uniformly curved beam}
\label{app:example}

The special case with zero denominator in 
(\ref{eee264}) needs to be treated separately. 
This happens only if $M_{ab}=0$
and $N_{ab}=-F_{ab}$, the latter condition leading to
$\sin(\varphi_a+\alpha)=0$. Inspection of the original
equation (\ref{e221}) with initial conditions $\varphi(0)=-\alpha$
and $\varphi'(0)=0$
shows that, in this particular case,
the solution is $\varphi(x)=-\alpha$, i.e., it is 
constant. The first derivative, $\varphi'(x)$,
identically vanishes, which explains why the 
general solution procedure developed in Section~\ref{app:anal1}
is not applicable to this particular case
(recall that the procedure started by multiplying both sides of the
original equations by function $\varphi'(x)$).

Another special case arises when $M_{ab}\ne 0$ and $X_{ab}=Z_{ab}=0$, because then
$c_1=\kappa_a^2$, $A=0$, $\alpha$ is undetermined and $k=0$. For $k=0$, 
functions ${\rm F_J}$ and 
``am'' reduce to identities, i.e., 
${\rm F_J}(\varphi,0)=\varphi$ and 
${\rm am}(x,0)=x$.
Formulae (\ref{ee267b})--(\ref{ee267c})  yield $a=\varphi_a/2$ and $b=\kappa_a/2$, and the rotation function
is according to (\ref{ee267}) given by
\beq
\varphi(x) = 2a+2bx = \varphi_a+ \kappa_a x
\eeq 
This is the case of a beam with constant curvature,
loaded at its end sections only by two moments of the same
magnitude but opposite orientations.

\subsection{Example -- cantilever loaded by an inclined force}\label{sec:examplecantilever}

\begin{figure}[h]
\centering
\begin{tabular}{cc}
(a) & (b)\\
\\
\includegraphics[width=0.3 \linewidth]{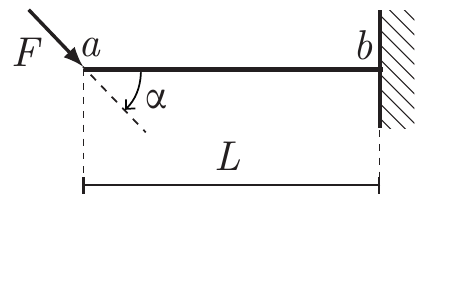} & \includegraphics[width=0.3 \linewidth]{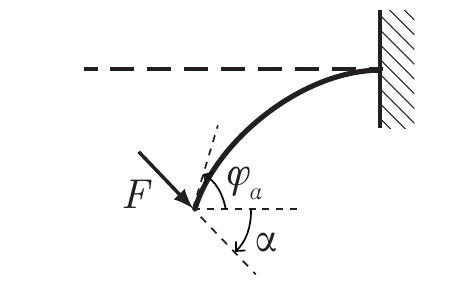}
\end{tabular}
\caption{Cantilever fixed at its right end and loaded at its left end by an inclined force: (a) Geometry in the undeformed state, (b) deformed state with $\alpha$ denoting the inclination angle of the applied force (positive clockwise) and $\varphi_{a}$ denoting the left-end rotation (positive anti-clockwise)} 
\label{f:cantilever}
\end{figure}

As a more challenging example, consider 
 a cantilever fixed at its
right end and loaded at its left end by a force of magnitude $F$ along an inclined line.
The horizontal and vertical components of the applied force
correspond to the left-end forces $X_{ab}$ and $Z_{ab}$,
and formula (\ref{ee5}) indicates that parameter $\alpha$ used by
the analytical solution is equal to the angle by which the
oriented direction of the force  deviates from the
$x$-axis, measured clockwise. Therefore, this angle can
be directly denoted as $\alpha$ and considered as a given quantity, as shown in Fig.~\ref{f:cantilever}.
The left-end moment, $M_{ab}$, is set equal to zero.

Without loss of generality, we can assume that $\alpha\in[0,\pi)$, and it is reasonable to expect that the
curvature $\kappa(x)$ has a negative sign for $x\in(0,L]$
and that $\alpha+\varphi_a\in[0,\pi]$.
Auxiliary constants
needed to construct the analytical solution 
are now evaluated as follows:
\bea 
%X_{ab} &=& F\cos\alpha \\
%Z_{ab} &=& F\sin\alpha \\
%\kappa_a &=& \frac{M_{ab}}{EI}=0 \\
%c_1 &=& -\frac{2F}{EI}\cos(\alpha+\varphi_a)\\
%c_2 &=& \frac{2F}{EI}\cos\alpha \\
%c_3 &=&  -\frac{2F}{EI}\sin\alpha\\
%A &=& \frac{2F}{EI}\\
\label{eq:ex284}
\tilde{k} &=&   \sin\frac{\alpha+\varphi_a}{2}
\\
\tilde a &=& {\rm F_J}\left(\arcsin 1,\tilde k\right) = {\rm F_J}(\pi/2,\tilde k) = {\rm K}(\tilde k)
\\ \label{eq:ex286}
\tilde b &=& -\sqrt{\frac{F}{EI}}
\eea

Note that $\alpha$ is a given angle while the value of the left-end rotation, $\varphi_a$, is related to the applied force $F$ and increases from 0 to positive values that never exceed
$\pi-\alpha$
(unless $\alpha=\pi$, in which case the beam is under axial tension and $\varphi_a$ remains equal to 0).
Parameter $\tilde k$ given by (\ref{eq:ex284}) never exceeds 1, and so
it is preferable to use elliptic functions and elliptic integrals with parameter $\tilde k$ and not $k$.
In the evaluation of parameter $\tilde b$ according
to formula (\ref{eq:tilb}), we have replaced ${\rm sgn}\;\kappa_a$,
which would be zero, by the signum of the curvature
in the neighborhood of the left end, which is expected
to be $-1$ if $\alpha\in(0,\pi)$. This explains the
negative sign in (\ref{eq:ex286}).

The analytical solution for the rotation function is given in (\ref{eq:284}), and substitution into boundary condition $\varphi(L)=0$ 
leads to
\beq 
2\,\arcsin(\tilde{k}\,{\rm sn}(\tilde{a}+\tilde{b}L,\tilde{k}))-\alpha=0 
\eeq 
which is satisfied if 
\beq\label{eq308} 
\tilde{a}+\tilde{b}L={\rm F_J}\left(\arcsin\left(\frac{1}{\tilde k}\sin\frac{\alpha}{2}\right),\tilde k\right)
\eeq 
Based on expressions (\ref{eq:ex284})--(\ref{eq:ex286}) for $\tilde k$, $\tilde a$ and $\tilde b$,
it is possible to rewrite (\ref{eq308}) as  an equation linking the 
applied force, $F$, to the left-end rotation, $\varphi_a$:
\beq\label{eq309} 
L\sqrt{\frac{F}{EI}}=K\left(\sin\frac{\alpha+\varphi_a}{2}\right)-{\rm F_J}\left(\arcsin\sqrt{\frac{1-\cos \alpha}{1-\cos(\alpha+\varphi_a)}},\sin\frac{\alpha+\varphi_a}{2}\right)
\eeq 
Instead of solving this  nonlinear equation numerically for each given value of $F$,
one can use a parametric description with
$\varphi_a$ considered as a control parameter
that varies in a suitable range. 
It is convenient to introduce auxiliary functions 
\bea 
%C(\varphi) &=& {\rm F_J}\left(\arcsin\sqrt{\frac{1-\cos \alpha}{1-\cos(\alpha+\varphi)}},\sin\frac{\alpha+\varphi}{2}\right)
%\\
\label{eq:B}
B(\varphi) &=& K\left(\sin\frac{\alpha+\varphi}{2}\right)-{\rm F_J}\left(\arcsin\sqrt{\frac{1-\cos \alpha}{1-\cos(\alpha+\varphi)}},\sin\frac{\alpha+\varphi}{2}\right)
\\
\label{eq:D}
D(\varphi) &=& {\rm E_J}\left(\arcsin\sqrt{\frac{1-\cos \alpha}{1-\cos(\alpha+\varphi)}},\sin\frac{\alpha+\varphi}{2}\right)
\eea 
Function $B(\varphi_a)$ represents the right-hand side of  (\ref{eq309}), and $D(\varphi_a)$
will prove to be useful in the formulae for displacements.

For a given series of values of $\varphi_a$, the corresponding forces 
\beq \label{eq:F}
F=\frac{EI}{L^2}\,B^2(\varphi_a)
\eeq 
can be evaluated from (\ref{eq309}).
The displacement functions are then given by 
(\ref{eq288})--(\ref{eq287}), 
in which integration constants $C_w$ and $C_u$ are determined from boundary conditions $u_s(L)=0$ and $w_s(L)=0$. 
The resulting displacements at the left end
of the cantilever turn out to be
\bea\label{eq:ua}
u_a &=& u_s(0) = L(1+\cos\alpha)-\frac{L\sin\alpha}{B(\varphi_a)}\sqrt{2\cos\alpha-2\cos(\alpha+\varphi_a)}-\frac{2L\cos\alpha}{B(\varphi_a)}\,\left({\rm E}\left(\sin\frac{\alpha+\varphi_a}{2}\right)-D(\varphi_a)\right)
\\
\label{eq:wa}
w_a &=& w_s(0) = L\sin\alpha+\frac{L\cos\alpha}{B(\varphi_a)}\sqrt{2\cos\alpha-2\cos(\alpha+\varphi_a)}-\frac{2L\sin\alpha}{B(\varphi_a)}\,\left({\rm E}\left(\sin\frac{\alpha+\varphi_a}{2}\right)-D(\varphi_a)\right)
\eea 
The end displacement projected onto the direction 
of applied force is then easily evaluated as
\beq \label{eq:a68}
u_F = u_a\cos\alpha+w_a\sin\alpha =
L(1+\cos\alpha)- \frac{2L}{B(\varphi_a)}\,\left({\rm E}\left(\sin\frac{\alpha+\varphi_a}{2}\right)-D(\varphi_a)\right)
\eeq 

As a special case, consider an {\bf axially compressed cantilever},
characterized by $\alpha=0$. For $\varphi>0$, the auxiliary functions defined in (\ref{eq:B})--(\ref{eq:D}) simplify to
\bea 
B(\varphi) &=& {\rm K}\left(\sin\frac{\varphi}{2}\right)
\\
D(\varphi) &=& 0
\eea 
and thus formulae (\ref{eq:F})--(\ref{eq:wa}) yield
\bea \label{eq:ex298}
F&=&\frac{EI}{L^2}K^2\left(\sin\frac{\varphi_a}{2}\right)
\\
\label{eq:ex299}
u_a &=&  2L\left(1-\frac{E\left(\sin\displaystyle\frac{\varphi_a}{2}\right)}{K\left(\sin\displaystyle\frac{\varphi_a}{2}\right)}\right)
\\ \label{eq:ex300}
w_a &=&  \frac{2L\sin\displaystyle\frac{\varphi_a}{2}}{K\left(\sin\displaystyle\frac{\varphi_a}{2}\right)}
\eea 
We have excluded the case of $\varphi_a=0$, which leads
to undetermined fractions in the definitions of $B$ and $D$.
However, the derived formulae (\ref{eq:ex298})--(\ref{eq:ex300}) have no singularity at $\varphi_a=0$. Since $E(0)=\pi/2$ and $K(0)=\pi/2$, we obtain from (\ref{eq:ex298}) the Euler critical force $F=EI\pi^2/(4L^2)$
and from (\ref{eq:ex299})--(\ref{eq:ex300}) zero displacements $u_a$ and $w_a$ at the onset of buckling of a perfectly straight cantilever.

The analytical solution derived for $\alpha=0$ is used as a reference in Section~\ref{sec:buckling}.
Other special cases, $\alpha=\pi/3$ and $\alpha=2\pi/3$, are exploited in the analysis of a periodic hexagonal cell in Section~\ref{sec:honeycomb}.

\subsection{Analytical solution with an inflexion point}
\label{app:inflex}

The solution described by (\ref{ee267}) or (\ref{eq:284}) is valid only 
as long as the sign of the curvature
does not change. 
If the sign changes inside the interval of interest $[0,L]$ representing a beam of length $L$,
these formulae  can be used only up to the
inflexion point of the deformed centerline. 
Let us denote
the initial coordinate of the inflexion point by $x_{in}$
and the corresponding value of rotation by $\varphi_{in}=\varphi(x_{in})$.
The curvature
sign changes from ${\rm sgn}\,\kappa_a$
to ${\rm sgn}\,\kappa_b=-{\rm sgn}\,\kappa_a$ at point
$x=x_{in}$ characterized by the condition 
$\varphi'(x_{in})=0$, which is the case if $2G(\varphi_{in})=C$
where $C$ is the integration
constant defined in (\ref{ee254}) and $G$ is the function 
defined in (\ref{ee251x}).

Suppose that $x_{in}$ and $\varphi_{in}$ are known.
Equations (\ref{ee263}) and (\ref{ee270})--(\ref{ee267}) are valid for $x\in[0,x_{in}]$.
By substituting $x=x_{in}$ and $\varphi=\varphi_{in}$
into (\ref{ee270}), we get
the identity
\beq \label{ee270a}
\frac{2\,{\rm sgn}\,\kappa_a}{\sqrt{c_1+A}}
\left({\rm F_J}\left(\frac{\varphi_{in}+\alpha}{2},k\right)-{\rm F_J}\left(\frac{\varphi_a+\alpha}{2},k\right)\right)= x_{in}
\eeq 
For $x\in[x_{in},L]$,
 equation (\ref{ee270}) is replaced by
\beq \label{ee270x}
-\frac{2\,{\rm sgn}\,\kappa_a}{\sqrt{c_1+A}}
\left({\rm F_J}\left(\frac{\varphi(x)+\alpha}{2},k\right)-{\rm F_J}\left(\frac{\varphi_{in}+\alpha}{2},k\right)\right)= x- x_{in}
\eeq 
The key point here is that the values of auxiliary constants
$c_1$, $A$, $k$ and $\alpha$ are the same as for $x\in[0,x_{in}]$.
Therefore, it is possible to eliminate $x_{in}$ by taking
the sum of (\ref{ee270a}) and (\ref{ee270x}), which leads to
\beq \label{ee270y}
\frac{2\,{\rm sgn}\,\kappa_a}{\sqrt{c_1+A}}
\left(2{\rm F_J}\left(\frac{\varphi_{in}+\alpha}{2},k\right)-{\rm F_J}\left(\frac{\varphi(x)+\alpha}{2},k\right)-{\rm F_J}\left(\frac{\varphi_a+\alpha}{2},k\right)\right)= x
\eeq 
The formal analytical solution valid for  $x\in[x_{in},L]$
thus reads
\beq \label{ee270z}
\varphi(x) =
2{\rm am}\left(2{\rm F_J}\left(\frac{\varphi_{in}+\alpha}{2},k\right)-{\rm F_J}\left(\frac{\varphi_a+\alpha}{2},k\right)-({\rm sgn}\,\kappa_a)\frac{\sqrt{c_1+A}}{2}x,k\right)-\alpha
\eeq 
This expression still contains the rotation at the inflexion
point, $\varphi_{in}$. We will now show how 
$\varphi_{in}$ could be determined, but at the same time 
it will turn out that its value is actually not needed,
because the integral that corresponds to ${\rm F_J}\left((\varphi_{in}+\alpha)/2,k\right)$
can be converted into a quantity that depends only on the 
given parameters.

As already mentioned, $\varphi_{in}$ satisfies
condition $2G(\varphi_{in})=C$, which can be rewritten as
\beq 
\frac{X_{ab}}{EI}(1-\cos\varphi_{in})+\frac{Z_{ab}}{EI}\sin\varphi_{in} = \frac{X_{ab}}{EI}(1-\cos\varphi_{a})+\frac{Z_{ab}}{EI}\sin\varphi_{a} + \half\kappa_a^2
\eeq 
and further transformed into
\bea 
%X_{ab}\cos\varphi_{in} - Z_{ab}\sin\varphi_{in} &=& X_{ab}\cos\varphi_{a} - Z_{ab}\sin\varphi_{a} - \half EI\kappa_a^2
%\\
%F_{ab}\cos\varphi_{in}\cos\alpha - F_{ab}\sin\varphi_{in}\sin\alpha &=& F_{ab}\cos\varphi_{a}\cos\alpha - F_{ab}\sin\varphi_{a}\sin\alpha - \half EI\kappa_a^2
%\\
\cos(\varphi_{in}+\alpha) &=& \cos(\varphi_{a}+\alpha)- \frac{EI}{2F_{ab}}\kappa_a^2
\eea 
Therefore, the value of $\varphi_{in}$, representing the rotation
at the inflexion point (i.e., a local extreme of the rotation),
can be expressed as
\beq 
\varphi_{in} = \pm \arccos \left(\cos(\varphi_{a}+\alpha)- \frac{EI}{2F_{ab}}\kappa_a^2 \right) - \alpha + 2n\pi
\eeq 
where the sign before $\arccos$ and the integer $n$ are selected
depending on the value of $\varphi_a$ and the sign of $\kappa_a$ such that ${\rm sgn}(\varphi_{in}-\varphi_a)={\rm sgn}\,\kappa_a$
and $|\varphi_{in}-\varphi_a|$ is minimized among all 
roots satisfying this constraint.

In fact, what matters more than the precise value of $\varphi_{in}$ is that if $\tilde\varphi$ is set to 
$(\varphi_{in}+\alpha)/2$, the denominator of the integral
in (\ref{ee263}) vanishes. In other words, $\tilde\varphi_{in}=(\varphi_{in}+\alpha)/2$
satisfies condition
\beq 
1-k^2\sin^2\tilde\varphi_{in} = 0
\eeq 
from which  
\beq 
 \tilde\varphi_{in} = \pm\arcsin \frac{1}{k}
\eeq 
Consequently, when we evaluate the analytical solution
(\ref{ee270z}),
the term that depends on $\varphi_{in}$ can be according to
() expressed as
\bea\nonumber
{\rm F_J}\left(\frac{\varphi_{in}+\alpha}{2},k\right)&=&{\rm F_J}\left(\tilde\varphi_{in},k\right)=
\frac{1}{k}{\rm F_J}\left(\arcsin\left(k\sin\tilde\varphi_{in}\right),\frac{1}{k}\right) = \frac{1}{k}{\rm F_J}\left(\pm \frac{\pi}{2},\frac{1}{k}\right)
= \pm\frac{1}{k}{\rm K}\left(\frac{1}{k}\right)
\label{ee285w}
\eea
The sign to be selected in the last expression in (\ref{ee285w})
corresponds to the sign of $\tilde\varphi_{in}$, which is the
same as the sign of $\kappa_a$.
Making use of (\ref{ee285w}) with the proper sign, formula (\ref{ee270z}) can be rewritten as
\beq \label{ee270w}
\varphi(x) =
2\,{\rm am}\left(\frac{2\,{\rm sgn}\,\kappa_a}{k}{\rm K}\left(\frac{1}{k}\right)-{\rm F_J}\left(\frac{\varphi_a+\alpha}{2},k\right)-({\rm sgn}\,\kappa_a)\frac{\sqrt{c_1+A}}{2}x,k\right)-\alpha, \hskip 5mm x_{in}\le x \le L
\eeq 
in which
\beq 
x_{in} = \frac{2\,{\rm sgn}\,\kappa_a}{\sqrt{c_1+A}}
\left(\frac{1}{k}{\rm K}\left(\frac{1}{k}\right)-{\rm F_J}\left(\frac{\varphi_a+\alpha}{2},k\right)\right)
\eeq

%%%%%%%%%%%%%%%%%%%%%%%%%%%%%%%%%%%%%%%%%%%%%%%%%%%%%%%%%%%%%%%

\section{Critical load for an axially compressible cantilever}
\label{appC}

Stability of an elastic equilibrium state can be evaluated based on the second variation of the potential energy functional. Formally, the second variation
of a functional $\Pi(\boldsymbol{u})$
can be defined as the second Gateaux derivative, i.e.,
as 
\beq \label{eq:delta2}
\delta^2\Pi(\boldsymbol{u},\delta\boldsymbol{u}) = \frac{{\rm d}^2\Pi(\boldsymbol{u}+h\,\delta\boldsymbol{u})}{{\rm d}h^2} \Big\vert_{h=0}
\eeq 
When this definition is applied to the potential energy of a beam $E_p$
given by (\ref{eq:Ep}), considered as a functional dependent
on centerline displacement functions $u_s$ and $w_s$,
a careful processing of formula (\ref{eq:delta2}) leads to
a relatively lengthy expression. However, we are primarily interested in stability of the solution that
corresponds to a beam that still remains straight
but is uniformly compressed. 

To be specific, consider a cantilever of length $L$, fixed at its right end and loaded at its left end by a compressive
force $P$.
The state of uniform compression is characterized by
displacement functions $u_s(x)=(\lambda_s-1)(L-x)$ and $w_s(x)=0$ where 
\beq\label{eq216} 
\lambda_s = 1-\frac{P}{EA}
\eeq
is a given positive constant that
represents the (uniform) stretch of the beam axis. 
The load $P>0$ is considered as compressive,
and so $0<\lambda_s<1$.
The second variation
of potential energy evaluated for such a particular state
turns out to be
\beq\label{eq215} 
\delta^2E_p(\delta u_s,\delta w_s) =
\int_0^L EA\left(\delta u_s'^2+\frac{\lambda_s-1}{\lambda_s}\,\delta w_s'^2\right)\dx + \int_0^L\frac{EI}{\lambda_s^2}\,\delta w_s''^2\,\dx
\eeq 
 The dependence on the state at
which the second variation is taken (i.e., on $u_s$
and $w_s$) is not marked explicitly on the left-hand
side of (\ref{eq215}), because the presented expression for the
second variation is not valid for general $u_s$ and
$w_s$ but only for the special case of a uniformly
compressed beam, which is fully described by the
scalar parameter $\lambda_s$.

If there exist admissible variations $\delta u_s$ and $\delta w_s$ for which $\delta^2E_p<0$, the considered equilibrium state (straight uniformly compressed beam)
is unstable. To find the critical load associated with
the onset of instability, we look for the 
minimum value of $P$ (and thus maximum value of $\lambda_s$) for which $\delta^2E_p\le 0$
for some nonzero combination of admissible
variations $\delta u_s$ and $\delta w_s$.
Since the contribution of $\delta u_s$
to the right-hand
side of (\ref{eq215}) is always non-negative,
the most ``dangerous'' case occurs when
$\delta u_s=0$. Also, since (\ref{eq215}) contains
only the first and second derivatives of function $\delta w$
but not the function itself, we can introduce function
$\delta\varphi=-\delta w_s'$ and then
search for 
nonzero $\delta\varphi$ that satisfies boundary
condition $\delta \varphi(L)=0$ (clamped right end) and the inequality
\beq\label{eq217} 
\int_0^L EA\frac{\lambda_s-1}{\lambda_s}\,\delta \varphi^2\,\dx + \int_0^L\frac{EI}{\lambda_s^2}\,\delta \varphi'^2\,\dx \le 0
\eeq 
It is clear that if $\lambda_s\ge 1$, 
the left-hand side of (\ref{eq217}) is positive for any nonzero $\delta\varphi$.
Therefore, stability cannot be lost in tension (for the
present model). The question is what happens in compression, when the factor $(\lambda_s-1)/\lambda_s$
multiplying $\delta \varphi^2$ becomes negative. 
Since we restrict attention to $\lambda_s>0$,
condition (\ref{eq217}) can be rewritten as
\beq\label{eq218} 
\frac{EA}{EI}\lambda_s(1-\lambda_s)   \ge  \frac{\int_0^L\delta \varphi'^2\,\dx}{\int_0^L \,\delta\varphi^2\,\dx}
\eeq 
and finally, based on (\ref{eq216}), it can be converted
into 
\beq\label{eq219} 
\frac{P}{EI}\left(1-\frac{P}{EA}\right)   \ge  \frac{\int_0^L\delta \varphi'^2\,\dx}{\int_0^L \,\delta \varphi^2\,\dx}
\eeq 

To find the critical value of $P$,
we need to minimize the right-hand side
of (\ref{eq219}) over the set of all nonzero functions $\delta \varphi$
that satisfy the boundary condition $\delta\varphi(L)=0$. 
Minimization of the fraction on the right-hand side
of (\ref{eq219}) can be replaced by minimization of the numerator subject to the constraint that the denominator 
be equal to 1. Introducing  a Lagrange multiplier
$\Lambda$ to enforce this constraint, we end up with the
differential eigenvalue problem
\beq 
\delta\varphi'' + \Lambda \,\delta\varphi = 0
\eeq 
supplemented by boundary conditions $\delta\varphi'(0)=0$ and $\delta\varphi(L)=0$. The smallest eigenvalue $\Lambda=\pi^2/(4L^2)$ then represents the minimum value
of the fraction on the right-hand side
of (\ref{eq219}), attained by setting $\delta\varphi(x)=\cos\frac{\pi x}{2L}$. 
Consequently, stability of the solution that corresponds to a uniformly compressed beam is lost if the applied load satisfies
condition
\beq\label{eq219y} 
\frac{P}{EI}\left(1-\frac{P}{EA}\right)   \ge  \frac{\pi^2}{4L^2}
\eeq 

In the limit of $EA\to\infty$, 
the left-hand side of (\ref{eq219y}) reduces to $P/EI$
and the smallest load for which the condition holds is the Euler critical load
\beq\label{eq220} 
P_E = \frac{EI\pi^2}{4L^2}
\eeq 
For a finite value of $EA$, the critical load 
\beq\label{eq222} 
P_{cr} = \frac{EA}{2}\left(1-\sqrt{1-\frac{EI\pi^2}{EAL^2}}\right)= EA\left(\frac{1}{2}-\sqrt{\frac{1}{4}-\frac{P_E}{EA}}\right) 
\eeq
is found
as the smaller root of the quadratic equation
\beq\label{eq221} 
\frac{P}{EI}\left(1-\frac{P}{EA}\right)   -  \frac{\pi^2}{4L^2} = 0
\eeq 
Typically, $EA\gg P_E$, and the ``exact'' expression
from formula (\ref{eq222}) can be approximated as follows:
\beq\label{eq223} 
P_{cr} =  \frac{EA}{2}\left(1-\sqrt{1-\frac{4P_E}{EA}}\right) \approx \frac{EA}{2}\left(1-\left(1-\frac{2P_E}{EA}-\frac{2P_E^2}{(EA)^2}\right)\right) = P_E\left(1+\frac{P_E}{EA}\right)
\eeq 
This confirms that $P_{cr}\to P_E$ as $EA\to\infty$,
which is not so obvious from (\ref{eq222}).

%For a cantilever characterized by kinematic boundary conditions $\delta w_s(L)=0$ and $\delta w_s'(L)=0$, the minimum value of the fraction on the right-hand side of (\ref{eq219}) is $\pi^2/(4L^2)$, and so formulae (\ref{eq222}) and (\ref{eq223}) remain valid provided that the actual length $L$ is replaced by the buckling length $L_b=2L$ and the Euler critical load is interpreted as $P_E=EI\pi^2/L_b^2$. 
It is worth noting that
the largest possible critical load $P_{cr}=EA/2$ is obtained from formula (\ref{eq222}) for $EA=4P_E$ and stability would never be lost if $EA<4P_E$, i.e., if
$EAL^2/EI<\pi^2$. However, this is already far from the range in which the adopted assumptions are physically meaningful. The objective here is to describe slender beams, which buckle at strains that can still be considered as small. This is true only if $P_{cr}\ll EA$
or, equivalently, $EI\pi^2\ll 4EAL^2$. 
The calculations in the paper have been done for
$EAL^2/EI=100$ and 10,000, which is indeed much larger than $\pi^2$.

%\noindent\textbf{Unnumbered figure}

%\begin{center}
%\includegraphics[width=7pc,height=8pc,draft]{empty}
%\end{center}

%== Figure 4 ==
%% Example for figure inside appendix
%\begin{figure}[t]
%\centerline{\includegraphics[height=10pc,width=78mm,draft]{empty}}
%\caption{This is an example for appendix figure.\label{fig5}}
%\end{figure}

%\begin{center}
%\begin{table}[b]%
%\centering
%\caption{This is an example of Appendix table showing food requirements of army, navy and airforce.\label{tab4}}%
%\begin{tabular*}{300pt}{@{\extracolsep\fill}lcc@{\extracolsep\fill}}%
%\toprule
%\textbf{col1 head} & \textbf{col2 head} & \textbf{col3 head} \\
%\midrule
%col1 text & col2 text & col3 text \\
%col1 text & col2 text & col3 text \\
%col1 text & col2 text & col3 text\\
%\bottomrule
%\end{tabular*}
%\end{table}
%\end{center}

%Example for an equation inside appendix
%\begin{equation}
%\mathcal{L}\quad \mathbf{\mathcal{L}} = i \bar{\psi} \gamma^\mu D_\mu \psi - \frac{1}{4} F_{\mu\nu}^a F^{a\mu\nu} - m \bar{\psi} \psi\label{eq25}
%\end{equation}

%\begin{center}
%\begin{tabular*}{300pt}{@{\extracolsep\fill}lcc@{\extracolsep\fill}}%
%\toprule
%\textbf{col1 head} & \textbf{col2 head} & \textbf{col3 head} \\
%\midrule
%col1 text & col2 text & col3 text \\
%col1 text & col2 text & col3 text \\
%col1 text & col2 text & col3 text \\
%\bottomrule
%\end{tabular*}
%\end{center}

%\nocite{*}% Show all bib entries - both cited and uncited; comment this line to view only cited bib entries;

\newpage
\bibliography{wileyNJD-AMA}%

%\clearpage

%\section*{Author Biography}

\end{document}